\documentclass[%
3p,%
preprint,
fleqn, 
times
]{elsarticle}

\makeindex

\usepackage{makeidx}         
\usepackage{graphicx}        
\usepackage{amsmath,amssymb}
\usepackage[english=american]{csquotes} 
\usepackage{dsfont} 
  \newcommand{\IP}{\ensuremath\mathds{P}}
  \newcommand{\IR}{\ensuremath\mathds{R}}
  \newcommand{\IN}{\ensuremath\mathds{N}}
\usepackage{tikz}
\usepackage{tabularx}             
  \newcolumntype{R}{>{\raggedleft\arraybackslash}X}  \newcolumntype{L}{>{\raggedright\arraybackslash}X}  \newcolumntype{C}{>{\centering\arraybackslash}X}
\usepackage{multirow} 
\usepackage{rotating} 
\usepackage{booktabs} 
\usepackage{bm}    
\usepackage{mathtools}
\usepackage{fixmath}
\usepackage{hyperref}
  \hypersetup{breaklinks=true,colorlinks=true,pdfborder={0 0 0}}
  \hypersetup{linkcolor=black, anchorcolor=black, citecolor=black, filecolor=black, urlcolor=black}
\usepackage{subcaption}
\newcommand*{\setT}{\ensuremath{\mathcal{T}}}
\newcommand*{\setE}{\ensuremath{\mathcal{E}}}
\newcommand*{\setV}{\ensuremath{\mathcal{V}}}
\newcommand*{\setA}{\ensuremath{\mathcal{A}}}
\newcommand*{\vphi}{\varphi}                                     
\newcommand*{\mapEE}{\hat{\vec{\vartheta}}}                      
\renewcommand*{\vec}[1]{{\boldsymbol{#1}}}                       
\DeclareMathAlphabet{\mathbfsf}{\encodingdefault}{\sfdefault}{bx}{n}
\newcommand*{\vecc}[1]{\mathbfsf{#1}}                            
\newcommand*{\grad}{\vec{\nabla}}                                
\renewcommand*{\div}{\vec{\nabla}\cdot}                          
\newcommand*{\dd}{\mathrm{d}}                                    
\newcommand*{\abs}[1]{\ensuremath{|#1|}}                         
\newcommand*{\card}[1]{\ensuremath{\##1}}                        
\newcommand*{\transpose}[1]{{#1}^\mathrm{T}}                     
\newcommand*{\invtrans}[1]{{#1}^\mathrm{-T}}                     
\DeclareMathOperator*{\diag}{diag}                               
\DeclareMathOperator*{\diam}{diam}                               

\newcommand{\II}{I\!I} \newcommand{\III}{I\!I\!I} \newcommand{\IV}{I\!V}   

\newcommand{\Matlab}{\mbox{MATLAB}}
\newcommand{\Octave}{GNU~\mbox{Octave}}
\newcommand{\MatOct}{\Matlab\,/\,\allowbreak\Octave} 

\usepackage{textcomp} 
\usepackage{listings}
  \definecolor{keywordcolor}{rgb}{0, 0.25, 0.5}
  \definecolor{commentcolor}{rgb}{0.2, 0.5, 0.2}
  \definecolor{stringcolor}{rgb}{0.5, 0.5, 0.2}
\newcommand*{\code}[1]{\mbox{\lstinline[basicstyle=\ttfamily\small]{#1}}}

\usepackage{caption}
\begin{document}

\title{FESTUNG: A~\MatOct~toolbox for the discontinuous Galerkin method. Part II: Advection operator and slope limiting}
\author[FAU]{Balthasar Reuter}
  \ead{reuter@math.fau.de}
\author[FAU]{Vadym Aizinger\corref{cor}}
  \ead{aizinger@math.fau.de}
\author[FAU]{Manuel Wieland}
  \ead{manuel.wieland@studium.uni-erlangen.de}
\author[Rice]{Florian Frank}
  \ead{florian.frank@rice.edu}
\author[FAU]{Peter Knabner}
  \ead{knabner@math.fau.de}
\address[FAU]{Friedrich--Alexander University of Erlangen--N\"urnberg, Department of Mathematics, 
Cauerstra{\ss}e~11, 91058~Erlangen, Germany}
\address[Rice]{Rice University, Department of Computational and Applied Mathematics, 
6100 Main Street -- MS~134, Houston, TX 77005-1892, USA}
\cortext[cor]{Corresponding author}

\begin{abstract}
This is the second in a~series of papers on implementing a~discontinuous Galerkin (DG) method 
as an~open source \MatOct~toolbox.
The intention of this ongoing project is to offer a rapid prototyping package for application
development using DG~methods. 
The implementation relies on fully vectorized matrix\,/\,vector operations and is comprehensively documented. 
Particular attention was paid to maintaining a~direct mapping between discretization terms and code routines as well as to supporting the full code functionality in \Octave.
The present work focuses on a~two-dimensional time-dependent linear advection equation with space\,/\,time-varying coefficients, and provides a general order implementation of several slope limiting schemes for the DG~method.
\end{abstract}


\maketitle

 
\section{Introduction}

The development milestones for the \MatOct~toolbox \textsl{FESTUNG} (\textsl{F}inite 
\textsl{E}lement \textsl{S}imulation \textsl{T}oolbox for \textsl{UN}structured \textsl{G}rids) 
available at~\cite{FESTUNG,FESTUNGGithub} run somewhat counter to the history of the development 
of the discontinuous Galerkin (DG) methods. 
Thus, our first paper in series~\cite{FESTUNG1} introduced a local discontinuous Galerkin 
discretization for a time-dependent diffusion equation using the numerical methods introduced
in~\cite{CockburnShu1998}. The current work, however, enhances the package with the functionality 
for purely hyperbolic equations---namely the original purpose of the DG~method proposed by 
Reed and Hill in~\cite{ReedHill1973} and analyzed by Johnson and Pitk\"aranta in~\cite{Johnson1986}.
The reason behind this time inversion is that the numerical and software development technology 
necessary to produce a fully functional DG~solver for hyperbolic equations has to include upwind 
fluxes and slope limiters---both tasks more complicated to solve in a~computationally efficient 
manner than those needed for a pure diffusion equation.
\par
The continued development of this toolbox still adheres to the same design principles declared in~\cite{FESTUNG1}:
\begin{enumerate}
\item Design a~general-purpose software package using the DG~method for a~range of standard applications 
and provide this toolbox as a~research and learning tool in the open source format~(cf.~\cite{FESTUNG}).
\item Supply a~well-documented, intuitive user-interface to 
ease adoption by a~wider community of application and engineering professionals.
\item Relying on the vectorization capabilities of \MatOct, optimize the computational performance 
of the toolbox components and demonstrate these software development strategies.
\item Maintain throughout full compatibility with \Octave~to support users of open source software.
\end{enumerate}
\par
We refer to~\cite{FESTUNG1} for a literature review on DG~methods and open source 
packages offering a~DG~capability. The present work expands the functionality of the 
numerical solver published in the first paper in series by adding linear advection terms and 
vertex-based slope limiters of general order. 
The latter development is particularly interesting, since,
to the best of our knowledge, no closed form description of vertex-based slope limiters for general
order discretizations are to be found in the literature, even less so implementations of such limiters. In addition to hierarchical vertex-based limiters of Kuzmin~\cite{Kuzmin2012}, this publication and the accompanying code includes an extension of the standard linear vertex-based slope limiter to general order discretizations and a~new scheme based on the hierarchical vertex-based limiter but using a stricter limiting strategy. Further additions in this work include a~selection of TVD~(total variation diminishing) Runge--Kutta methods of orders one, two, and three employed for time discretization instead of a~simple implicit Euler method used in the first paper.
\par
The rest of this paper is organized as follows:
We introduce the model problem in the remainder of this section and describe its
discretization using the DG~method in Sec.~\ref{sec:discretization}. Section~\ref{sec:slope-lim}
introduces slope limiting algorithms, first, for linear DG~discretizations followed by the 
general order case.
Implementation specific details such as reformulation and assembly of matrix
blocks as well as numerical results are given in Sec.~\ref{sec:implementation}.
All routines mentioned in this work are listed and documented in Sec.~\ref{sec:routines}.
Section~\ref{sec:conclusion} concludes the work and gives future perspectives.

\subsection{Model problem} Let~$J\coloneqq\,(0,t_\mathrm{end})\,$ be a~finite time interval and~$\Omega\subset\IR^2$ 
a~polygonally bounded domain with boundary~$\partial\Omega$.
We consider the \emph{advection equation} in conservative form
\begin{subequations}\label{eq:model}
\begin{equation}\label{eq:model:c}
\partial_t c(t, \vec{x})  +  \div \big( \vec{u}(t, \vec{x})\, c(t, \vec{x}) \big) \;=\; f(t, \vec{x}) \qquad~~\text{in}~J\times\Omega
\end{equation}
with time\,/\,space-varying coefficients~$\vec{u}:J\times\Omega\rightarrow \IR^2$ and $f:J\times\Omega\rightarrow \IR$.
A~prototype application of~\eqref{eq:model:c} is the advective transport in fluids, i.\,e., the movement of a~solute due to the bulk movement of the fluid,
in which case the primary unknown~$c$ denotes the solute concentration, $\vec{u}$~the velocity of the fluid, 
and~$f$ accounts for generation or degradation of~$c$, e.\,g., by chemical reactions.  
Equation~\eqref{eq:model:c} is complemented by the following boundary and initial conditions:
\begin{align}
c  &\;=\; c_\mathrm{D}                         &&\text{on}~J\times{\partial\Omega}_\mathrm{in}(t)\;, &&\label{eq:model:in}\\
c   &\;=\; c^0                                 &&\text{on}~\{0\}\times\Omega&&
\end{align}
\end{subequations}
with inflow boundary~$\partial\Omega_\mathrm{in}(t) \coloneqq \{ \vec{x} \in \partial\Omega \,|\, \vec{u}(t, \vec{x}) \cdot \vec{\nu}(\vec{x}) < 0 \}$ and $\vec{\nu}(\vec{x})$ denoting the outward unit normal. The outflow boundary~$\partial\Omega_\mathrm{out}(t)$ is defined as $\partial\Omega_\mathrm{out}(t) \coloneqq \partial\Omega \setminus \partial\Omega_\mathrm{in}(t)$; $c^0:\Omega\rightarrow\IR^+_0$ and~$c_\mathrm{D}:J\times \partial\Omega_\mathrm{in}(t)\rightarrow\IR_0^+$ are the given initial and Dirichlet boundary data, respectively.

\section{Discretization}\label{sec:discretization}
\subsection{Notation}
Before describing the DG~scheme for~\eqref{eq:model} we introduce some notation; an overview can be found in the Section~\enquote{Index of notation}.
Let $\setT_h=\{T\}$ be a~regular family of non-overlapping partitions of~$\Omega$ into~$K$ closed 
triangles~$T$ of characteristic size~$h$ such that $\displaystyle \overline{\Omega}=\cup T$.
For~$T\in\setT_h$, let $\vec{\nu}_T$ denote the unit normal on~$\partial T$ exterior to~$T$.  
Let $\setE_\Omega$ denote the set of interior edges, $\setE_{\partial\Omega}$ the set of boundary edges, 
and $\setE \coloneqq \setE_\Omega\cup\setE_{\partial\Omega}=\lbrace E\rbrace$ 
the set of all edges (the subscript~$h$ is suppressed here). For an interior edge $E\in\setE_\Omega$ shared by triangles $T^-$ and $T^+$, we define the one-sided values of a~scalar quantity~$w=w(\vec{x})$ on~$E$ by
\begin{equation*}
w^-(\vec{x})\;\coloneqq\;\lim_{\varepsilon \to 0^{+}} w(\vec{x} - \varepsilon\,\vec{\nu}_{T^-})
\qquad\text{and}\qquad
w^+(\vec{x})\;\coloneqq\;\lim_{\varepsilon \to 0^{+}} w(\vec{x} - \varepsilon\,\vec{\nu}_{T^+})\;,
\end{equation*}
respectively.  For a~boundary edge~$E\in\setE_{\partial\Omega}$, only the definition on the left is meaningful. 

\subsection{Variational formulation}
Because of the local nature of the DG~method, we can formulate the variational system of equations on a~triangle-by-triangle basis. 
To do that, we multiply~\eqref{eq:model:c} by
a~smooth test function~$w:T\rightarrow\IR$ and integrate by parts
over element~$T \in \setT_h$. This gives us
\begin{equation*}
\int_{T} w\,\partial_t c(t)\,\dd\vec{x} - \int_{T}\grad w \cdot  \vec{u}(t)\,c(t) \,\dd\vec{x} + \int_{\partial T} w\, \vec{u}(t)\,c(t) \cdot\vec{\nu}_T \,\dd s
\;=\; \int_{T} w\,f(t)\,\dd\vec{x}\;.
\end{equation*}

\subsection{Semi-discrete formulation}\label{sec:semidiscreteformulation}
We denote by~$\IP_p(T)$ the space of complete polynomials of degree at most~$p$ on~$T\in\setT_h$. 
Let
\begin{equation*}
\IP_p(\setT_h) \;\coloneqq\; \Big\{ w_h:\overline{\Omega}\rightarrow \IR\,;~\forall T\in\setT_h,~ {w_h}|_T\in\IP_p(T)\Big\}
\end{equation*}
denote the broken polynomial space on the triangulation~$\setT_h$.
For the semi-discrete formulation, we assume that the coefficient functions (for $t\in J$ fixed) are approximated as:
$\vec{u}_h \in [\IP_p(\setT_h)]^2$ and $f_h(t), c^0_h \in \IP_p(\setT_h)$. 
A~specific way to compute these approximations was given in the first paper of the series~\cite{FESTUNG1}; 
here we use the standard $L^2$-projection into $\IP_p(T)$, therefore the accuracy of this approximation improves with increasing polynomial order~$p$. 
Choosing the same polynomial space for all functions simplifies the implementation and is done in preparation for later applications, in which~$\vec{u}_h$ might be part of the solution of a~coupled system.
Incorporating the boundary condition~\eqref{eq:model:in}, the semi-discrete formulation reads:
\par
Seek $c_h(t)\in\IP_p(\setT_h)$ such that the following holds for 
$t\in J$ and $\forall T^-\in\setT_h,\, \forall w_h\in \IP_p(\setT_h)\,$:
\begin{equation}\label{prob:semidiscrete}
\int_{T^-} w_h\,\partial_t c_h(t)\,\dd\vec{x}
\;- \int_{T^-}\grad w_h \cdot\vec{u}_h(t)\,c_h(t)\,\dd\vec{x}
\;+ \int_{\partial T^-}w_h^-\,\Big(\vec{u}(t)\cdot\vec{\nu}_{T^-}\Big)\,\hat{c}_h(t)\,\dd s
\;=\; \int_{T^-} w_h\,f_h(t)\,\dd\vec{x}\,,
\end{equation}
where the boundary integral is calculated using the upwind-sided value
\begin{equation*}
\hat{c}_h(t,\vec{x})\big|_{\partial T^-} = \left\{
\begin{aligned}
  c_h^-(t,\vec{x})        & \quad\text{if}\quad \vec{u}(t,\vec{x}) \cdot \vec{\nu}_{T^-} \ge 0                                                && \mbox{(outflow from~$T^-$)} \\
  c_h^+(t,\vec{x})        & \quad\text{if}\quad \vec{u}(t,\vec{x}) \cdot \vec{\nu}_{T^-} < 0~\wedge~\vec{x} \notin \partial\Omega_\mathrm{in} && \mbox{(inflow into~$T^-$ from~$T^+$)}\\
  c_\mathrm{D}(t,\vec{x}) & \quad\text{if}\quad  \vec{x} \in \partial\Omega_\mathrm{in}   && \mbox{(inflow into~$T^-$ over~$\partial\Omega_\mathrm{in}$)}
\end{aligned}
\right\}.
\end{equation*}
Note that we did not use the approximate representation of the velocity~$\vec{u}_h$ in the boundary integral. 
This is due to the fact that the $L^2$-projection on elements may have poor approximation quality on edges and generally produces different values on both sides of the edge ultimately leading to different upwind-sided values and inconsistent flux approximations. 
Instead we evaluate the normal velocity $\vec{u} \cdot \vec{\nu}_T$ in each quadrature point analytically and use the result for both the numerical integration and the determination of the upwind direction as will be demonstrated in Sec.~\ref{sec:defBlocks:E}.
\par
Thus far, we used an~algebraic indexing style.
In the remainder, we switch to a~mixture of algebraic and numerical style: for instance, 
$E_{kn}\in\partial T_k\cap\setE_\Omega$ means \emph{all possible} combinations of element indices~$k\in\{1,\ldots,K\}$ 
and local edge indices $n\in\{1,2,3\}$ such that $E_{kn}$ lies in~$\partial T_k\cap\setE_\Omega$. 
This implicitly fixes the numerical indices which accordingly can be used to index matrices or arrays.
\par
We use a~bracket notation followed by a~subscript to index matrices and multidimensional arrays.  
Thus, for an $n$-dimensional array~$\vecc{X}$, the symbol $[\vecc{X}]_{i_1,\ldots,i_n}$ 
stands for the component of~$\vecc{X}$ with index~$i_l$ in the $l$-th dimension.  
As in \MatOct, a~colon is used to abbreviate all indices within a~single dimension.  
For example, $[\vecc{X}]_{:,:,i_3,\ldots,i_n}$ is a~two-dimensional array\,/\,matrix.

\subsubsection{Local basis representation}\label{sec:basisfunctions}
In contrast to globally continuous basis functions mostly used by the continuous finite element method, 
the DG~basis functions have no continuity constraints across triangle boundaries.  
Thus a~standard DG~\emph{basis function}~$\vphi_{ki}: \overline{\Omega}\rightarrow \IR$ is only supported 
on the triangle~$T_k\in\setT_h$ (i.\,e., $\vphi_{ki}=0$ on $\overline{\Omega}\smallsetminus T_k$) 
and can be defined arbitrarily while ensuring
\begin{equation}\label{eq:defOfIPpAndN}
\forall k\in\{1,\ldots,K\}\,,\quad  \IP_p(T_k)\;=\;\mathrm{span}\,\big\{ \vphi_{ki} \big\}_{i\in\{1,\ldots,N_p\}}\;,
\qquad
\text{where}
\quad
N_p \;\coloneqq\; \frac{(p+1)(p+2)}{2} \;=\; \begin{pmatrix}p+2\\p\end{pmatrix}
\end{equation}
is the number of \emph{local degrees of freedom}.  
Note that~$N_p$ may in general vary from triangle to triangle, but, for simplicity, we assume here a~uniform polynomial degree~$p$ for every triangle and abbreviate~$N\coloneqq N_p$.
Clearly, the number of global degrees of freedom equals~$KN$. 
Closed-form expressions for orthonormal basis functions on the reference triangle~$\hat{T}$ (cf.~Sec.~\ref{sec:transformationtoThat}) employed in our implementation up to order two can be found in our first paper~\cite{FESTUNG1}.
The basis functions up to order four are provided in the routine~\code{phi} and their gradients in~\code{gradPhi}.  
Bases of even higher order can be constructed, e.\,g., with the Gram--Schmidt algorithm or by using a~three-term recursion relation---the latter is unfortunately not trivial to derive in the case of triangles. 
Note that these so-called \emph{modal} basis functions~$\hat{\vphi}_i$ do \emph{not} posses interpolation properties at nodes unlike Lagrangian\,/\,nodal basis functions, which are often used by the continuous finite element or nodal DG~methods.
\par
The local concentration~$c_h$ and the local velocity~$\vec{u}_h$ on~$T_k\in\setT_h$ can be represented in terms of the local basis~$\{\varphi_{ki}\}_{i\in\{1,\ldots,N\}}$:
\begin{equation*}
c_h(t,\vec{x})      \big|_{T_k}\eqqcolon \sum_{j=1}^N C_{kj}(t)\,\vphi_{kj}(\vec{x})\,, 
\qquad 
\vec{u}_h(t,\vec{x})  \big|_{T_k}\eqqcolon \sum_{j=1}^N \sum_{m=1}^2 U_{kj}^m(t)\,\vec{e}_m\vphi_{kj}(\vec{x})\,,
\end{equation*}
where~$\vec{e}_m$ denotes the $m$-th~unit vector~in~$\IR^2$. 
We condense the coefficients associated with unknowns into two-dimensional 
arrays~$\vecc{C}(t)$ such that~$C_{kj}(t) \allowbreak\coloneqq  [\vecc{C}(t)]_{k,j}$, etc. 
The symbol~$[\vecc{C}]_{k,:}$ is called \emph{local representation matrix}
of~$c_h$ on~$T_k$ with respect to the basis~$\big\{ \vphi_{ki} \big\}_{i\in\{1,\ldots,N\}}$.
In~a~similar way, we express the coefficient functions as linear combinations of the basis functions: 
On~$T_k\in\setT_h$, we use the local representation matrices~$[\vecc{C}^0]_{k,:}$~for 
$c_h^0$ and~$[\vecc{F}]_{k,:}$ for~$f_h$. 

\subsubsection{System of equations}
Testing~\eqref{prob:semidiscrete} with 
$w_h = \vphi_{ki}$ for $i\in\{1,\ldots,N\}$ yields a~\emph{time-dependent system of equations} 
whose contribution from~$T_k$ (identified with $T_{k^-}$ in boundary integrals) reads
\begin{equation}
\begin{multlined}
\underbrace{\sum_{j=1}^N\partial_t C_{kj}(t)\int_{T_k}\vphi_{ki}\,\vphi_{kj}\,\dd\vec{x}}_{I}
-\underbrace{\sum_{j=1}^N C_{kj}(t) \sum_{l=1}^N\sum_{m=1}^2U_{kl}^m(t)\int_{T_k} \partial_{x^m} \vphi_{ki}\,\vphi_{kl}\,\vphi_{kj}\,\dd\vec{x}}_{\II} \\
+ \underbrace{\int_{\partial T_{k^-}} \vphi_{k^-i}\, \Big( \vec{u}(t) \cdot \vec{\nu}_{k^-} \Big)
\left\{\begin{aligned}
\sum_{j=1}^N C_{k^-j}(t)\,\vphi_{k^-j}&~~\text{if}~~\vec{u}(t)\cdot\vec{\nu}_{k^-}\ge0 \\
\sum_{j=1}^N C_{k^+j}(t)\,\vphi_{k^+j}&~~\text{if}~~\vec{u}(t)\cdot\vec{\nu}_{k^-}<0\,\wedge\, \vec{x} \notin \partial\Omega_\mathrm{in} \\
  c_\mathrm{D}(t)                     &~~\text{if}~~ \vec{x} \in \partial\Omega_\mathrm{in} 
\end{aligned}\right\}\,\dd s}_{\III}
\;=\;
\underbrace{\sum_{l=1}^N F_{kl}(t) \int_{T_k}\vphi_{ki}\,\vphi_{kl}\,\dd\vec{x}}_{\IV}\;,
\end{multlined}
\label{eq:spacediscretesystem}
\end{equation}
where we abbreviated $\vec{\nu}_{T_k}$ by $\vec{\nu}_k$. 
Written in matrix form, system~\eqref{eq:spacediscretesystem} is then given by
\begin{equation}\label{eq:timeDepSystem}
\vecc{M}\,\partial_t\vec{C}
+ \underbrace{\left(-\vecc{G}^1 - \vecc{G}^2 + \vecc{R}\right)}_{\eqqcolon\;\vecc{A}(t)} \vec{C}
= \underbrace{\vec{L} - \vec{K}_\mathrm{D}}_{\eqqcolon\;\vec{V}(t)}
\end{equation}
with the representation vector
\begin{align*}
\vec{C}(t) &\;\coloneqq\;
\transpose{\begin{bmatrix}
C_{11}(t) & \cdots & C_{1N}(t) & \cdots &\cdots & C_{K1}(t) & \cdots & C_{KN}(t)
\end{bmatrix}}\;.
\end{align*}
The block matrices and the right-hand side vectors of~\eqref{eq:timeDepSystem} are described in~Sections~\ref{sec:defBlocks:T} and~\ref{sec:defBlocks:E}. 
Note that all blocks except for the mass matrix $\vecc{M}$ are time-dependent (we have suppressed the time~arguments here).

\subsubsection[Contributions from area terms I, II, IV]{Contributions from area terms~$I$, $\II$, $\IV$}\label{sec:defBlocks:T}
The matrices in the remainder of this section have sparse block structure; by giving definitions 
for non-zero blocks we tacitly assume a~zero fill-in.
The \emph{mass matrix}~$\vecc{M}\in\IR^{KN\times KN}$ in term 
$I$ is defined component-wise as
\begin{equation*}
[\vecc{M}]_{(k-1)N+i, (k-1)N+j} \;\coloneqq\; \int_{T_k}\vphi_{ki}\,\vphi_{kj}{\color{red}\,\dd\vec{x}}\;.
\end{equation*}
Since the basis functions~$\vphi_{ki}$, $i\in\{1,\ldots,N\}$ are supported only on~$T_k$, 
$\vecc{M}$ has a~block-diagonal structure 
\begin{equation}\label{eq:globMlocM}
\vecc{M} \;=\;
\begin{bmatrix}
\vecc{M}_{T_1} &          & \\
               & ~\ddots~ & \\
               &          & \vecc{M}_{T_K}
\end{bmatrix}
\qquad\text{with}\qquad
\vecc{M}_{T_k} \;\coloneqq\;
\int_{T_k}\begin{bmatrix}
\vphi_{k1}\,\vphi_{k1} & \cdots & \vphi_{k1}\,\vphi_{kN} ~\\
 \vdots                & \ddots & \vdots \\
\vphi_{kN}\,\vphi_{k1} & \cdots & \vphi_{kN}\,\vphi_{kN} 
\end{bmatrix}\,\dd\vec{x}\;,
\end{equation}
i.\,e., it consists of $K$~\emph{local mass matrices}~$\vecc{M}_{T_k}\in\IR^{N\times N}$.  
Henceforth, we write~$\vecc{M} = \diag \big(\vecc{M}_{T_1},\ldots,\vecc{M}_{T_K}\big)$.
\par
The block matrices~$\vecc{G}^m\in\IR^{KN\times KN}, \;m\in\{1,2\}$ from term~$\II$ are given by
\begin{align*}
[\vecc{G}^m]_{(k-1)N+i, (k-1)N+j} \;\coloneqq\; 
 \sum_{l=1}^N U^m_{kl}(t) \int_{T_k} \partial_{x^m} \vphi_{ki}\,\vphi_{kl} \,\vphi_{kj}\,\dd\vec{x}\,.
\end{align*}
Similarly to~$\vecc{M}$, the matrices~$\vecc{G}^m = \diag \big(\vecc{G}^m_{T_1},\ldots,\vecc{G}^m_{T_K}\big)$ are block-diagonal with local matrices
\begin{equation}\label{eq:locGm}
\vecc{G}^m_{T_k} \;\coloneqq\;
\sum_{l=1}^N U^m_{kl}(t)
\int_{T_k}\begin{bmatrix}
\partial_{x^m}\vphi_{k1}\,\vphi_{kl}\,\vphi_{k1} & \cdots & \partial_{x^m}\vphi_{k1}\,\vphi_{kl}\,\vphi_{kN}\\
 \vdots                & \ddots & \vdots \\
\partial_{x^m}\vphi_{kN}\,\vphi_{kl}\,\vphi_{k1} & \cdots & \partial_{x^m}\vphi_{kN}\,\vphi_{kl}\,\vphi_{kN} 
\end{bmatrix}\,\dd\vec{x}\;.
\end{equation}
\par
Vector~$\vec{L}(t)$ resulting from $\IV$ is obtained 
by multiplication of the representation vector of~$f_h(t)$ to the global mass matrix:
\begin{equation*}
\vec{L}(t) \;=\; \vecc{M}\, \transpose{\begin{bmatrix}
F_{11}(t) & \cdots & F_{1N}(t) & \cdots &\cdots & F_{K1}(t) & \cdots & F_{KN}(t)
\end{bmatrix}}\;.
\end{equation*}

\subsubsection[Contributions from edge term III]{Contributions from edge term~$\III$}\label{sec:defBlocks:E}


\begin{figure}[t!]\centering
\includegraphics{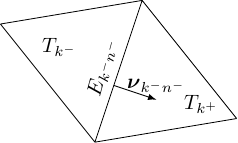}
\captionsetup{justification=raggedright,singlelinecheck=false}
\caption{Two triangles adjacent to edge~$E_{k^-n^-}$. It holds:~$E_{k^-n^-} = E_{k^+n^+}$ and $\vec{\nu}_{k^-n^-} = -\vec{\nu}_{k^+n^+}$.}
\label{fig:T1T2}
\end{figure}
\paragraph{Interior Edges $\setE_\Omega$}
In this section, we consider a~fixed triangle $T_k=T_{k^-}$ with an~interior 
edge~$E_{k^-n^-} \in\partial T_{k^-}\cap\setE_\Omega = \partial T_{k^-} \cap \partial T_{k^+}$ 
shared by an adjacent triangle $T_{k^+}$ and associated with fixed local edge indices $n^-,n^+ \in \{1,2,3\}$~(cf.~Fig.~\ref{fig:T1T2}).
\par
For a~fixed index~$i\in\{1,\ldots,N\}$, we have a~contribution for $\vphi_{k^-i}$ in a~block matrix $\vecc{R}_\Omega\in\IR^{KN\times KN}$
\begin{equation*}
\int_{E_{k^-n^-}} \vphi_{k^-i} \, \Big(\vec{u}(t) \cdot \vec{\nu}_{k^-n^-}\Big) \,
\begin{Bmatrix}
\sum_{j=1}^N C_{k^-j}(t)\,\vphi_{k^-j} &\text{if}& \vec{u}\cdot\vec{\nu}_{k^-n^-}\ge0\\
\sum_{j=1}^N C_{k^+j}(t)\,\vphi_{k^+j} &\text{if}& \vec{u}\cdot\vec{\nu}_{k^-n^-}<0
\end{Bmatrix}\,\dd s\;.
\end{equation*}
This means that, depending on the direction of the velocity field~$\vec{u}(t,\vec{x})$, we obtain entries in the diagonal or off-diagonal blocks of $\vecc{R}_\Omega$.
Entries in diagonal blocks are then component-wise given by
\begin{subequations}
\begin{equation}\label{eq:globRInterior:diag}
[\vecc{R}_\Omega]_{(k-1)N+i,(k-1)N+j} \;\coloneqq
\sum_{E_{kn}\in\partial T_k\cap\setE_\Omega} \int_{E_{kn}} \vphi_{ki}\,\vphi_{kj}\,
\Big(\vec{u}\cdot\vec{\nu}_{kn}\Big) \,\delta_{\vec{u}\cdot\vec{\nu}_{kn}\ge0}\,\dd s
\quad\text{with}\quad
\delta_{\vec{u}\cdot\vec{\nu}_{kn}\ge0}\,(t,\vec{x}) \;\coloneqq\; \begin{Bmatrix}
1 &\text{if}& \vec{u}(t,\vec{x})\cdot\vec{\nu}_{kn} \ge 0 \\
0 &\text{if}&\vec{u}(t,\vec{x})\cdot\vec{\nu}_{kn} < 0
\end{Bmatrix}\,.
\end{equation}
Entries in off-diagonal blocks in~$\vecc{R}$ are possibly non-zero only for pairs of triangles $T_{k^-}$, $T_{k^+}$ with $\partial T_{k^-}\cap\partial T_{k^+}\neq\emptyset$ and read
\begin{equation}\label{eq:globRInterior:offdiag}
[\vecc{R}_\Omega]_{(k^--1)N+i,(k^+-1)N+j} \;\coloneqq\;
\int_{E_{k^-n^-}} \vphi_{k^-i}\,\vphi_{k^+j}\,
\Big(\vec{u}\cdot\vec{\nu}_{k^-n^-}\Big) \,\delta_{\vec{u}\cdot\vec{\nu}_{k^-n^-}<0}\,\dd s 
\quad\text{with}\quad
\delta_{\vec{u}\cdot\vec{\nu}_{k^-n^-}<0} \;\coloneqq\;
1 - \delta_{\vec{u}\cdot\vec{\nu}_{k^-n^-}\ge0}\;.
\end{equation}
\end{subequations}

\paragraph{Boundary Edges $\setE_{\partial\Omega}$}
Similarly to interior edges we have contributions for a~boundary edge $E_{kn}\in\partial T_k \cap \setE_{\partial\Omega}$
\begin{equation*}
\int_{E_{kn}} \vphi_{ki} \, \Big(\vec{u}(t) \cdot \vec{\nu}_{kn}\Big) \,
\begin{Bmatrix}
\sum_{j=1}^N C_{kj}(t) \, \vphi_{kj} &\text{if}& \vec{u}(t)\cdot\vec{\nu}_{kn} \ge 0 \\
c_\mathrm{D}(t)                      &\text{if}& \vec{u}(t)\cdot\vec{\nu}_{kn} < 0
\end{Bmatrix}\,\dd s\;.
\end{equation*}
These consist of entries in the block diagonal matrix~$\vecc{R}_{\partial\Omega}\in\IR^{KN\times KN}$ 
\begin{equation}\label{eq:globRBoundary:diag}
[\vecc{R}_{\partial\Omega}]_{(k-1)N+i,(k-1)N+j} \;\coloneqq
\sum_{E_{kn}\in\partial T_k\cap\setE_{\partial\Omega}} \int_{E_{kn}} \vphi_{ki}\,\vphi_{kj}\,
\Big(\vec{u}\cdot\vec{\nu}_{kn}\Big) \,\delta_{\vec{u}\cdot\vec{\nu}_{kn}\ge0}\,\dd s \;,
\end{equation}
and in the right-hand side vector $\vec{K}_\mathrm{D}\in\IR^{KN}$
\begin{equation}\label{eq:globKD}
[\vec{K}_\mathrm{D}]_{(k-1)N+i} \;\coloneqq\;
\sum_{E_{kn}\in\partial T_k\cap\setE_{\partial\Omega}} \int_{E_{kn}} \vphi_{ki}\,c_\mathrm{D}(t)
\Big(\vec{u}\cdot\vec{\nu}_{kn}\Big) \,\delta_{\vec{u}\cdot\vec{\nu}_{kn}<0}\,\dd s \;.
\end{equation}
We combine the block matrices $\vecc{R}_\Omega, \vecc{R}_{\partial\Omega}$ into a~block matrix $\vecc{R}\in\IR^{KN\times KN}$
\begin{equation}\label{eq:globR}
  \vecc{R} \;\coloneqq\; \vecc{R}_\Omega + \vecc{R}_{\partial\Omega}\,.
\end{equation}
Since the definition of entries in the diagonal blocks in Eqns.~\eqref{eq:globRInterior:diag},~\eqref{eq:globRBoundary:diag} is the same for both matrices differing only in the set of edges included in the sum, we can disregard the fact whether they are interior or boundary edges and simply assemble the entries for all $E_{kn} \in \partial T_k$.

\subsection{Time discretization}\label{sec:timediscretization}
The system~\eqref{eq:timeDepSystem} is equivalent to
\begin{equation} \label{eq:fullSystem}
\vecc{M}\partial_t\vec{C}(t) \;=\; \vec{V}(t) - \vecc{A}(t)\,\vec{C}(t) \;\eqqcolon\; \vec{S}\Big(\vec{C}(t), t\Big)
\end{equation}
with $\vecc{A}(t) \in \IR^{KN\times KN}$ and 
right-hand-side vector~$\vec{V}(t) \in \IR^{KN}$ as defined in~\eqref{eq:timeDepSystem}.
\par
We discretize system~\eqref{eq:fullSystem} in time using TVD (total variation diminishing) Runge--Kutta methods~\cite{GottliebShu1998} of orders one, two, and three, which are representatives of the class of SSP~(strong stability preserving) Runge--Kutta methods~\cite{GottliebShu2001}. 
The advantage of using a~time stepping algorithm of such type lies in the guaranteed preservation of the monotonicity of the solution if the~DG~discretization is also post processed by a~slope limiting method.
\par
Let $0= t^1<t^2<\ldots< t_\mathrm{end}$ be a~not necessarily equidistant decomposition 
of the time interval~$J$ and let $\Delta t^n \coloneqq t^{n+1} -t^n$ denote the time step size.
The update scheme of the $s$-step Runge--Kutta method is given by
\begin{equation}\label{eq:SSP-RK}
\begin{array}{lll}
\vec{C}^{(0)} &=\; \vec{C}^{n} \,, \\
\vec{C}^{(i)} &=\; \omega_i\,\vec{C}^{n} + (1-\omega_i)\,\left( \vec{C}^{(i-1)} 
+ \Delta t^n\,\vecc{M}^{-1} \vec{S}^{n + \delta_i} \right) \,, \quad \text{for}~i = 1,\ldots,s\,,\\
\vec{C}^{n+1} &=\; \vec{C}^{(s)} \,,
\end{array}
\end{equation}
where we abbreviated~$\vec{C}^n\coloneqq \vec{C}(t^n)$ and $\vec{S}^{n + \delta_i}\coloneqq \vec{S}(\vec{C}^{(i-1)}, t^n + \delta_i \Delta t^n)$ with coefficients
\begin{align*}
s = 1\,: &&&\omega_1 = 0\,,                                          && \delta_1 = 0 \,.\\
s = 2\,: &&&\omega_1 = 0\,,\; \omega_2 = 1/2 \,,                     && \delta_1 = 0 \,,\; \delta_2 = 1 \,.\\
s = 3\,: &&&\omega_1 = 0\,,\; \omega_2 = 3/4 \,,\; \omega_3 = 1/3\,, && \delta_1 = 0 \,,\; \delta_2 = 1 \,,\; \delta_3 = 1/2 \,.
\end{align*}
When possible, we choose the order of the time-discretization to be~$p+1$, with~$p$ being the spatial approximation order, in order to avoid the temporal discretization error dominating the spatial one.
{\color{red} The chosen SSP~Runge--Kutta methods are optimal in the sense that they achieve $p$-th order with $p$ stages.}
Unfortunately, no optimal SSP~Runge--Kutta methods higher than order three are known~\cite{GottliebShu2001} {\color{red}(there exist, however, non-optimal higher order schemes)}, which is why we restrict ourselves to orders one to three for the time discretization.

\section{Slope limiting}\label{sec:slope-lim}
Slope limiters are a~technique to prevent the onset of spurious oscillations that violate the monotonicity preserving property of the piecewise constant part of a~DG~solution by means of restricting some of the degrees of freedom (generally linear and superlinear) to certain bounds and thus eliminating over- and undershoots.
All limiting procedures utilize the fact that the lowest order (piecewise constant) part of a~DG~solution 
in explicit TVD time stepping schemes is guaranteed 
to preserve the monotonicity of the solution and produce no spurious extrema.
Using this physically consistent but numerically not very accurate solution part, all slope limiters attempt to modify the full higher order DG~solution in a~suitable way---on the one hand, to prevent any oscillations and, on the other hand, to preserve as much of the accuracy as possible.
The key differences in slope limiters affect the limiting stencil used (edge neighbors, node neighbors, neighbors of the neighbors, etc.), presence of ad hoc parameters, the amount of the introduced numerical diffusion, a~strict or less strict preservation of the monotonicity, and the degree of solution degradation in smooth extrema.
\par
Whereas a large literature on slope limiting for piecewise linear DG~discretizations exists~\cite{CockburnShuRKDG21989, Krivodonova2004, Tu2005}, the limiting of DG~solutions with $p \ge 2$ is a~much less explored area.
The traditional approach to dealing with superlinear DG~solutions~\cite{Michoski2011} has been based on ignoring all higher order degrees of freedom on elements on which linear limiting is active.
Other methods require a much larger stencil~\cite{Yang2009, Zhang2012} to provide enough information for the reconstruction of higher order derivatives.
The hierarchical vertex-based limiters of Kuzmin~\cite{Kuzmin2010, Kuzmin2012} represent a~computationally efficient scheme easily extendable to any discretization order and supporting fully unstructured meshes.
These limiters do not guarantee the strict monotonicity of the DG solution, but the violations are small and may be further reduced by simple modifications described in Sec.~\ref{sec:slope-lim:strict}.

\subsection{Taylor basis representation}\label{sec:taylor:basis}
Many limiting procedures rely on some fundamental properties of a~certain choice of basis, in our case the 2D~Taylor basis, which we introduce in a~way similar to Kuzmin~\cite{Kuzmin2010}.
Consider the 2D~Taylor series expansion of a~local solution~$c_h \in \IP_p(T_k)$,
\begin{equation} \label{eq:taylor:exp}
c_h(\vec{x}) \;= \sum_{0\le|\vec{a}|\le p} 
\partial^\vec{a} c_h(\vec{x}_{k\mathrm{c}}) 
\frac{(\vec{x}-\vec{x}_{k\mathrm{c}})^\vec{a}}{\vec{a}!}  \qquad \text{on}~T_k\in\setT_h
\end{equation}
about the centroid $\vec{x}_{k\mathrm{c}} = \transpose{[x_{k\mathrm{c}}^1,x_{k\mathrm{c}}^2]}$ of~$T_k\in\setT_h$ with a~two-dimensional multi-index~$\vec{a} = \transpose{[a^1, a^2]} \in \IN_0^2$,
where we use some standard notation for multi-indices $\vec{a},\vec{b} \in \IN_0^2$ and $\vec{x}\in\IR^2$:
\begin{align*}
\vec{a} \pm \vec{b} &= \transpose{[a^1 \pm b^1, a^2 \pm b^2]} \;, &
|\vec{a}|           &\coloneqq a^1 + a^2                  \;, &
\vec{a}!            &\coloneqq a^1! a^2!                  \;, \\
\vec{x}^\vec{a}     &\coloneqq (x^1)^{a^1} (x^2)^{a^2} \;, &
\partial^\vec{a}    &\coloneqq \partial^{|\vec{a}|} \big/ \partial(x^1)^{a^1}\,\partial(x^2)^{a^2} \;. & &
\end{align*}
For $v:T_k\rightarrow \IR$, let $\overline{v} \coloneqq \frac{1}{|T_k|}\int_{T_k} v(\vec{x})\,\dd\vec{x}$ denote the integral mean of~$v$ on~$T_k$.
We express~\eqref{eq:taylor:exp} in the equivalent form~\cite{YangWang2009,Michalak2008,Luo2008}
\begin{equation}\label{eq:taylor:exp2}
c_h(\vec{x}) \;=\; \overline{c}_h 
+ \frac{\partial c_h}{\partial x^1}(\vec{x}_{k\mathrm{c}})  (x^1 - x_{k\mathrm{c}}^1)
+ \frac{\partial c_h}{\partial x^2}(\vec{x}_{k\mathrm{c}})  (x^2 - x_{k\mathrm{c}}^2)
+ \sum_{2\le|\vec{a}|\le p} 
\partial^\vec{a} c_h(\vec{x}_{k\mathrm{c}}) 
\frac{(\vec{x}-\vec{x}_{k\mathrm{c}})^\vec{a} - \overline{(\vec{x}-\vec{x}_{k\mathrm{c}})^\vec{a}}}{\vec{a}!}
\quad\text{on}~T_k\in\setT_h\,.
\end{equation}
Note that varying any terms in~\eqref{eq:taylor:exp2} except $\overline{c}_h$ does not affect the mean of $c_h(\vec{x})$ over $T_k$.
\par
To be able to identify each term in expansions~\eqref{eq:taylor:exp},\,\eqref{eq:taylor:exp2} by a~consecutive index, we introduce a~linear index mapping $I: \IN_0^2 \rightarrow \IN$ corresponding to any two-dimensional multi-index $\vec{a}\in\IN_0^2$ as
\begin{equation}\label{eq:multiindex:lin}
  I(\vec{a}) \;=\; N_{|\vec{a}|-1}+a^2+1  \;=\;  \frac{|\vec{a}|(|\vec{a}|+1)}{2}+a^2+1
\end{equation}
with $N_p = \dim \IP_p(T)$ as defined in~\eqref{eq:defOfIPpAndN}.
We implicitly define $\vec{a}_j$ such that $\forall j\in\IN$, $I(\vec{a}_j) = j$. 
The linear indices, polynomial degrees, and corresponding multi-indices up to order four are listed in Table~\ref{tab:multi-index}.
\begin{table}[t!]
\setlength{\tabcolsep}{1.45mm}
\centering
\begin{tabularx}{\linewidth}{@{}c|c|cc|ccc|cccc|ccccc@{}}
\toprule
$p$            & 0       &\multicolumn{2}{c|}{1} & \multicolumn{3}{c|}{2}      & \multicolumn{4}{c|}{3}             & \multicolumn{5}{c}{4} \\
$I(\vec{a}_i)$ & 1       & 2       & 3           & 4       & 5       & 6       & 7       & 8       & 9     & 10  & 11  & 12  & 13  & 14 & 15   \\
$\transpose{\vec{a}_i}$    & $[0,0]$ & $[1,0]$ & $[0,1]$   & $[2,0]$ & $[1,1]$ & $[0,2]$ & $[3,0]$ & $[2,1]$ & $[1,2]$ & $[0,3]$ & $[4,0]$ & $[3,1]$ & $[2,2]$ & $[1,3]$ & $[0,4]$ \\ 
\bottomrule
\end{tabularx}
\caption{Multi-indices \emph{(bottom)}, linear indices \emph{(middle)}, and corresponding polynomial degrees in the Taylor basis~\eqref{eq:taylor:basis} \emph{(top)}.}
\label{tab:multi-index}
\end{table}
This leads to the following definition of the local Taylor basis~\cite{YangWang2009}:
\begin{align}\label{eq:taylor:basis}
\phi_{k1} &= 1 \;, &
\phi_{k2} &= \frac{x_k^1 - x_{k\mathrm{c}}^1}{\Delta (x_k^1)} \;, &
\phi_{k3} &= \frac{x_k^2 - x_{k\mathrm{c}}^2}{\Delta (x_k^2)} \;, &
\phi_{ki} &= \frac{(\vec{x} - \vec{x}_{k\mathrm{c}})^{\vec{a}_i} - \overline{(\vec{x} - \vec{x}_{k\mathrm{c}})^{\vec{a}_i}}}{\vec{a}_i! \, (\Delta \vec{x}_k)^{\vec{a}_i}} \quad\text{for}\quad i \ge 4 \,.
\end{align}
As opposed to the DG~basis~$\{\varphi_{kj}\}$ (cf.~\cite{FESTUNG1}, Sec.~2.4.1), the basis~$\{\phi_{kj}\}$ cannot be defined on a~reference element~$\hat{T}$.
The scaling by $\Delta\vec{x}_k = \transpose{\big[\Delta (x_k^1), \Delta (x_k^2)\big]}$ with $\Delta (x_k^j) \coloneqq (x_{k,\mathrm{max}}^j - x_{k,\mathrm{min}}^j)/2$, where $x_{k,\mathrm{max}}^j \coloneqq \max_{i\in\{1,2,3\}} x_{ki}^j$ and $x_{k,\mathrm{min}}^j \coloneqq \min_{i\in\{1,2,3\}} x_{ki}^j$ are the minimum and maximum values of the corresponding spatial coordinates on~$T_k$, is introduced to obtain a~better conditioned operator~\cite{YangWang2009}.
The Taylor degrees of freedom are now proportional to the cell mean values $\overline{c}_h$ and derivatives of $c_h$ at the centroid $\vec{x}_{k\mathrm{c}}$
\begin{equation}\label{eq:taylor:representation}
c_h(\vec{x}) = \overline{c}_h\,\phi_{k1} 
+ \left(\frac{\partial c_h}{\partial x^1}(\vec{x}_{k\mathrm{c}}) \, \Delta\left(x_k^1\right)\right)\,\phi_{k2}(\vec{x})
+ \left(\frac{\partial c_h}{\partial x^2}(\vec{x}_{k\mathrm{c}}) \, \Delta\left(x_k^2\right)\right)\,\phi_{k3}(\vec{x})
+ \sum_{i=4}^{N_p} \Big( \partial^{\vec{a}_i}c_h(\vec{x}_{k\mathrm{c}})\, \left(\Delta\vec{x}_k\right)^{\vec{a}_i} \Big) \, \phi_{ki}(\vec{x}) 
\quad\text{on}~T_k\in\setT_h\,.
\end{equation}
\par
Note that the Taylor basis is non-orthogonal on triangular meshes~\cite{Kuzmin2010}, but the cell means are still decoupled from the other degrees of freedom since
\begin{equation*}
  \int_{T_k} \phi_{k1}^2\,\dd\vec{x} = |T_k| \;,\qquad
  \int_{T_k} \phi_{k1}\,\phi_{kj}\,\dd\vec{x} = 0 \quad\text{for}\quad j>1\,.
\end{equation*}
\par
To transform a~function $c_h$ from the modal basis representation with representation matrix $\vecc{C}(t) \in \IR^{K\times N}$---as described in Sec.~\ref{sec:basisfunctions}---into a~Taylor basis representation with representation matrix $\vecc{C}^\mathrm{Taylor}(t) \in \IR^{K\times N}$, we employ the $L^2$-projection defined locally for $T_k \in \setT_h$ by
\begin{equation*}
\forall w_h \in \IP_p(T_k) \,,
\quad
\int_{T_k} w_h \left(\sum_{j=1}^N C_{kj}(t)\,\vphi_{kj}\right)\,\dd\vec{x} =
\int_{T_k} w_h \left(\sum_{j=1}^N C_{kj}^\mathrm{Taylor}(t)\,\phi_{kj}\right)\,\dd\vec{x} \,.
\end{equation*}
Choosing $w_h = \vphi_{ki}$ for $i \in \{1,\ldots,N\}$ we obtain
\begin{equation*}
\sum_{j=1}^N C_{kj}(t) \int_{T_k} \vphi_{ki} \,\vphi_{kj}\, \dd\vec{x} =
\sum_{j=1}^N C_{kj}^\mathrm{Taylor}(t) \int_{T_k} \vphi_{ki}\, \phi_{kj}\, \dd\vec{x}
\quad\Leftrightarrow\quad 
\vecc{M}_{T_k} \left[\vecc{C}\right]_{k,:} =
\vecc{M}^\mathrm{DG,Taylor}_{T_k} \left[\vecc{C}^\mathrm{Taylor}\right]_{k,:}
\end{equation*}
with the local mass matrix $\vecc{M}_{T_k}$ as defined in Eq.~\ref{eq:globMlocM} and the local basis transformation matrix
\begin{equation}\label{eq:MDGTaylor}
\vecc{M}_{T_k}^\mathrm{DG,Taylor} \;\coloneqq\;
\int_{T_k}\begin{bmatrix}
\vphi_{k1}\,\phi_{k1} & \cdots & \vphi_{k1}\,\phi_{kN} ~\\
\vdots                & \ddots & \vdots \\
\vphi_{kN}\,\phi_{k1} & \cdots & \vphi_{kN}\,\phi_{kN} 
\end{bmatrix}\,\dd\vec{x}\;.
\end{equation}
Using $\vecc{M}^\mathrm{DG,Taylor} \coloneqq \mathrm{diag}\left(\vecc{M}_{T_1}^\mathrm{DG,Taylor},\ldots,\vecc{M}_{T_K}^\mathrm{DG,Taylor}\right)$ and representation \emph{vectors} $\vec{C},\vec{C}^\mathrm{Taylor} \in \IR^{KN}$ we obtain a~linear system of equations
\begin{equation}\label{eq:taylor:trafo}
\vecc{M}\,\vec{C} = \vecc{M}^\mathrm{DG,Taylor}\vec{C}^\mathrm{Taylor}\,,
\end{equation}
which can be employed to transform back-and-forth between both bases.

\subsection{Linear vertex-based limiter}\label{sec:slope-lim:linear}
\begin{figure}
  \centering
  \includegraphics{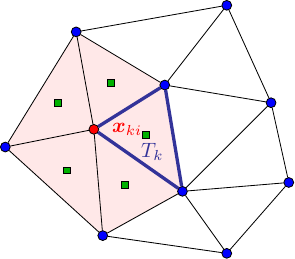}
  \caption{The neighborhood of a~vertex~$\vec{x}_{ki}\in T_k$ \emph{(red circle)} considered in Eq.~\eqref{eq:slope-lim:lin:cond} consists of the patch of elements containing $\vec{x}_{ki}$ \emph{(red area)}. The bounds $c_{ki}^\mathrm{min}$, $c_{ki}^\mathrm{max}$ from Eq.~\eqref{eq:slope-lim:lin:bounds} are determined from all centroid values \emph{(green squares)} within this neighborhood.}
  \label{fig:vertex-lim}
\end{figure}

Kuzmin~\cite{Kuzmin2012,Kuzmin2010} and Aizinger~\cite{Aizinger2011} described the vertex-based limiter, which is based on the Barth--Jespersen limiter~\cite{BarthJespersen1989} and improved it further by taking the bounds from all elements containing the vertex instead of taking only edge neighbors of the cell.
The goal is to determine the maximum admissible slope for a~linear reconstruction of the form
\begin{equation*}
c_h(\vec{x}) = c_{k\mathrm{c}} + \alpha_{ke} \, \grad c_h (\vec{x}_{k\mathrm{c}}) \cdot (\vec{x} - \vec{x}_{k\mathrm{c}}) \,,
\qquad 0\le \alpha_{ke}\le 1\,,\quad \vec{x} \in T_k\,,
\end{equation*}
where we abbreviated the function value $c_{k\mathrm{c}} \coloneqq c_h(\vec{x}_{k\mathrm{c}})$ in the centroid $\vec{x}_{k\mathrm{c}}$.
The correction factor~$\alpha_{ke}$ is chosen such that above reconstruction is bounded in all vertices $\vec{x}_{ki} \in T_k$ by the minimum and maximum centroid values of all elements containing~$\vec{x}_{ki}$, that is
\begin{equation}\label{eq:slope-lim:lin:cond}
\forall T_k \in \setT_h \,, \forall i\in\{1,2,3\} \,,\quad c_{ki}^\mathrm{min} \le c_h(\vec{x}_{ki}) \le c_{ki}^\mathrm{max}
\end{equation}
with
\begin{equation}\label{eq:slope-lim:lin:bounds}
c_{ki}^\mathrm{min} \coloneqq \min_{\{ T_l \in \setT_h \,|\, \vec{x}_{ki} \in T_l \}} c_{l\mathrm{c}} \,,\quad
c_{ki}^\mathrm{max} \coloneqq \max_{\{ T_l \in \setT_h \,|\, \vec{x}_{ki} \in T_l \}} c_{l\mathrm{c}} 
\end{equation}
(cf.~Fig.~\ref{fig:vertex-lim} for an~illustration).
To enforce~\eqref{eq:slope-lim:lin:cond}, the correction factor~$\alpha_{ke}$ is defined as~\cite{Kuzmin2010}
\begin{equation}\label{eq:slope-lim:lin:corr}
\forall T_k\in\setT_h\,,\quad  \alpha_{ke} \coloneqq \min_{i\,\in\, \{1, 2, 3\}} \left\{
\begin{array}{llr}
  (c_{ki}^\mathrm{max} - c_{k\mathrm{c}})\big/(c_{ki} - c_{k\mathrm{c}})&\text{if}\;&c_{ki} > c_{ki}^\mathrm{max}\\
  1                  &\text{if}\;&c_{ki}^\mathrm{min}\le c_{ki}\le c_{ki}^\mathrm{max}\\
  (c_{ki}^\mathrm{min} - c_{k\mathrm{c}})\big/(c_{ki} - c_{k\mathrm{c}})&\text{if}\;&c_{ki} < c_{ki}^\mathrm{min}
\end{array}\right\}\;,
\end{equation}
where $c_{ki} \coloneqq c_{k\mathrm{c}} + \grad c_h (\vec{x}_{k\mathrm{c}}) \cdot (\vec{x}_{ki} - \vec{x}_{k\mathrm{c}})$ is the unconstrained linear reconstruction in $\vec{x}_{ki}$.
The limited counterpart of the DG~solution~\eqref{eq:taylor:representation} becomes then
\begin{equation*}
c_h(\vec{x}) = \overline{c}_h\,\phi_{k1} 
+ \alpha_{ke} \left[
  \left(\frac{\partial c_h}{\partial x^1}(\vec{x}_{k\mathrm{c}}) \, \Delta\left(x_k^1\right)\right)\,\phi_{k2}(\vec{x})
  + \left(\frac{\partial c_h}{\partial x^2}(\vec{x}_{k\mathrm{c}}) \, \Delta\left(x_k^2\right)\right)\,\phi_{k3}(\vec{x})
\right]
\qquad\text{on}~T_k\in\setT_h\,,
\end{equation*}
i.\,e., the linear degrees of freedom are scaled by $\alpha_{ke}$, and any degrees of freedom associated with higher polynomial degrees are set to zero. 
In the particular case of a~linear DG~approximation, the limiting can be performed using any hierarchical basis as opposed to higher-order DG~solutions that require a Taylor basis representation (see Sec. \ref{sec:slope-lim:kuzmin}).
In case of $\alpha_{ke}=1$, both linear and superlinear degrees of freedom remain unchanged.

\subsection{Hierarchical vertex-based limiter}\label{sec:slope-lim:kuzmin}
Further improvements by Kuzmin combine the vertex-based or standard Barth--Jespersen limiter with the higher order limiting scheme of Yang and Wang~\cite{YangWang2009}, who limit the numerical solution by multiplying all derivatives of order $q$ by a common correction factor $\alpha_{ke}^{(q)}$ instead of applying the correction only to the linear terms and dropping all higher degrees of freedom.
Kuzmin~\cite{Kuzmin2010} described this scheme in detail for quadratic representations, and here we offer a~closed form expression of this limiting procedure for DG~discretizations of arbitrary orders.
\par
Let $\setA_q\coloneqq\{\vec{a}\in\IN_0^2\,\big|\,|\vec{a}|=q\}$ be the set of all two-dimensional multi-indices of order $q$.
We determine the correction factor~$\alpha_{ke}^{(q)}$ for each order~$q \le p$ by computing correction factors~\eqref{eq:slope-lim:lin:corr} using the linear vertex-based limiter for all linear reconstructions of derivatives of order $q-1$,
\begin{equation}\label{eq:slope-lim:kuzmin:c_i}
\forall\,\vec{a} \in \setA_{q-1}\,, \quad  c_{k,\vec{a},i} \coloneqq C_{k,I(\vec{a})}^\mathrm{Taylor}\,\phi_{k1}(\vec{x}_{ki})
+ C_{k,I(\vec{a}+\transpose{[1,0]})}^\mathrm{Taylor}\,\phi_{k2}(\vec{x}_{ki})
+ C_{k,I(\vec{a}+\transpose{[0,1]})}^\mathrm{Taylor}\,\phi_{k3}(\vec{x}_{ki})
\quad\text{on}~T_k\in\setT_h\,,
\end{equation}
where the indices of the corresponding degrees of freedom are given by $I(\vec{a})$ and the first $x^1$-derivatives (identified by $I(\vec{a}+\transpose{[1,0]})$) and $x^2$-derivatives (identified by $I(\vec{a}+\transpose{[0,1]})$). 
Formally, the correction factor $\alpha_{ke}^{(q)}$ is defined as
\begin{equation}\label{eq:slope-lim:kuzmin:corr}
\alpha_{ke}^{(q)} = \min_{\vec{a}\in\setA_{q-1}} \, \alpha_{k\vec{a}}^{(q)} \,,
\qquad\text{with}\quad
\alpha_{k\vec{a}}^{(q)} \;\coloneqq\; \min_{i\in\{1,2,3\}} \left\{
\begin{array}{llr}
  (c_{k,\vec{a},i}^\mathrm{max} - c_{k,\vec{a},c})\big/(c_{k,\vec{a},i} - c_{k,\vec{a},c})&\text{if}\;&c_{\vec{a},i} > c_{\vec{a},i}^\mathrm{max}\\
  1                  &\text{if}\;&c_{k,\vec{a},i}^\mathrm{min}\le c_{k,\vec{a},i}\le c_{k,\vec{a},i}^\mathrm{max}\\
  (c_{k,\vec{a},i}^\mathrm{min} - c_{k,\vec{a},c})\big/(c_{k,\vec{a},i} - c_{k,\vec{a},c})&\text{if}\;&c_{k,\vec{a},i} < c_{k,\vec{a},i}^\mathrm{min}
\end{array}\right\}\;,
\end{equation}
where $c_{k,\vec{a},i}^\mathrm{min}, c_{k,\vec{a},i}^\mathrm{max}$ are defined as in~\eqref{eq:slope-lim:lin:bounds}.
To avoid the loss of accuracy at smooth extrema, the lower order derivatives should be limited by a factor not exceeding that of the higher order derivatives, since lower orders are typically smoother.
Beginning with the highest-order degrees of freedom, we compute the correction factors
\begin{equation}\label{eq:slope-lim:kuzmin:hierarch}
\forall q \ge 1\,\quad  \alpha_{ke}^{(q)} \;\coloneqq\; \max_{q \le d \le p} \alpha_{ke}^{(d)}  \,.
\end{equation}
Once the correction factor $\displaystyle \alpha_{ke}^{(q)}$ becomes equal to one for some $d>1$, no further limiting on this element is necessary.  The limited solution becomes
\begin{equation*}
c_h(\vec{x}) = \overline{c}_h\,\phi_{k1} 
+ \alpha_{ke}^{(1)} \left(\frac{\partial c_h}{\partial x^1}(\vec{x}_{k\mathrm{c}}) \, \Delta\left(x_k^1\right)\right)\,\phi_{k2}(\vec{x})
+ \alpha_{ke}^{(1)} \left(\frac{\partial c_h}{\partial x^2}(\vec{x}_{k\mathrm{c}}) \, \Delta\left(x_k^2\right)\right)\,\phi_{k3}(\vec{x})
+ \sum_{i=4}^{N} \alpha_{ke}^{(|\vec{a}_i|)}\left(\partial^{\vec{a}_i}c_h(\vec{x}_{k\mathrm{c}}) \, \left(\Delta\vec{x}_k\right)^{\vec{a}_i}\right) \, \phi_{ki}(\vec{x}) \,.
\end{equation*}

\subsection{Stricter form of the vertex-based limiter}\label{sec:slope-lim:strict}
Our numerical experiments showed that implicitly assuming that higher order derivatives are always smoother than lower order derivatives results in limiting procedures that do not guarantee strict fulfillment of condition~\eqref{eq:slope-lim:lin:cond}, especially at discontinuities in the solution.
We modified two key components of the limiter presented in the previous section and obtained a~limiter that exhibited slightly stronger peak clipping but turned out to be always effective:
\begin{enumerate}
\item 
Instead of employing only the linear reconstruction as given in Eq.~\eqref{eq:slope-lim:kuzmin:c_i}, we replace $c_{k,\vec{a},i}$ in the computation of the correction factor in Eq.~\eqref{eq:slope-lim:kuzmin:corr} by the full reconstruction
\begin{equation*}
\forall\,\vec{a}\in\setA_{q-1}\,,\quad 
c_{k,\vec{a},i} \;\coloneqq\; \sum_{0 \le |\vec{b}| < p - q} C
_{k,I(\vec{a}+\vec{b})}^\mathrm{Taylor} \, \phi_{k,I(\vec{b})}(\vec{x}_{ki})\,,
\end{equation*}
where $p$ is the polynomial degree of the DG~solution $c_h(\vec{x})$.
\item 
Again, we begin with the highest-order derivative but drop the hierarchical limiting condition~\eqref{eq:slope-lim:kuzmin:hierarch} and instead apply each correction coefficient $\alpha_{ke}^{(q)}$ immediately to all coefficients corresponding to polynomial degree~$q$ or higher.
These limited coefficients are then used to compute the next correction coefficient $\alpha_{ke}^{(q-1)}$.
\end{enumerate}
The result of our stricter limiter is
\begin{equation*}
c_h(\vec{x}) = \overline{c}_h\,\phi_{k1} 
+ \alpha_{ke}^{(1)} \left(\frac{\partial c_h}{\partial x^1}(\vec{x}_{k\mathrm{c}}) \, \Delta\left(x_k^1\right)\right)\,\phi_{k2}(\vec{x})
+ \alpha_{ke}^{(1)} \left(\frac{\partial c_h}{\partial x^2}(\vec{x}_{k\mathrm{c}}) \, \Delta\left(x_k^2\right)\right)\,\phi_{k3}(\vec{x})
+ \sum_{i=4}^{N} \left(\alpha_{ke}^{(1)} \cdots \alpha_{ke}^{(|\vec{a}_i|)}\right)
\Big(\partial^{\vec{a}_i}c_h (\vec{x}_{k\mathrm{c}}) \left(\Delta\vec{x}_k\right)^{\vec{a}_i}\Big) \, \phi_{ki}(\vec{x}) \,.
\end{equation*}

\subsection{Slope limiting in time-dependent problems}\label{sec:slope-lim:time}
For time-dependent problems, the slope limiting procedure is applied to each intermediate solution $\vec{C}^{(i)}$ in the update scheme~\eqref{eq:SSP-RK}.
However, due to the fact that the Taylor basis is non-orthogonal on triangles---as discussed in Sec.~\ref{sec:taylor:basis}---an~implicit coupling between the spatial derivatives is present and leads to non-smooth spatial variations in the time derivatives of $c_h$.
For that reason, Kuzmin~\cite{Kuzmin2012} applied the slope limiter not only to the solution of each Runge--Kutta stage but also to the time derivative $\dot{c} \coloneqq \partial_t c$ and used in addition a~filtering procedure that can be interpreted as \textit{selective mass lumping}.
\par
We describe this technique first for a~discretization in Taylor basis and extend it then to arbitrary basis representations. Let $\vecc{\Phi}^\mathrm{Taylor}$ denote the slope limiting operator that applies any of the above slope limiting procedures to a~global representation vector $\vec{C}^\mathrm{Taylor}(t)$ of a~solution $c_h(t)$ in Taylor basis representation.
The semi-discrete system~\eqref{eq:fullSystem} written in a~Taylor basis
\begin{equation*}
\vecc{M}_\mathrm{C} \partial_t \vec{C}^\mathrm{Taylor}(t) 
= \vec{S}^\mathrm{Taylor}\left(\vec{C}^\mathrm{Taylor}(t), t\right)
\end{equation*}
is replaced by
\begin{equation*}
\vecc{M}_\mathrm{L} \partial_t \vec{C}^\mathrm{Taylor}(t)
= \vec{S}^\mathrm{Taylor}\left(\vecc{\Phi}^\mathrm{Taylor} \vec{C}^\mathrm{Taylor}(t), t\right) 
+ \left( \vecc{M}_\mathrm{L} - \vecc{M}_\mathrm{C} \right) \, \vecc{\Phi}^\mathrm{Taylor} \partial_t \vec{C}^\mathrm{Taylor} (t)\,,
\end{equation*}
where $\vecc{M}_\mathrm{C} = \{m_{ij}\}$, $\vecc{M}_\mathrm{L} \coloneqq \mathrm{diag}\{m_{ii}\}$ denote the full and the lumped mass matrices in Taylor basis (to improve readability, we drop the superscript `Taylor' here).
Note that for the case $\vecc{\Phi}^\mathrm{Taylor}=\vecc{I}$ both formulations are identical.
Consequently, update scheme~\eqref{eq:SSP-RK} is modified replacing
\begin{equation*}
  \vec{C}^{\mathrm{Taylor},(i)} =\; \omega_i\,\vec{C}^{\mathrm{Taylor},n} + (1-\omega_i)\,\left( \vec{C}^{\mathrm{Taylor},(i-1)} 
+ \Delta t^n\,\dot{\vec{C}}^{\mathrm{Taylor},(i)} \right)
\quad \text{with} \quad
\dot{\vec{C}}^{\mathrm{Taylor},(i)} \coloneqq \vecc{M}_\mathrm{C}^{-1} \vec{S}^{\mathrm{Taylor},n + \delta_i}
\end{equation*}
by the selectively lumped and limited update 
\begin{equation*}
  \vec{C}^{\mathrm{Taylor},(i)} =\; \vecc{\Phi}^\mathrm{Taylor} \left[ \omega_i\,\vec{C}^{\mathrm{Taylor},n} + (1-\omega_i)\,\left( \vec{C}^{\mathrm{Taylor},(i-1)} 
+ \Delta t^n\, \tilde{\vec{C}}^{\mathrm{Taylor},(i)} \right) \right]
\end{equation*}
with
\begin{equation*}
\tilde{\vec{C}}^{\mathrm{Taylor},(i)} = \vecc{M}_\mathrm{L}^{-1} \left[ \left(\vecc{M}_\mathrm{L} - \vecc{M}_\mathrm{C}\right) \vecc{\Phi}^\mathrm{Taylor} \dot{\vec{C}}^{\mathrm{Taylor},(i)} + \vec{S}^{\mathrm{Taylor},n+\delta_i} \right] \,.
\end{equation*}
Although the mass matrix~$\vecc{M}$ in the modal basis is diagonal, the implicit coupling between the spatial $x^1$- and $x^2$-derivatives is still present.
Only the Taylor basis has vectors coinciding with the coordinate directions, hence, the lumping technique cannot be directly applied to representations in other bases.
\par
To get rid of this implicit coupling of the spatial derivatives in the time derivative of the modal~DG~basis, we reformulate the lumped time derivative in Taylor basis as
\begin{align*}
\tilde{\vec{C}}^{\mathrm{Taylor},(i)} 
=\;&\vecc{M}_\mathrm{L}^{-1} \vecc{M}_\mathrm{L} \vecc{\Phi}^\mathrm{Taylor} \dot{\vec{C}}^{\mathrm{Taylor},(i)}
-\vecc{M}_\mathrm{L}^{-1} \vecc{M}_\mathrm{C} \vecc{\Phi}^\mathrm{Taylor} \dot{\vec{C}}^{\mathrm{Taylor},(i)}
+\vecc{M}_\mathrm{L}^{-1} \underbrace{\vecc{M}_\mathrm{C}\vecc{M}_\mathrm{C}^{-1}}_{=\,\vecc{I}} \vec{S}^{\mathrm{Taylor},n+\delta_i} \\
=\;&\vecc{\Phi}^\mathrm{Taylor} \dot{\vec{C}}^{\mathrm{Taylor},(i)}
+\vecc{M}_\mathrm{L}^{-1} \vecc{M}_\mathrm{C} \Big(
  \underbrace{\vecc{M}_\mathrm{C}^{-1} \vec{S}^{\mathrm{Taylor},n+\delta_i}}_{=\dot{\vec{C}}^{\mathrm{Taylor},(i)}} 
  -\vecc{\Phi}^\mathrm{Taylor} \dot{\vec{C}}^{\mathrm{Taylor},(i)}
\Big) \\
=\;&\vecc{\Phi}^\mathrm{Taylor} \dot{\vec{C}}^{\mathrm{Taylor},(i)}
+\vecc{M}_\mathrm{L}^{-1} \vecc{M}_\mathrm{C} \Big(
  \dot{\vec{C}}^{\mathrm{Taylor},(i)}-\vecc{\Phi}^\mathrm{Taylor} \dot{\vec{C}}^{\mathrm{Taylor},(i)}
\Big) \,.
\end{align*}
Using the time derivative in modal basis $\dot{\vec{C}}^{(i)} = \vecc{M}^{-1} \vec{S}^{n+\delta_i}$ and transformation~\eqref{eq:taylor:trafo}, we obtain
\begin{equation}\label{eq:slope-lim:lumped-update}
\tilde{\vec{C}}^{(i)} 
=\; \vecc{M}^{-1} \vecc{M}^\mathrm{DG,Taylor} \tilde{\vec{C}}^{\mathrm{Taylor},(i)} 
=\; \vecc{M}^{-1} \vecc{M}^\mathrm{DG,Taylor} \left[
  \vecc{\Phi}^\mathrm{Taylor} \dot{\vec{C}}^{\mathrm{Taylor},(i)}
  +\vecc{M}_\mathrm{L}^{-1} \vecc{M}_\mathrm{C} \Big(
    \dot{\vec{C}}^{\mathrm{Taylor},(i)}-\vecc{\Phi}^\mathrm{Taylor}\dot{\vec{C}}^{\mathrm{Taylor},(i)}
  \Big)
\right]\,,
\end{equation} 
where
$\dot{\vec{C}}^{\mathrm{Taylor},(i)} \;\coloneqq\; 
\left(\vecc{M}^\mathrm{DG,Taylor}\right)^{-1} \vecc{M} \dot{\vec{C}}^{(i)}$.
Thus, the fully modified version of update scheme~\eqref{eq:SSP-RK} reads as
\begin{equation}\label{eq:slope-lim:SSP-RK}
\begin{array}{lll}
\vec{C}^{(0)} &=\; \vec{C}^{n} \,, \\
\vec{C}^{(i)} &=\; \vecc{\Phi} \left[ \omega_i\,\vec{C}^{n} + (1-\omega_i)\,\left( \vec{C}^{(i-1)} 
+ \Delta t^n\,\tilde{\vec{C}}^{(i)}  \right) \right] \quad\text{for}~ \forall i = 1,\ldots,s\,,\\
\vec{C}^{n+1} &=\; \vec{C}^{(s)} \,,
\end{array}
\end{equation}
with $\tilde{\vec{C}}^{(i)} $ as given in~\eqref{eq:slope-lim:lumped-update} and the slope limiting operator $\vecc{\Phi}$ formally defined as
\begin{equation}\label{eq:slope-lim:op}
\vecc{\Phi} \coloneqq \vecc{M}^{-1} \vecc{M}^\mathrm{DG,Taylor} \vecc{\Phi}^\mathrm{Taylor} \left(\vecc{M}^\mathrm{DG,Taylor}\right)^{-1} \vecc{M}\,.
\end{equation}

\subsection{Boundary conditions}
Problems can occur when computing the bounds $c_{ki}^\mathrm{min}, c_{ki}^\mathrm{max}$ in Eq.~\eqref{eq:slope-lim:lin:cond} for control points $\vec{x}_{ki} \in \partial\Omega_\mathrm{in}$ on the Dirichlet boundary.
To account for the boundary data in the limiting procedure we include on those control points the boundary value~$c_\mathrm{D}(t^n + \delta_i \Delta t^n, \vec{x}_{ki})$ in the minimum/maximum-operation in Eq.~\eqref{eq:slope-lim:lin:bounds} for $\alpha_{ke}^{(1)}$ when applying the slope limiting operator~$\vecc{\Phi}$ in Eq.~\eqref{eq:slope-lim:SSP-RK}.

\section{Implementation}\label{sec:implementation}
An~extensive documentation on our data structures, grid, etc. can be found in the first paper of the series~\cite{FESTUNG1}.
These explanations are not reproduced here; greater detail is provided for routines first introduced in the present work.

\subsection{Backtransformation to the reference triangle}\label{sec:transformationtoThat}
\begin{figure}[t!]
\centering%
\includegraphics{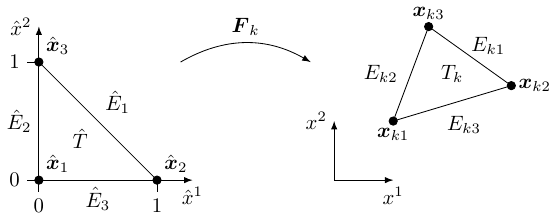}
\caption{The affine mapping~$\vec{F}_k$ transforms the reference triangle~$\hat{T}$ with vertices~$\hat{\vec{x}}_1 = \transpose{[0,\,0]}$, $\hat{\vec{x}}_2 = \transpose{[1,\,0]}$, $\hat{\vec{x}}_3 = \transpose{[0,\,1]}$ to the physical triangle~$T_k$ with counter-clockwise-ordered vertices~$\vec{x}_{ki}$, $i\in\{1,2,3\}$.}%
\label{fig:referencetriangle}
\end{figure}

We are using a back transformation to the reference triangle $\hat{T} = \{(0,0),(1,0),(0,1)\}$ defined by an affine mapping
\begin{equation}\label{eq:affinemappings}
\vec{F}_k :\quad  \hat{T} \ni \hat{\vec{x}}\mapsto \vecc{B}_k \hat{\vec{x}} + \vec{x}_{k1} = \vec{x}\in T_k\,,
\qquad\text{with}\quad
\IR^{2\times2} \ni \vecc{B}_k \coloneqq 
\left[ \vec{x}_{k2}-\vec{x}_{k1} \,\big|\, \vec{x}_{k3}-\vec{x}_{k1} \right]\,,
\end{equation}
for any triangle $T_k = \{\vec{x}_{k1}, \vec{x}_{k2}, \vec{x}_{k3}\} \in \setT_h$ (see Fig.~\ref{fig:referencetriangle}).
It holds $0 < \det \vecc{B}_k = 2|T_k|$.
Any function~$\,w:T_k\rightarrow \IR\,$ implies $\hat{w}:\hat{T}\rightarrow \IR\,$ by $\,\hat{w}=w\circ \vec{F}_k\,$, i.\,e., $\,w(\vec{x}) = \hat{w}(\hat{\vec{x}})\,$.  The transformation of the gradient is obtained by the chain rule:
\begin{equation}\label{eq:affinegradtrafo}
  \grad \;=\; \invtrans{\big(\hat{\grad}\vec{F}_k \big)}\,\hat{\grad}\,,
\end{equation}
where we abbreviated~$\hat{\grad} = \transpose{[\partial_{\hat{x}^1},\partial_{\hat{x}^2} ]}$.
This results in transformation formulas for integrals over an~element~$T_k$ or an~edge~$E_{kn} \subset T_k$ for a~function $w: \Omega \rightarrow \IR$
\begin{subequations}
\begin{align}
\label{eq:trafoRule:T}
\int_{T_k}w(\vec{x})\,\dd\vec{x}
&\;=\;\frac{\abs{T_k}}{\abs{\hat{T}}} \int_{\hat{T}}w\circ\vec{F}_k(\hat{\vec{x}})\,\dd\hat{\vec{x}}
\;=\;2\abs{T_k} \int_{\hat{T}}w\circ\vec{F}_k(\hat{\vec{x}})\,\dd\hat{\vec{x}}
\;=\;2\abs{T_k} \int_{\hat{T}}\hat{w}(\hat{\vec{x}})\,\dd\hat{\vec{x}}
\;,\\
\label{eq:trafoRule:E}
\int_{E_{kn}} w(\vec{x})\,\dd\vec{x} 
&\;=\; \frac{\abs{E_{kn}}}{\abs{\hat{E}_n}} \int_{\hat{E}_n} w\circ\vec{F}_k(\hat{\vec{x}})\,\dd\hat{\vec{x}}
\;=\;\frac{\abs{E_{kn}}}{\abs{\hat{E}_n}} \int_{\hat{E}_n} \hat{w}(\hat{\vec{x}})\,\dd\hat{\vec{x}}
\;.
\end{align}
\end{subequations}

\subsection{Numerical integration}\label{sec:quadrature}
As an~alternative to the symbolic integration functions provided by \Matlab, we implemented a~quadrature
integration functionality for triangle and edge integrals. 
\par
Since we transform all integrals on~$T_k\in\setT_h$ to the reference triangle~$\hat{T}$ (cf.~Sec.~\ref{sec:transformationtoThat}), it is sufficient to define the quadrature rules on~$\hat{T}$ (which, of course, can be rewritten to apply for every physical triangle~$T = \vec{F}_T(\hat{T})$):
\begin{equation}\label{eq:quadrature}
\int_{\hat{T}} \hat{g}(\hat{\vec{x}})\,\dd\hat{\vec{x}}\;\approx\; \sum_{r=1}^R \omega_r\, \hat{g}(\hat{\vec{q}}_r)
\end{equation}
with $R$~\emph{quadrature points}~$\hat{\vec{q}}_r\in\hat{T}$ and \emph{quadrature weights}~$\omega_r\in\IR$.  The \emph{order} of a~quadrature rule is the largest integer~$s$ such that~\eqref{eq:quadrature} is \emph{exact} for polynomials~$g\in\IP_s(\hat{T})$.   
Note that we exclusively rely on quadrature rules with positive weights and quadrature points located strictly in the interior of~$\hat{T}$ and not on $\partial{\hat{T}}$. 
The rules used in the implementation are found in the routine~\code{quadRule2D}. 
An~overview of quadrature rules on triangles can be found in the \enquote{Encyclopaedia of Cubature Formulas}~\cite{Cools2003}. 
For edge integration, we rely on standard Gauss quadrature rules of required order.
\par
The integrals in~\eqref{eq:spacediscretesystem} contain integrands that are polynomials of maximum order~$3p-1$ on triangles and of maximum order~$3p$ on edges. 
Using quadrature integration, one could choose rules that integrate all such terms exactly; however, sufficient accuracy can be achieved with quadrature rules that are exact for polynomials of order~$2p$ on triangles and $2p+1$ on edges (cf.~\cite{CockburnShu1998b}).

\subsection{Assembly}
The aim of this section is to transform the terms required to build the block matrices in~\eqref{eq:timeDepSystem} to the reference triangle~$\hat{T}$ and then to evaluate those either via numerical quadrature or analytically.  
The assembly of block matrices from local contributions is then performed in vectorized operations.
\par
For the implementation, we need the explicit form for the components of the mappings~$\vec{F}_k:\hat{T}\rightarrow T_k$ and their inverses~$\vec{F}_k^{-1}: T_k\rightarrow \hat{T}$ as defined in~\eqref{eq:affinemappings}.
Recalling that $0<\det  \vecc{B}_k = 2\abs{T_k}$ (cf.~Sec.~\ref{sec:transformationtoThat}) we obtain
\begin{equation*}
\vec{F}_k(\hat{\vec{x}}) \;=\;
\begin{bmatrix}
B_k^{11}\,\hat{x}^1 + B_k^{12}\,\hat{x}^2 + a_{k1}^1\\
B_k^{21}\,\hat{x}^1 + B_k^{22}\,\hat{x}^2 + a_{k1}^2
\end{bmatrix}
\qquad\text{and}\qquad
\vec{F}_k^{-1}(\vec{x}) \;=\;
\frac{1}{2\,\abs{T_k}}
\begin{bmatrix}
B_k^{22}\,(x^1 - a_{k1}^1) - B_k^{12}\,(x^2 - a_{k1}^2)\\
B_k^{11}\,(x^2 - a_{k1}^2) - B_k^{21}\,(x^1 - a_{k1}^1)
\end{bmatrix}\;.
\end{equation*}
From \eqref{eq:affinegradtrafo} we obtain the component-wise rule for the gradient in~$\vec{x}\in T_k$:
\begin{equation}\label{eq:rule:gradient}
\begin{bmatrix}
\partial_{x^1}\\
\partial_{x^2}
\end{bmatrix}
\;=\;
\frac{1}{2\,\abs{T_k}}
\begin{bmatrix}
B_k^{22}\,\partial_{\hat{x}^1} - B_k^{21}\,\partial_{\hat{x}^2}\\
B_k^{11}\,\partial_{\hat{x}^2} - B_k^{12}\,\partial_{\hat{x}^1}
\end{bmatrix}.
\end{equation}
Similarly to the first paper in series~\cite{FESTUNG1}, we make extensive use of the Kronecker product $\vecc{A} \otimes \vecc{B}$ of two matrices $\vecc{A} = [a_{ij}] \in \IR^{m_a\times n_a}$, $\vecc{B} = [b_{kl}] \in \IR^{m_b\times n_b}$ defined as
\begin{equation}\label{eq:kron}
\vecc{A} \otimes \vecc{B} \coloneqq \left[ a_{ij} \vecc{B}\right] \in \IR^{m_a m_b\times n_a n_b}\,.
\end{equation}
\par
In the following, we present the necessary transformation for all blocks of system~\eqref{eq:timeDepSystem} and name the corresponding \MatOct~routines that can be found in Sec.~\ref{sec:routines}.

\subsubsection[Assembly of M]{Assembly of $\vecc{M}$}\label{sec:assembly:globM}
Using the transformation rule~\eqref{eq:trafoRule:T}, the following holds for the local mass matrix~$\vecc{M}_{T_k}$ as defined in~\eqref{eq:globMlocM}:
\begin{equation}\label{eq:hatM}
\vecc{M}_{T_k} \;=\; 2\abs{T_k}\,\hat{\vecc{M}}
\qquad\text{with}\qquad
\hat{\vecc{M}}\;\coloneqq\;
\int_{\hat{T}}\,\begin{bmatrix}
\hat{\vphi}_{1}\,\hat{\vphi}_{1} & \cdots & \hat{\vphi}_{1}\,\hat{\vphi}_{N} ~\\
 \vdots                & \ddots & \vdots \\
\hat{\vphi}_{N}\,\hat{\vphi}_{1} & \cdots & \hat{\vphi}_{N}\,\hat{\vphi}_{N} 
\end{bmatrix}\;,
\end{equation}
where~$\hat{\vecc{M}}\in\IR^{N\times N}$ is the representation of the local mass matrix on the reference triangle~$\hat{T}$.  
With~\eqref{eq:globMlocM} we see that the global mass matrix~$\vecc{M}$ can be expressed as a~Kronecker product of a~matrix containing the areas~$\abs{T_k}$ and the local matrix~$\hat{\vecc{M}}$:
\begin{equation*}
\vecc{M}
\;=\;
\begin{bmatrix}
\vecc{M}_{T_1} &          & \\
               & ~\ddots~ & \\
               &          & \vecc{M}_{T_K}
\end{bmatrix}
\;=\;
2 \begin{bmatrix}
\abs{T_1} &          & \\
               & ~\ddots~ & \\
               &          & \abs{T_K}
\end{bmatrix} \otimes \hat{\vecc{M}}\;.
\end{equation*}
In the corresponding assembly routine~\code{assembleMatElemPhiPhi}, the sparse block-diagonal matrix is generated using the command~\code{spdiags} with the list~\code{g.areaT} (cf.~\cite{FESTUNG1}).

\subsubsection[Assembly of Gm]{Assembly of $\vecc{G}^m$}\label{sec:assembly:globG}
Application of the product rule, \eqref{eq:trafoRule:T}, and \eqref{eq:rule:gradient} gives us
\begin{equation*}
\int_{T_k}  \partial_{x^1}\vphi_{ki}\,\vphi_{kl}\,\vphi_{kj} \;=\;
\phantom{-}B_k^{22}\,[\hat{\vecc{G}}]_{i,j,l,1} - B_k^{21}\,[\hat{\vecc{G}}]_{i,j,l,2}\;,
\qquad
\int_{T_k}  \partial_{x^2}\vphi_{ki}\,\vphi_{kl}\,\vphi_{kj} \;=\;
-B_k^{12}\,[\hat{\vecc{G}}]_{i,j,l,1} + B_k^{11}\,[\hat{\vecc{G}}]_{i,j,l,2}
\end{equation*}
with a~multidimensional array~$\hat{\vecc{G}}\in\IR^{N\times N\times N\times 2}$ representing the transformed integral on the reference triangle~$\hat{T}$:
\begin{equation}\label{eq:hatG}
[\hat{\vecc{G}}]_{i,j,l,m}\;\coloneqq\;\int_{\hat{T}} \partial_{\hat{x}^m} \hat{\vphi}_i\, \hat{\vphi}_j\, \hat{\vphi}_l\, , \quad \mbox{for }m\in\{1,2\}\,.
\end{equation}
Now we can express the local matrix~$\vecc{G}^1_{T_k}$ from~\eqref{eq:locGm} as
\begin{align*}
\vecc{G}^1_{T_k} &= \sum_{l=1}^N U^1_{kl}(t)
\left(B_k^{22}\int_{\hat{T}}\begin{bmatrix}
\partial_{\hat{x}^1}\hat{\vphi}_1\hat{\vphi}_1\hat{\vphi}_l & \cdots & \partial_{\hat{x}^1}\hat{\vphi}_1\hat{\vphi}_N\hat{\vphi}_l \\
 \vdots                & \ddots & \vdots \\
 \partial_{\hat{x}^1}\hat{\vphi}_N\hat{\vphi}_1\hat{\vphi}_l & \cdots & \partial_{\hat{x}^1}\hat{\vphi}_N\hat{\vphi}_N\hat{\vphi}_l
\end{bmatrix}
-B_k^{21}\int_{\hat{T}}\begin{bmatrix}
\partial_{\hat{x}^2}\hat{\vphi}_1\hat{\vphi}_1\hat{\vphi}_l & \cdots & \partial_{\hat{x}^2}\hat{\vphi}_1\hat{\vphi}_N\hat{\vphi}_l \\
 \vdots                & \ddots & \vdots \\
 \partial_{\hat{x}^2}\hat{\vphi}_N\hat{\vphi}_1\hat{\vphi}_l & \cdots & \partial_{\hat{x}^2}\hat{\vphi}_N\hat{\vphi}_N\hat{\vphi}_l
\end{bmatrix}
\right)
\\
&= \sum_{l=1}^N U^1_{kl}(t)
\left(B_k^{22} [\hat{\vecc{G}}]_{:,:,l,1}-B_k^{21}  [\hat{\vecc{G}}]_{:,:,l,2}\right)
\end{align*}
and analogously~$\vecc{G}^2_{T_k}$.  
With $\vecc{G}^m=\diag(\vecc{G}_{T_1}^m,\ldots,\vecc{G}_{T_K}^m)$ we can vectorize over all triangles using the Kronecker product as done in the routine \code{assembleMatElemDphiPhiFuncDiscVec}.
We would like to point out that this is identical to the assembly of~$\vecc{G}^m$ in our first paper~\cite{FESTUNG1} except for the vectorial coefficient function; however, the respective section had some typos, above is the corrected version.
Additionally, in the assembly routine the wrong component of the normal vector was used for the assembly of $\vecc{G}^m$.
Corrected and up-to-date versions of the code can be found in our Github-repository~\cite{FESTUNGGithub}.

\subsubsection[Assembly of R]{Assembly of~$\vecc{R}$}\label{sec:assembly:globR}
To ease the assembly of~$\vecc{R}$ we split the global matrix as given in~\eqref{eq:globR} into a~block-diagonal part and a~remainder so that~$\vecc{R} = \vecc{R}^{\mathrm{diag}} + \vecc{R}^{\mathrm{offdiag}}$ holds.
We first consider the block-diagonal entries of $\vecc{R}$ given in Eqns.~\eqref{eq:globRInterior:diag},~\eqref{eq:globRBoundary:diag} and transform the integral terms to the $n$-th edge of the reference triangle $\hat{E}_n$:
\begin{align}
\int_{E_{kn}} \vphi_{ki}\,\vphi_{kj}\,
\Big( \vec{u}\cdot\vec{\nu}_{kn} \Big) \,
\delta_{\vec{u}\cdot\vec{\nu}_{kn}\ge0}\,\dd\vec{x}
&= \frac{|E_{kn}|}{|\hat{E}_n|} \int_{\hat{E}_n} \hat{\vphi}_i\,\hat{\vphi}_j\,
\Big( \left( \vec{u}\circ \vec{F}_k(\hat{\vec{x}}) \right) \cdot \vec{\nu}_{kn} \Big)\,
\delta_{\vec{u}\cdot\vec{\nu}_{kn}\ge0} \, \dd \hat{\vec{x}}\nonumber\\
&= \frac{|E_{kn}|}{|\hat{E}_n|} \int_0^1
\hat{\vphi}_i\circ\hat{\vec{\gamma}}_n(s)\,\hat{\vphi}_j\circ\hat{\vec{\gamma}}_n(s)\,
\Big( \underbrace{\left( \vec{u}\circ\vec{F}_k\circ\hat{\vec{\gamma}}_n(s) \right)}_{\eqqcolon\; \hat{\vec{u}}_{kn}(s)} \cdot \vec{\nu}_{kn} \Big)\,
\delta_{\hat{\vec{u}}_{kn}\cdot\vec{\nu}_{kn}\ge0} \,|\hat{\vec{\gamma}}_n'(s)| \,\dd s\nonumber\\
&\approx |E_{kn}| \sum_{r=1}^R \omega_r
\underbrace{\hat{\vphi}_i\circ\hat{\vec{\gamma}}_n(q_r)\,\hat{\vphi}_j\circ\hat{\vec{\gamma}}_n(q_r)\,}_{\eqqcolon \left[ \hat{\vecc{R}}^\mathrm{diag} \right]_{i,j,n,r}}
\Big( \hat{\vec{u}_{kn}(q_r)} \cdot \vec{\nu}_{kn} \Big)\,
\delta_{\hat{\vec{u}}_{kn}\cdot\vec{\nu}_{kn}\ge0} \,,
\label{eq:hatRdiag}
\end{align}
where we used transformation rule~\eqref{eq:trafoRule:E}, quadrature rule~\eqref{eq:quadrature}, and $\abs{\hat{\vec{\gamma}}_n'(s)}=\abs{\hat{E}_n}$. 
The explicit forms of the mappings~$\hat{\vec{\gamma}}_n:[0,1]\rightarrow \hat{E}_n$ can be easily derived:
\begin{equation}\label{eq:gammaMap}
\hat{\vec{\gamma}}_1(s) \;\coloneqq\;
\begin{bmatrix}
1-s\\s
\end{bmatrix}\,,\qquad
\hat{\vec{\gamma}}_2(s) \;\coloneqq\;
\begin{bmatrix}
0\\1-s
\end{bmatrix}\,,\qquad
\hat{\vec{\gamma}}_3(s) \;\coloneqq\;
\begin{bmatrix}
s\\0
\end{bmatrix}\,.
\end{equation}
Thus, we can assemble the global matrix using the Kronecker product
\begin{equation*}
\vecc{R}^{\mathrm{diag}} \coloneqq \sum_{n=1}^3 \sum_{r=1}^R \omega_r \,
\begin{bmatrix}
|E_{1n}| & & \\
& \ddots & \\
& & |E_{Kn}|
\end{bmatrix}
\circ
\begin{bmatrix}
\left(\hat{\vec{u}}_{1n}(q_r)\cdot\vec{\nu}_{1n}\right)\delta_{\hat{\vec{u}}_{1n}\cdot\vec{\nu}_{1n}\ge0} & & \\
& \ddots & \\
& & \left(\hat{\vec{u}}_{Kn}(q_r)\cdot\vec{\nu}_{Kn}\right)\delta_{\hat{\vec{u}}_{Kn}\cdot\vec{\nu}_{Kn}\ge0}
\end{bmatrix}
\otimes
\left[\hat{\vecc{R}}^\mathrm{diag} \right]_{:,:,n,r} \,,
\end{equation*}
where `$\circ$' denotes the Hadamard product.
\par
Next, we consider the off-diagonal blocks of~$\vecc{R}$ stored in~$\vecc{R}^\mathrm{offdiag}$.  
For an~interior edge~ $E_{k^-n^-}=E_{k^+n^+}\in\partial T_{k^-}\cap \partial T_{k^+}, \; n^-,n^+\in\{1,2,3\}$, we obtain analogously:
\begin{align*}
&\int_{E_{k^-n^-}} \vphi_{k^-i}\,\vphi_{k^+j}\,
\Big( \vec{u}\cdot\vec{\nu}_{k^-n} \Big) \,
\delta_{\vec{u}\cdot\vec{\nu}_{k^-n}<0}\,\dd\vec{x}\\
&\qquad= \frac{|E_{k^-n^-}|}{|\hat{E}_{n^-}|} \int_{\hat{E}_{n^-}} \hat{\vphi}_i\,
\vphi_{k^+j} \circ \overbrace{\vec{F}_{k^+} \circ \vec{F}_{k^+}^{-1}}^{=\,\vecc{I}} \circ \vec{F}_{k^-}(\hat{\vec{x}})\,
\Big( \left( \vec{u}\circ \vec{F}_{k^-}(\hat{\vec{x}}) \right) \cdot \vec{\nu}_{k^-n^-} \Big)\,
\delta_{\vec{u}\cdot\vec{\nu}_{k^-n^-}<0} \, \dd \hat{\vec{x}}\\
&\qquad= \frac{|E_{k^-n^-}|}{|\hat{E}_{n^-}|} \int_{\hat{E}_{n^-}} \hat{\vphi}_i\,
\hat{\vphi}_j \circ \vec{F}_{k^+}^{-1} \circ \vec{F}_{k^-}(\hat{\vec{x}})\,
\Big( \left( \vec{u}\circ \vec{F}_{k^-}(\hat{\vec{x}}) \right) \cdot \vec{\nu}_{k^-n^-} \Big)\,
\delta_{\vec{u}\cdot\vec{\nu}_{k^-n^-}<0} \, \dd \hat{\vec{x}}\\
&\qquad= |E_{k^-n^-}| \int_0^1
\hat{\vphi}_i\circ\hat{\vec{\gamma}}_{n^-}(s)\,
\hat{\vphi}_j \circ \vec{F}_{k^+}^{-1} \circ \vec{F}_{k^-} \circ \hat{\vec{\gamma}}_{n^-}(s)\,
\Big( \hat{\vec{u}}_{k^-n^-}(s) \cdot \vec{\nu}_{k^-n^-} \Big)\,
\delta_{\hat{\vec{u}}_{k^-n^-}\cdot\vec{\nu}_{k^-n^-}<0}  \,\dd s\;.
\end{align*}
Note that~$\vec{F}_{k^+}^{-1}\circ\vec{F}_{k^-}$ maps from~$\hat{T}$ to~$\hat{T}$.  
Since we compute a~line integral, the integration domain is further restricted to an edge~$\hat{E}_{n^-}$, $n^-\in\{1,2,3\}$ and its co-domain to an~edge~$\hat{E}_{n^+}$, $n^+\in\{1,2,3\}$.
As a~result, this integration can be boiled down to nine possible maps between the sides of the reference triangle expressed as
\begin{equation*}
\mapEE_{n^-n^+}:\quad \hat{E}_{n^-}\ni\hat{\vec{x}}\mapsto \mapEE_{n^-n^+}(\hat{\vec{x}})\,=\,\vec{F}_{k^+}^{-1}\circ\vec{F}_{k^-}(\hat{\vec{x}})\in \hat{E}_{n^+}
\end{equation*}
for an~arbitrary index pair~$\{k^-,k^+\}$ as described above. The closed-form expressions of the nine cases are given in our first paper~\cite{FESTUNG1}.
We apply quadrature rule~\eqref{eq:quadrature} and define $\hat{\vecc{R}}^\mathrm{offdiag} \in \IR^{N\times N\times 3 \times 3\times R}$ by
\begin{equation}\label{eq:hatRoffdiag}
\left[\hat{\vecc{R}}^\mathrm{offdiag}\right]_{i,j,n^-,n^+,r} \;\coloneqq\;
\hat{\vphi}_i\circ\hat{\vec{\gamma}}_{n^-}(q_r)\,\hat{\vphi}_j\circ\mapEE_{n^-n^+}\circ\hat{\vec{\gamma}}_{n^-}(q_r)
\end{equation}
and thus arrive at
\begin{multline*}
\vecc{R}^\mathrm{offdiag} \coloneqq \sum_{n^-=1}^3\sum_{n^+=1}^3\sum_{r=1}^R
\omega_r\,
\begin{bmatrix}
0                          & \delta_{E_{1n^-}=E_{2n^+}} & \hdots                & \hdots                & \delta_{E_{1n^-}=E_{Kn^+}} \\
\delta_{E_{2n^-}=E_{1n^+}} & 0                          &   \ddots              &                       & \textstyle\vdots  \\ 
\vdots                     &         \ddots             & \ddots                &          \ddots       & \vdots   \\
\vdots                     & {}                         &    \ddots             & 0                     & \delta_{E_{(K-1)n^-}=E_{Kn^+}} \\
\delta_{E_{Kn^-}=E_{1n^+}}  &  \hdots                   & \hdots                & \delta_{E_{Kn^-}=E_{(K-1)n^+}}  & 0
\end{bmatrix}
\\
\circ
\begin{bmatrix}
|E_{1n^-}|
\left(\hat{\vec{u}}_{1n^-}\cdot\vec{\nu}_{1n^-}\right)\delta_{\hat{\vec{u}}_{1n^-}\cdot\vec{\nu}_{1n^-}<0} &
\cdots &
|E_{1n^-}|
\left(\hat{\vec{u}}_{1n^-}\cdot\vec{\nu}_{1n^-}\right)\delta_{\hat{\vec{u}}_{1n^-}\cdot\vec{\nu}_{1n^-}<0} \\
\vdots & & \vdots \\
|E_{Kn^-}|
\left(\hat{\vec{u}}_{Kn^-}\cdot\vec{\nu}_{Kn^-}\right)\delta_{\hat{\vec{u}}_{Kn^-}\cdot\vec{\nu}_{Kn^-}<0} &
\cdots &
|E_{Kn^-}|
\left(\hat{\vec{u}}_{Kn^-}\cdot\vec{\nu}_{Kn^-}\right)\delta_{\hat{\vec{u}}_{Kn^-}\cdot\vec{\nu}_{Kn^-}<0}
\end{bmatrix}
\otimes [\hat{\vecc{R}}^\mathrm{offdiag}]_{:,:,n^-,n^+,r}\;.
\end{multline*}
The sparsity structure for off-diagonal blocks depends on the numbering of mesh entities and is given for each combination of $n^-$ and $n^+$ by the list~\code{markE0TE0T}.
Due to the upwind flux in the edge integrals, it is not possible here to include the quadrature rule directly in the element-blocks~$\hat{\vecc{R}}^\mathrm{offdiag}$ as opposed to the assembly of element integrals and the assembly routines for edges in the first paper~\cite{FESTUNG1}.
Implementing the assembly with above formulation is possible but expensive, since the global matrix would have to to be built for every quadrature point.
Instead we make use of the fact that the sparsity structure is the same for every quadrature point as we solely rely on the element numbering and do not account for the upwind direction when determining the structure of~$\vecc{R}^\mathrm{offdiag}$.
We define a~tensor of block vectors~$\tilde{\vecc{R}}^\mathrm{offdiag} \in \IR^{KN\times N\times 3\times 3}$
\begin{equation*}
\left[ \tilde{\vecc{R}}^\mathrm{offdiag} \right]_{(k-1)N+1:kN,:,n-,n+} \coloneqq\;
\sum_{r=1}^R \omega_r 
\begin{bmatrix}
|E_{1n^-}|\left(\hat{\vec{u}}_{1n^-}\cdot\vec{\nu}_{1n^-}\right)\delta_{\hat{\vec{u}}_{1n^-}\cdot\vec{\nu}_{1n^-}<0} \\
\vdots \\
|E_{Kn^-}|\left(\hat{\vec{u}}_{Kn^-}\cdot\vec{\nu}_{Kn^-}\right)\delta_{\hat{\vec{u}}_{Kn^-}\cdot\vec{\nu}_{Kn^-}<0}
\end{bmatrix}
\otimes [\hat{\vecc{R}}^\mathrm{offdiag}]_{:,:,n^-,n^+,r}\;.
\end{equation*}
For $m_b = r m_a, r\in\IN$, let
\begin{equation}\label{eq:kronVec}
\cdot\,\otimes_\mathrm{V}\,\cdot\,:\quad  \IR^{m_a \times n_a}\times \IR^{m_b \times n_b}\ni (\vecc{A},\vecc{B})\mapsto\vecc{A} \otimes_\mathrm{V}\vecc{B}\,\coloneqq\,\left[[\vecc{A}]_{i,j} \left[\vecc{B}\right]_{(i-1)r\,:\,ir, :} \right]\in \IR^{m_b \times n_a n_b}
\end{equation}
be an~operator which can be interpreted as a~Kronecker product (cf.~Eq.~\eqref{eq:kron}) that takes a~different right-hand side for every row of the left-hand side and is implemented in the routine \code{kronVec}.
This allows us to write
\begin{equation*}
\vecc{R}^\mathrm{offdiag} = \sum_{n^-=1}^3\sum_{n^+=1}^3
\begin{bmatrix}
0                          & \delta_{E_{1n^-}=E_{2n^+}} & \hdots                & \hdots                & \delta_{E_{1n^-}=E_{Kn^+}} \\
\delta_{E_{2n^-}=E_{1n^+}} & 0                          &   \ddots              &                       & \textstyle\vdots  \\ 
\vdots                     &         \ddots             & \ddots                &          \ddots       & \vdots   \\
\vdots                     & {}                         &    \ddots             & 0                     & \delta_{E_{(K-1)n^-}=E_{Kn^+}} \\
\delta_{E_{Kn^-}=E_{1n^+}}  &  \hdots                   & \hdots                & \delta_{E_{Kn^-}=E_{(K-1)n^+}}  & 0
\end{bmatrix}
\otimes_\mathrm{V} \left[ \tilde{\vecc{R}}^\mathrm{offdiag} \right]_{:,:,n-,n+}\,,
\end{equation*}
which omits the expensive assembly of the global matrix in every quadrature point.
\par
The routine~\code{assembleMatEdgePhiPhiValUpwind} assembles the matrices $\vecc{R}^\mathrm{diag}$ and $\vecc{R}^\mathrm{offdiag}$ directly into $\vecc{R}$ with a~code very similar to the formulation above.
To avoid repeated computation of the normal velocity $\vec{u}\cdot\vec{\nu}_{kn}$, we evaluate it once in all quadrature points on each edge using the globally continuous function $\vec{u}(t,\vec{x})$ and store it in a~dedicated variable~\code{vNormalOnQuadEdge}, which is then employed in the decision of the upwind direction using $\delta_{\vec{u}\cdot\vec{\nu}_{kn}\ge0}$ and $\delta_{\vec{u}\cdot\vec{\nu}_{kn}<0}$ as well as the normal velocity in the assembly of~$\vecc{R}$ and~$\vec{K}_\mathrm{D}$.

\subsubsection[Assembly of KD]{Assembly of~$\vec{K}_{\mathrm{D}}$}\label{sec:assembly:globKD}
The entries of $\vec{K}_{\mathrm{D}}$ in~\eqref{eq:globKD} are transformed using transformation rule~\eqref{eq:trafoRule:E} and mapping~\eqref{eq:gammaMap}
\begin{align*}
\left[\vec{K}_\mathrm{D}\right]_{(k-1)N+i} 
&\;=\;\sum_{E_{kn}\in\partial T_k\cap\setE_{\partial\Omega}} \int_{E_{kn}}
\vphi_{ki}\,c_{\mathrm{D}}(t) 
\left(\vec{u} \cdot \nu_{kn}^m\right) \delta_{\vec{u}\cdot\nu_{kn}^m<0} \,\dd\vec{x} \\
&\;=\; \sum_{E_{kn}\in\partial T_k\cap\setE_{\partial\Omega}} \frac{|E_{kn}|}{|\hat{E}_n|} \int_{\hat{E}_n}
\hat{\vphi}_i\,c_{\mathrm{D}}(t, \vec{F}_k(\hat{\vec{x}})) 
\left( \left(\vec{u} \circ \vec{F}_k(\hat{\vec{x}}) \right) \cdot \nu_{kn}^m\right)
\delta_{\vec{u}\cdot\nu_{kn}^m<0}\, \dd\hat{\vec{x}} \\
&\;=\; \sum_{E_{kn}\in\partial T_k\cap\setE_{\partial\Omega}}|E_{kn}| \int_0^1
\hat{\vphi}_i\circ\hat{\vec{\gamma}}_n(s) \,c_{\mathrm{D}}(t, \vec{F}_k\circ\hat{\vec{\gamma}}_n(s)) 
\left( \hat{\vec{u}}_{kn}(s) \cdot \nu_{kn}^m\right)
\delta_{\hat{\vec{u}}_{kn}\cdot\nu_{kn}^m<0}\, \dd s \,,
\end{align*}
where we again implicitly assumed the application of $\vec{F}_k$ and $\hat{\vec{\gamma}}_n$ to $\delta_{\vec{u}\cdot\vec{\nu}_{kn}\ge0}$.
This integral is then approximated using a~1D~quadrature rule~\eqref{eq:quadrature} on the reference interval~$(0,1)$
\begin{equation*}
\left[\vec{K}_\mathrm{D}\right]_{(k-1)N+i} \approx 
\sum_{E_{kn}\in\partial T_k\cap\setE_{\partial\Omega}}|E_{kn}| \sum_{r=1}^R \omega_r
\hat{\vphi}_i\circ\hat{\vec{\gamma}}_n(q_r) \,c_{\mathrm{D}}(t, \vec{F}_k\circ\hat{\vec{\gamma}}_n(q_r)) 
\left( \hat{\vec{u}}_{kn}(q_r) \cdot \nu_{kn}^m\right)
\delta_{\hat{\vec{u}}_{kn}\cdot\nu_{kn}^m<0}
\end{equation*}
allowing to vectorize the computation over all triangles and resulting in the routine~\code{assembleVecEdgePhiIntFuncContVal}.

\subsection{Slope limiters}\label{sec:implementation:slope-lim}
For the implementation of the slope limiters described in Sec.~\ref{sec:slope-lim}, three parts must be considered:
\begin{enumerate}
\item The assembly of the transformation matrix~$\vecc{M}^\mathrm{DG,Taylor}$ from Eq.~\eqref{eq:MDGTaylor};
\item the slope limiters themselves;
\item the selective mass lumping in the limiting of time-derivatives, as explained in Sec.~\ref{sec:slope-lim:time}.
\end{enumerate}

\subsubsection[Assembly of DG-Taylor transformation matrix]{Assembly of $\vecc{M}^\mathrm{DG,Taylor}$}
The entries in~$\vecc{M}^\mathrm{DG,Taylor}$ are transformed using transformation rule~\eqref{eq:trafoRule:T}, and the integral is then approximated using a~2D quadrature rule~\eqref{eq:quadrature} on the reference triangle~$\hat{T}$
\begin{equation*}
\left[\vecc{M}^\mathrm{DG,Taylor}\right]_{(k-1)N+i,(k-1)N+j} \;=\;
\int_{T_k} \vphi_{ki} \,\phi_{kj} \,\dd\vec{x} \;=\;
2|T_k| \int_{\hat{T}} \hat{\vphi}_i \,\phi_{kj} \circ \vec{F}_k(\hat{\vec{x}})\, \dd\hat{\vec{x}} \;\approx\;
2 |T_k| \sum_{r=1}^R \hat{\vphi}_i(\vec{q}_r) \, \phi_{kj} \circ \vec{F}_k(\vec{q}_r)\,.
\end{equation*}
Recall that we cannot define the Taylor basis~\eqref{eq:taylor:basis} on the reference triangle~$\hat{T}$ since the basis functions depend directly on the physical coordinates of the element.
The assembly is vectorized over all triangles resulting in the routine \code{assembleMatElemPhiDiscPhiTaylor}.
Matrices~$\vecc{M}^\mathrm{DG,Taylor}$ and~$\vecc{M}$ assembled in Sec.~\ref{sec:assembly:globM} are then used in the routines \code{projectDataDisc2DataTaylor} and \code{projectDataTaylor2DataDisc} to transform the representation matrix between the modal and Taylor basis---as explained in Sec.~\ref{sec:taylor:basis}---by solving system~\eqref{eq:taylor:trafo}.

\subsubsection{Slope limiting operators}
The slope limiters themselves are implemented in a generic manner according to Sec.~\ref{sec:slope-lim}.
Beginning with the highest-order derivatives, we evaluate the linear or the full reconstruction of these derivatives at all control points~$\vec{x}_{ki}$.
These values along with the centroid values are used to determine the minimum and maximum values~$c_{ki}^\mathrm{min}$ and~$c_{ki}^\mathrm{max}$ (see Eq.~\eqref{eq:slope-lim:lin:cond}) for each control point and, from these, the element-wise correction factors~$\alpha_e$ according to Eq.~\eqref{eq:slope-lim:lin:corr} or~\eqref{eq:slope-lim:kuzmin:corr} are calculated.
This is implemented in routine \code{computeVertexBasedLimiter}.
Depending on the limiter type, this value is then applied to a~certain subset of degrees of freedom, and the computation is repeated for the next lower-order derivatives.
The full slope limiting operators $\vecc{\Phi}$, $\vecc{\Phi}^\mathrm{Taylor}$ are provided by the functions \code{applySlopeLimiterDisc} or \code{applySlopeLimiterTaylor}, respectively.
\par
To overcome numerical problems for cases where~$c_{ki} \approx c_{ki}^\mathrm{max}$ or~$c_{ki} \approx c_{ki}^\mathrm{min}$, we modify Eqs.~\eqref{eq:slope-lim:lin:corr} and~\eqref{eq:slope-lim:kuzmin:corr} so that the condition for cases~1 and~3 are modified to $c_{ki} > c_{ki}^\mathrm{max} - \epsilon$ and $c_{ki} < c_{ki}^\mathrm{min} + \epsilon$ for a~small $0 < \epsilon \in \IR$.
Additionally, we increase the absolute value of the denominators in cases 1 and 3 by~$\epsilon$, i.\,e., add or subtract~$\epsilon$, respectively.
This makes our limiter slightly more strict than the one given by condition~\eqref{eq:slope-lim:lin:cond}; $\epsilon = 10^{-8}$ has been found to be a~suitable value for double precision computations.
Moreover, we enforce $0\le\alpha_{ke}\le 1$ to overcome cases where the division by a~number close to zero might lead to values $\alpha_{ke}>1$.
To reduce execution time, we perform the necessary evaluation of the Taylor basis functions in all vertices only once and store the result in a~global variable that is used in the slope limiting routines.
This is implemented in \code{computeTaylorBasesV0T}.

\subsubsection{Limiting time-derivatives}
To obtain the selectively lumped time-derivative~$\tilde{\vec{C}}^{(i)}$ required for update scheme~\eqref{eq:slope-lim:SSP-RK}, we compute the discrete time derivative $\dot{\vec{C}}^{(i)}$ as given in Sec.~\ref{sec:slope-lim:time} and transform it to a~Taylor basis representation using Eq.~\eqref{eq:taylor:trafo} (implemented in routine \code{projectDataDisc2DataTaylor}).
We compute the stationary matrix
\begin{equation*}
\vecc{M}^\mathrm{corr} \coloneqq \vecc{M}_\mathrm{L}^{-1} \vecc{M}
\end{equation*}
in the beginning (which is computationally cheap, since $\vecc{M}_\mathrm{L}$ is a~diagonal matrix) and calculate the selectively lumped time-derivative~$\tilde{\vec{C}}^{\mathrm{Taylor},(i)}$ from Equation~\eqref{eq:slope-lim:lumped-update}.
Backtransformation to the modal DG~basis is again cheap, since it requires only the inverse of the diagonal matrix~$\vecc{M}$ (implemented in \code{projectDataTaylor2DataDisc}) to produce the selectively lumped time-derivative~$\tilde{\vec{C}}^{(i)}$.
Each intermediate solution~$\vec{C}^{(i)}$ in update scheme~\eqref{eq:slope-lim:SSP-RK} is then obtained by applying the slope limiting operator~$\vecc{\Phi}$ implemented in \code{applySlopeLimiterDisc}.

\subsection{Numerical results}\label{sec:codeverification}

\subsubsection{Analytical convergence test}\label{sec:codeverification:analytical}
The code is verified by showing that the numerically estimated orders of convergences match the analytically predicted ones for prescribed smooth solutions. 
To verify the spatial discretization, we restrict ourselves to the stationary version of~\eqref{eq:model} and investigate the impact of different slope limiting schemes on the convergence order.
The effectiveness of the slope limiters is verified in the next section.
\par
We choose the exact solution~$c(\vec{x}) \coloneqq \cos(7x^1)\,\cos(7x^2)$ and velocity field~$\vec{u}(\vec{x}) \coloneqq [\exp((x^1+x^2)/2), \exp((x^1-x^2)/2)]^\mathrm{T}$ on the domain~$\Omega\coloneqq (0,1)^2$.
The data~$c_\mathrm{D}$ and~$f$ are derived analytically by inserting~$c$ and~$\vec{u}$ into~\eqref{eq:model}.  
We then compute the solution~$c_{h_j}$ for a~sequence of increasingly finer meshes with element sizes~$h_j$, where the coarsest grid~$\setT_{h_0}$ covering~$\Omega$ is an~irregular grid, and each finer grid is obtained by regular refinement of its predecessor.
The \emph{discretization error}~$\|c_h(t)-c(t)\|_{L^2(\Omega)}$ at time~$t\in J$ is computed as the $L^2$-norm of the difference between the discrete solution~$c_h(t)$ and the analytical solution~$c(t)$, as it was described in detail in our first paper~\cite{FESTUNG1}.
Table~\ref{tab:verification:stationary} contains the results demonstrating the experimental order of convergence~$\alpha$ estimated using
\begin{equation*}
\alpha\;\coloneqq\;\ln\Bigg(\frac{\|c_{h_{j-1}} - c\|_{L^2(\Omega)}}{\|c_{h_j} - c\|_{L^2(\Omega)}}\Bigg)\Bigg/\ln\Bigg(\frac{h_{j-1}}{h_j}\Bigg)\;.
\end{equation*}
\begin{table}[ht!]
\small
\begin{tabularx}{\linewidth}{@{}clcCcCcCcCcc@{}}\toprule
& $p$ & 0 & 0 & 1 & 1 & 2 & 2 & 3 & 3 & 4 & 4\\
limiter & $j$ & $\|c_h-c\|$ & $\alpha$ & $\|c_h-c\|$ & $\alpha$ & $\|c_h-c\|$ & $\alpha$ & $\|c_h-c\|$ & $\alpha$ & $\|c_h-c\|$ & $\alpha$\\\midrule
\multirow{6}{*}{\begin{sideways}none\end{sideways}} %
& 0  & 3.64e--1 &   --- & 5.44e--1 &   --- &  6.98e--1 &   --- &  2.73e--1 &   --- &  2.21e--1 &   --- \\
& 1  & 5.64e--1 &--0.63 & 4.25e--1 &  0.36 &  1.32e--1 &  2.40 &  8.51e--2 &  1.68 &  1.46e--2 &  3.92 \\
& 2  & 3.49e--1 &  0.69 & 1.28e--1 &  1.73 &  2.82e--2 &  2.23 &  6.10e--3 &  3.80 &  1.25e--3 &  3.54 \\
& 3  & 2.49e--1 &  0.49 & 3.65e--2 &  1.81 &  3.15e--3 &  3.16 &  3.87e--4 &  3.98 &  4.58e--5 &  4.77 \\
& 4  & 1.73e--1 &  0.53 & 8.29e--3 &  2.14 &  3.68e--4 &  3.10 &  2.39e--5 &  4.02 &  1.51e--6 &  4.92 \\
& 5  & 1.09e--1 &  0.66 & 1.93e--3 &  2.11 &  4.50e--5 &  3.03 &  1.50e--6 &  4.00 &  4.80e--8 &  4.97 \\
& 6  & 6.35e--2 &  0.78 & 4.71e--4 &  2.03 &  5.59e--6 &  3.01 &  9.36e--8 &  4.00 &  1.51e--9 &  4.99 \\
\midrule
\multirow{6}{*}{\begin{sideways}linear\end{sideways}} %
& 0  & --- & --- &  5.15e--1 &   --- &  7.21e--1 &   --- &  4.60e--1 &   --- &  5.27e--1 &   --- \\
& 1  & --- & --- &  4.53e--1 &  0.19 &  3.55e--1 &  1.02 &  3.54e--1 &  0.38 &  3.59e--1 &  0.55 \\
& 2  & --- & --- &  1.64e--1 &  1.47 &  1.34e--1 &  1.41 &  1.51e--1 &  1.23 &  1.47e--1 &  1.29 \\
& 3  & --- & --- &  4.53e--2 &  1.85 &  3.12e--2 &  2.10 &  3.12e--2 &  2.27 &  3.12e--2 &  2.23 \\
& 4  & --- & --- &  8.76e--3 &  2.37 &  4.48e--3 &  2.80 &  4.46e--3 &  2.81 &  4.46e--3 &  2.81 \\
& 5  & --- & --- &  1.98e--3 &  2.14 &  6.20e--4 &  2.85 &  6.19e--4 &  2.85 &  6.19e--4 &  2.85 \\
& 6  & --- & --- &  4.74e--4 &  2.07 &  7.63e--5 &  3.02 &  7.61e--5 &  3.02 &  7.61e--5 &  3.02 \\
\midrule
\multirow{6}{*}{\begin{sideways}hier. vert. based\end{sideways}} %
& 0  & --- & --- &  5.15e--1 &   --- &  6.98e--1 &   --- &  2.73e--1 &   --- &  2.21e--1 &   --- \\
& 1  & --- & --- &  4.53e--1 &  0.19 &  2.47e--1 &  1.50 &  2.57e--1 &  0.09 &  1.33e--1 &  0.73 \\
& 2  & --- & --- &  1.64e--1 &  1.47 &  1.06e--1 &  1.22 &  5.89e--2 &  2.12 &  5.10e--2 &  1.38 \\
& 3  & --- & --- &  4.53e--2 &  1.85 &  7.87e--3 &  3.76 &  4.02e--3 &  3.87 &  3.17e--3 &  4.01 \\
& 4  & --- & --- &  8.76e--3 &  2.37 &  6.39e--4 &  3.62 &  6.16e--5 &  6.03 &  6.53e--6 &  8.92 \\
& 5  & --- & --- &  1.98e--3 &  2.14 &  5.32e--5 &  3.59 &  2.62e--6 &  4.56 &  1.08e--7 &  5.92 \\
& 6  & --- & --- &  4.74e--4 &  2.07 &  5.94e--6 &  3.16 &  1.17e--7 &  4.48 &  5.48e--9 &  4.30 \\
\midrule
\multirow{6}{*}{\begin{sideways}strict\end{sideways}} %
& 0  & --- & --- &  5.15e--1 &   --- &  6.84e--1 &   --- &  3.28e--1 &   --- &  3.56e--1 &   --- \\
& 1  & --- & --- &  4.53e--1 &  0.19 &  3.79e--1 &  0.85 &  3.69e--1 &--0.17 &  3.54e--1 &  0.00 \\
& 2  & --- & --- &  1.64e--1 &  1.47 &  1.14e--1 &  1.74 &  6.74e--2 &  2.45 &  6.15e--2 &  2.53 \\
& 3  & --- & --- &  4.53e--2 &  1.85 &  1.23e--2 &  3.21 &  1.20e--2 &  2.49 &  1.16e--2 &  2.40 \\
& 4  & --- & --- &  8.76e--3 &  2.37 &  1.15e--3 &  3.41 &  8.99e--4 &  3.74 &  8.85e--4 &  3.72 \\
& 5  & --- & --- &  1.98e--3 &  2.14 &  1.19e--4 &  3.27 &  1.33e--4 &  2.76 &  1.33e--4 &  2.74 \\
& 6  & --- & --- &  4.74e--4 &  2.07 &  6.64e--6 &  4.17 &  1.09e--6 &  6.93 &  1.13e--6 &  6.88 \\
\bottomrule
\end{tabularx}
\caption{Discretization errors measured in~$L^2(\Omega)$ and estimated orders of convergences for different polynomial degrees and limiter types. We have~$h_j = \frac{1}{3\cdot2^j}$ and $K=36\cdot4^j$ triangles in the \mbox{$j$th} refinement level.}
\label{tab:verification:stationary}
\end{table}

\par
The interesting points to compare in Table \ref{tab:verification:stationary} are the errors for different limiting schemes.
Whereas the linear approximation does not get limited at all on finer grids---indicating that the vertex-based limiters are able to distinguish between smooth and non-smooth solutions given sufficient mesh resolution, the case of higher order approximation warrants a closer look. 
First of all, the linear limiter takes a heavy toll on the higher order degrees of freedom: quadratic, cubic, and quartic solutions produce virtually the same error, thus negating the effect of a more accurate DG~solution. 
The strict limiter seems to flatten out after the cubic approximation, thus gaining one order of convergence on the linear limiter. 
The hierarchical limiter of Kuzmin, however, appears to perform very well for all tested meshes and approximation orders having a very small effect on the errors of the analytical test case.

\subsubsection{Solid body rotation}\label{sec:codeverification:rotating}
\begin{figure}[ht!]
\centering
\begin{subfigure}[t]{.4\textwidth}
\includegraphics[width=\textwidth]{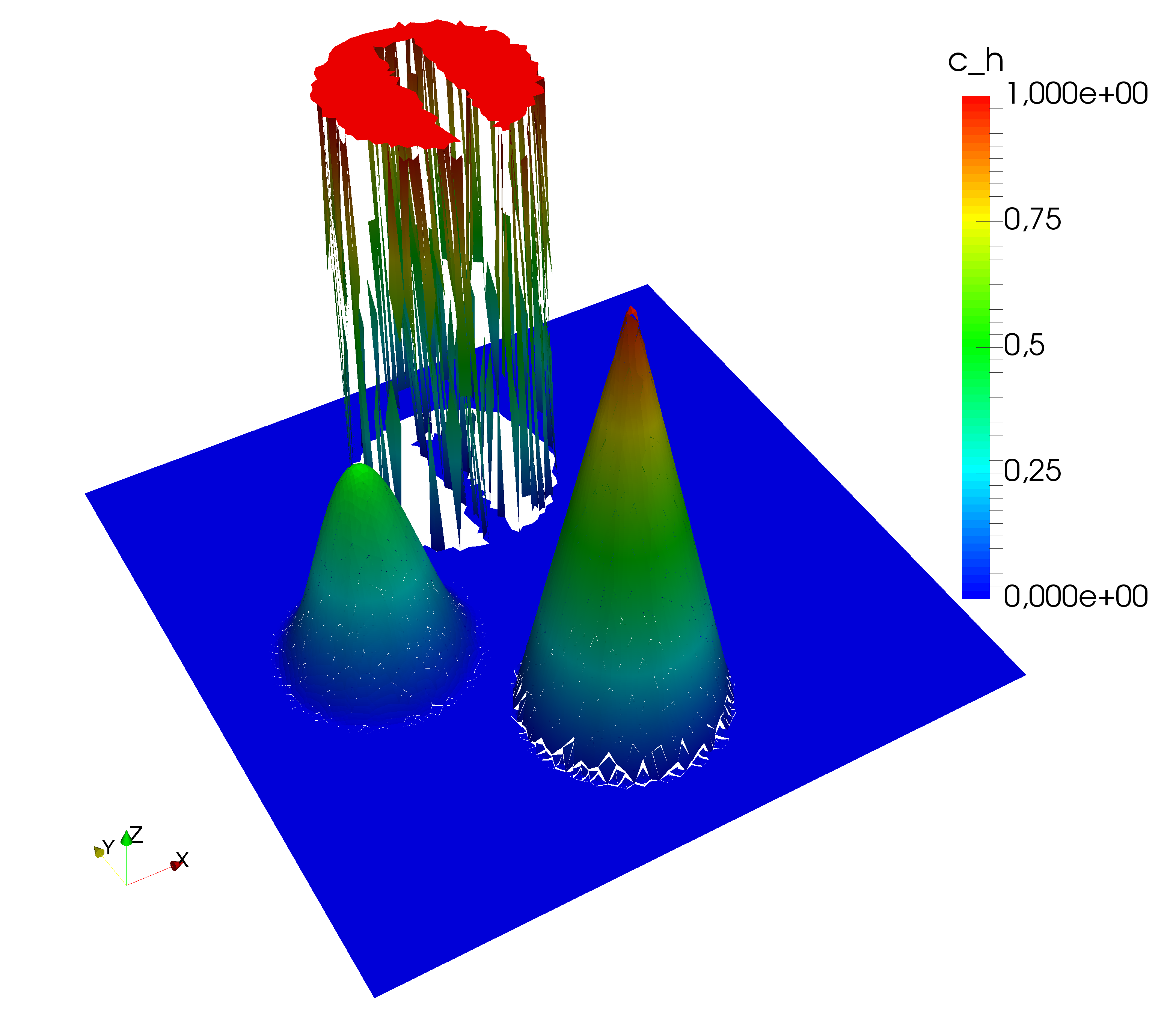}
\captionsetup{justification=raggedright,singlelinecheck=false}
\caption{Projected and limited initial data.}
\label{fig:solid-body:3d:initial}
\end{subfigure}%
~%
\begin{subfigure}[t]{.4\textwidth}
\includegraphics[width=\textwidth]{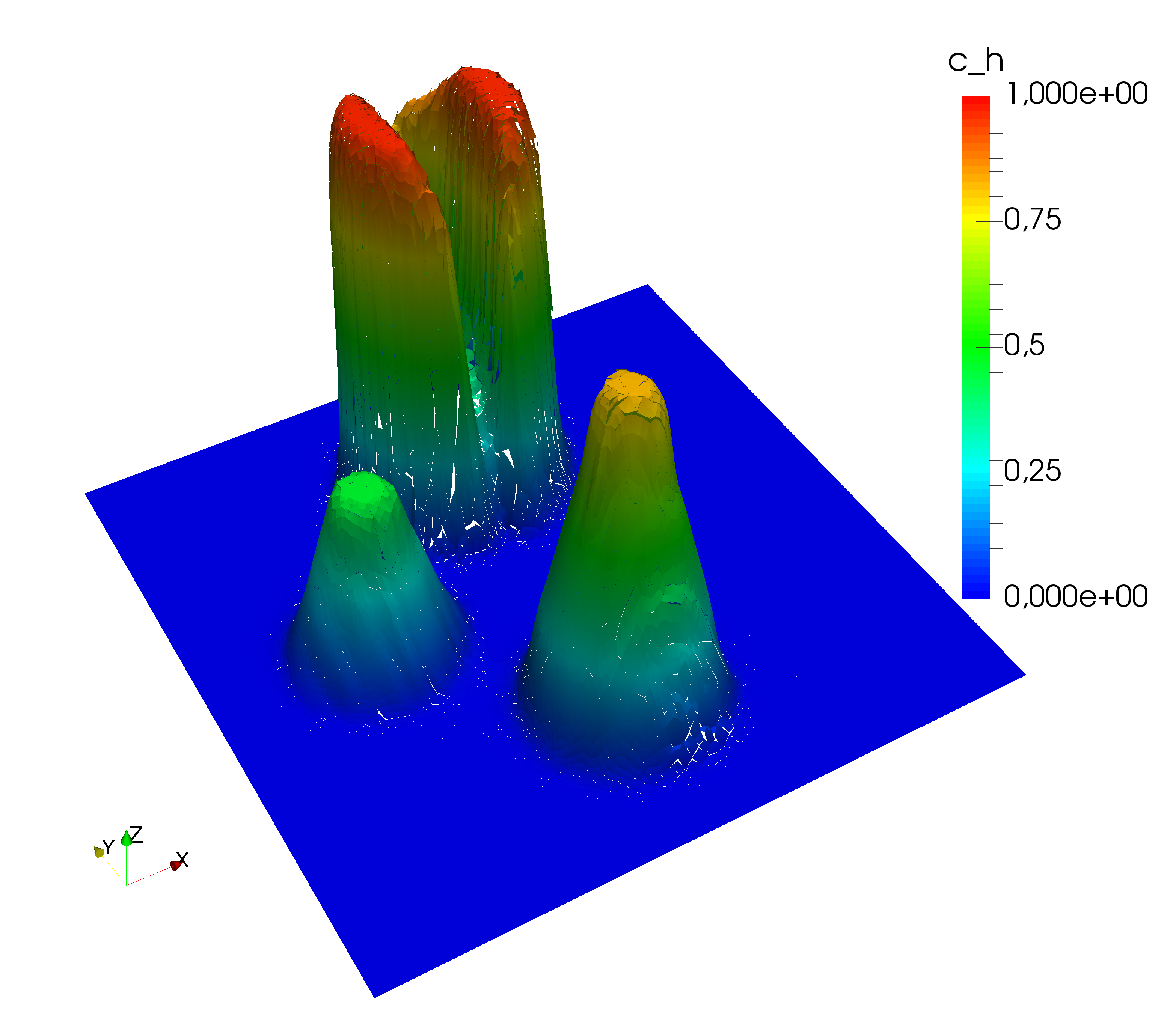}
\caption{Linear limiter, no lumped time-derivative.}
\label{fig:solid-body:3d:linear-non-lumped}
\end{subfigure}%
\\
\begin{subfigure}[t]{.4\textwidth}
\includegraphics[width=\textwidth]{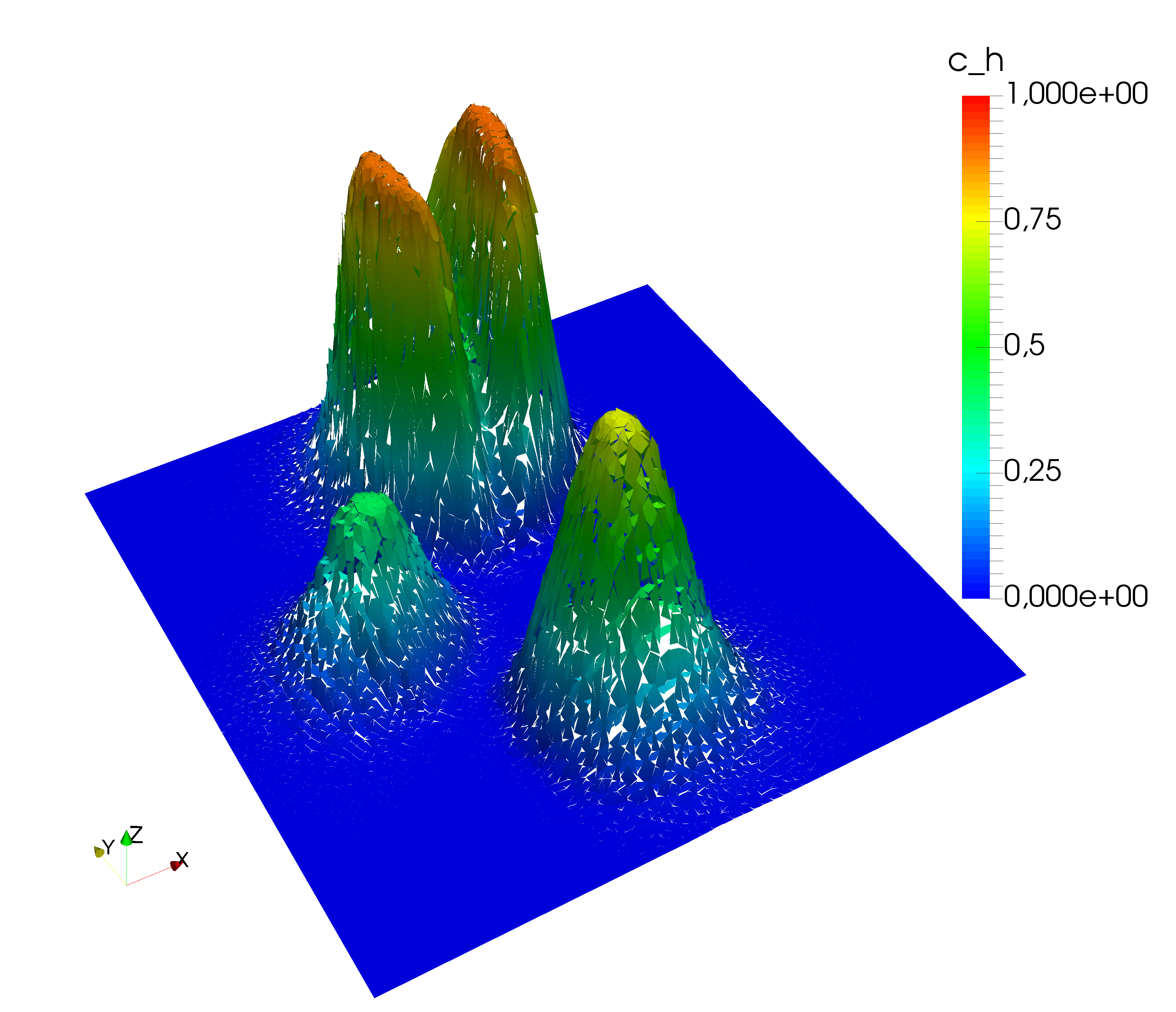}
\caption{Hier. vert.-based lim., no lumped time-derivative.}
\label{fig:solid-body:3d:kuzmin-non-lumped}
\end{subfigure}%
~%
\begin{subfigure}[t]{.4\textwidth}
\includegraphics[width=\textwidth]{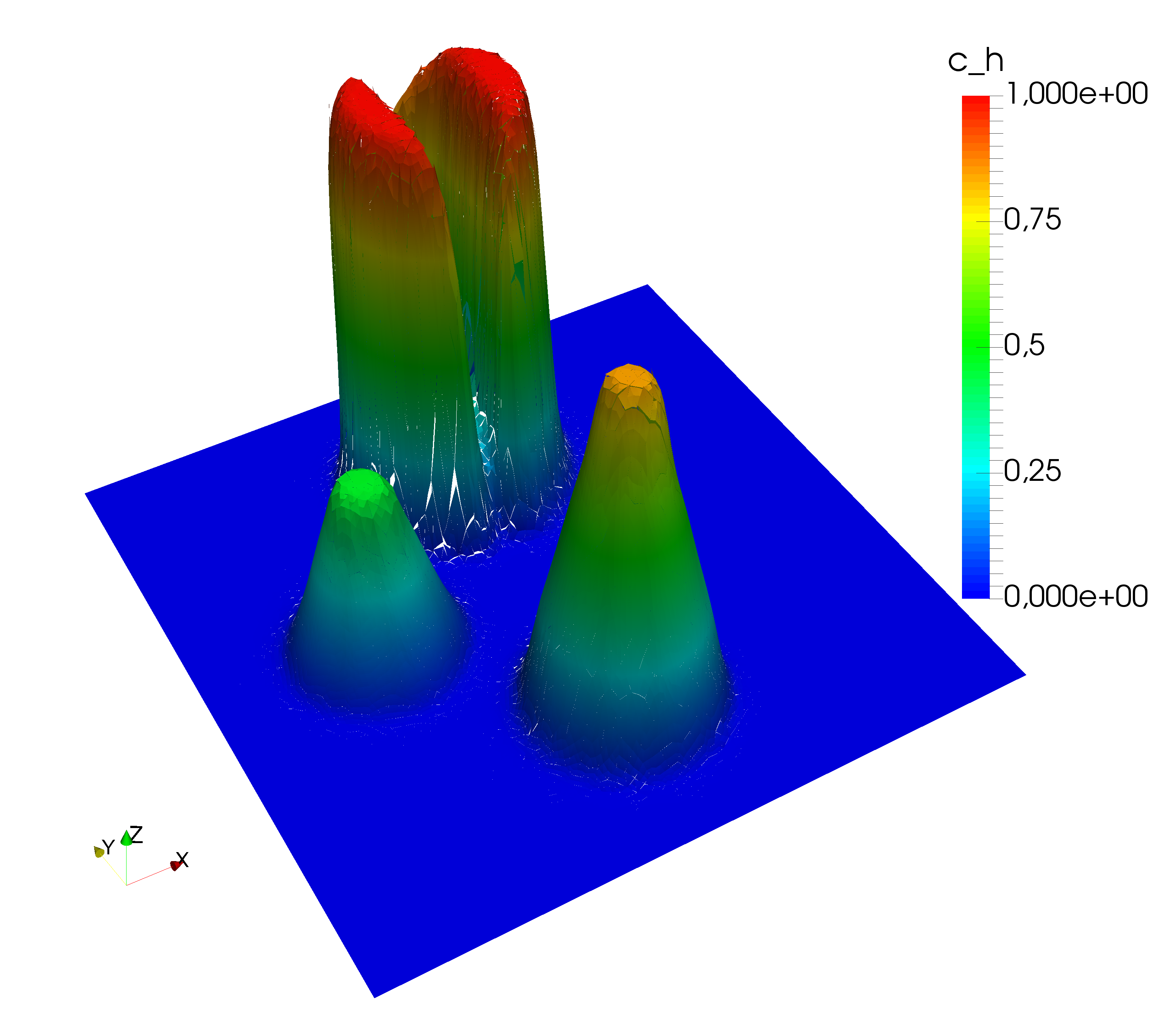}
\caption{Linear limiter, lumped time derivative.}
\label{fig:solid-body:3d:linear}
\end{subfigure}%
\\
\begin{subfigure}[t]{.4\textwidth}
\includegraphics[width=\textwidth]{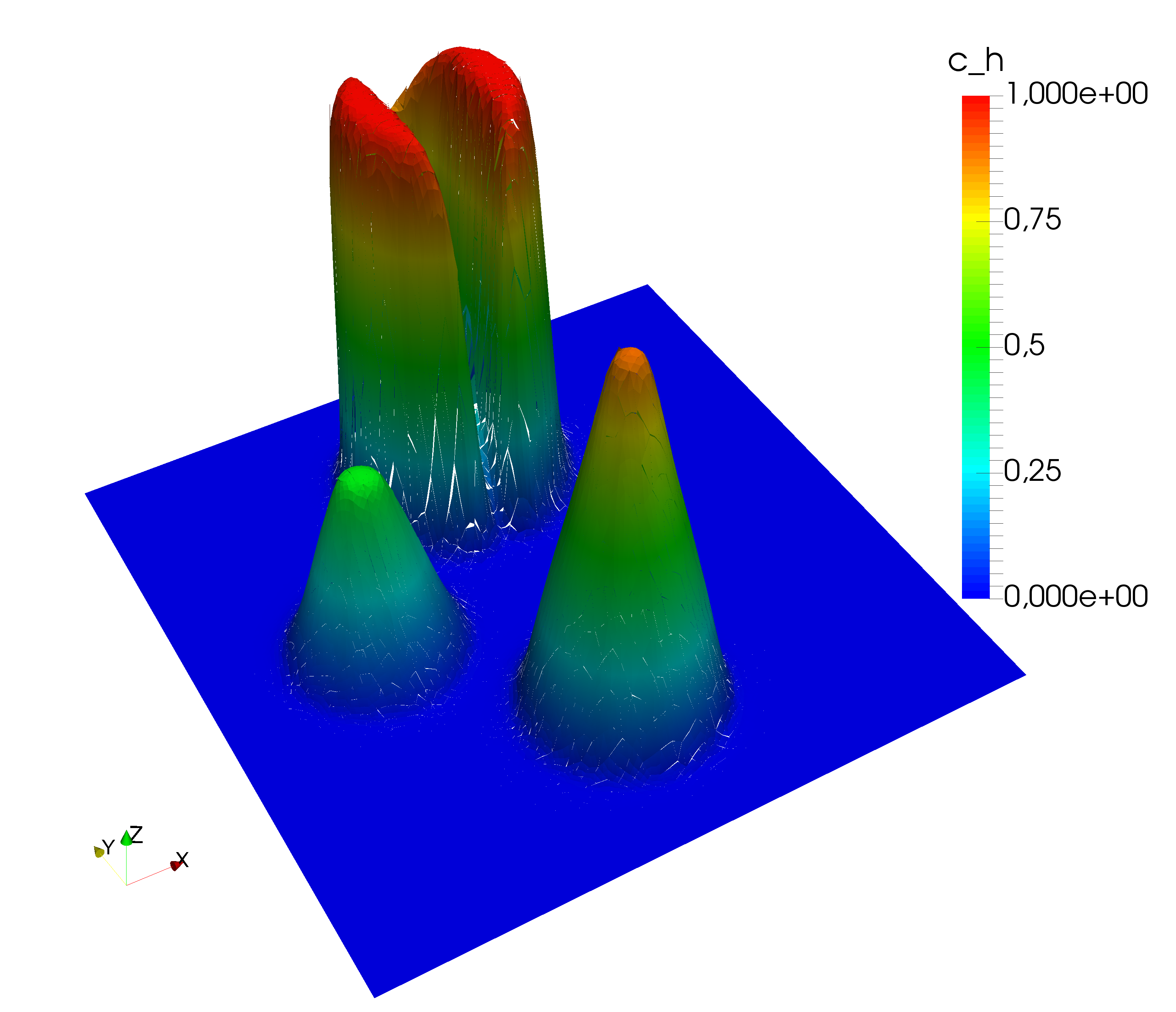}
\caption{Hier. vert.-based limiter, lumped time derivative.}
\label{fig:solid-body:3d:kuzmin}
\end{subfigure}%
~%
\begin{subfigure}[t]{.4\textwidth}
\includegraphics[width=\textwidth]{images/Rotation_h2e-6_p2_linear_lumped}
\caption{Strict limiter, lumped time derivative.}
\label{fig:solid-body:3d:strict}
\end{subfigure}%
\captionsetup{justification=raggedright,singlelinecheck=false}
\caption{Visualization of the DG~solutions with $p=2$ at end time~$t_\mathrm{end}=2\pi$.}
\label{fig:solid-body:3d}
\end{figure}

As a~benchmark problem, we use solid body rotation test proposed by  LeVeque~\cite{LeVeque1996}, which is often used to investigate limiter performance~\cite{Kuzmin2012,Michoski2011}.
It consists of a~slotted cylinder, a~sharp cone, and a~smooth hump (see Figure~\ref{fig:solid-body:3d:initial}) that are placed in a square domain~$\Omega = [0,1]^2$ and transported by a~time-independent velocity field~$\vec{u}(\vec{x})=\transpose{[0.5 - x^2, x^1 - 0.5]}$ in a~counterclockwise rotation over~$J=(0,2\pi)$.
With $r = 0.0225$ and
$G(\vec{x},\vec{x}_0) \coloneqq \frac{1}{0.15} \|\vec{x}-\vec{x}_0\|_2$,
we choose initial data satisfying
\begin{equation*}
c^0(\vec{x}) = \left\{
\begin{array}{lll}
  1                                               & \quad\text{if}\quad \parbox{.36\textwidth}{$(x^1 - 0.5)^2 + (x^2 - 0.75)^2 \le r$\\$\land\;(x^1\le0.475 \lor x^1\ge0.525 \lor x^2\ge0.85)$} & \mbox{(slotted cylinder)} \\
  1-G(\vec{x},\transpose{[0.5,0.25]})                       & \quad\text{if}\quad (x^1 - 0.5)^2 + (x^2 - 0.25)^2 \le r & \mbox{(sharp cone)}       \\
  \frac{1}{4}(1+\cos(\pi G(\vec{x},\transpose{[0.25,0.5]})) & \quad\text{if}\quad (x^1 - 0.25)^2 + (x^2 - 0.5)^2 \le r & \mbox{(smooth hump)}      \\
  0                                               & \quad\text{otherwise} & 
\end{array}
\right\}
\end{equation*}
and zero boundary $c_\mathrm{D} = 0$ and right-hand side $f = 0$.
\par
\begin{figure}[ht!]
\centering
\includegraphics[width=.8\textwidth]{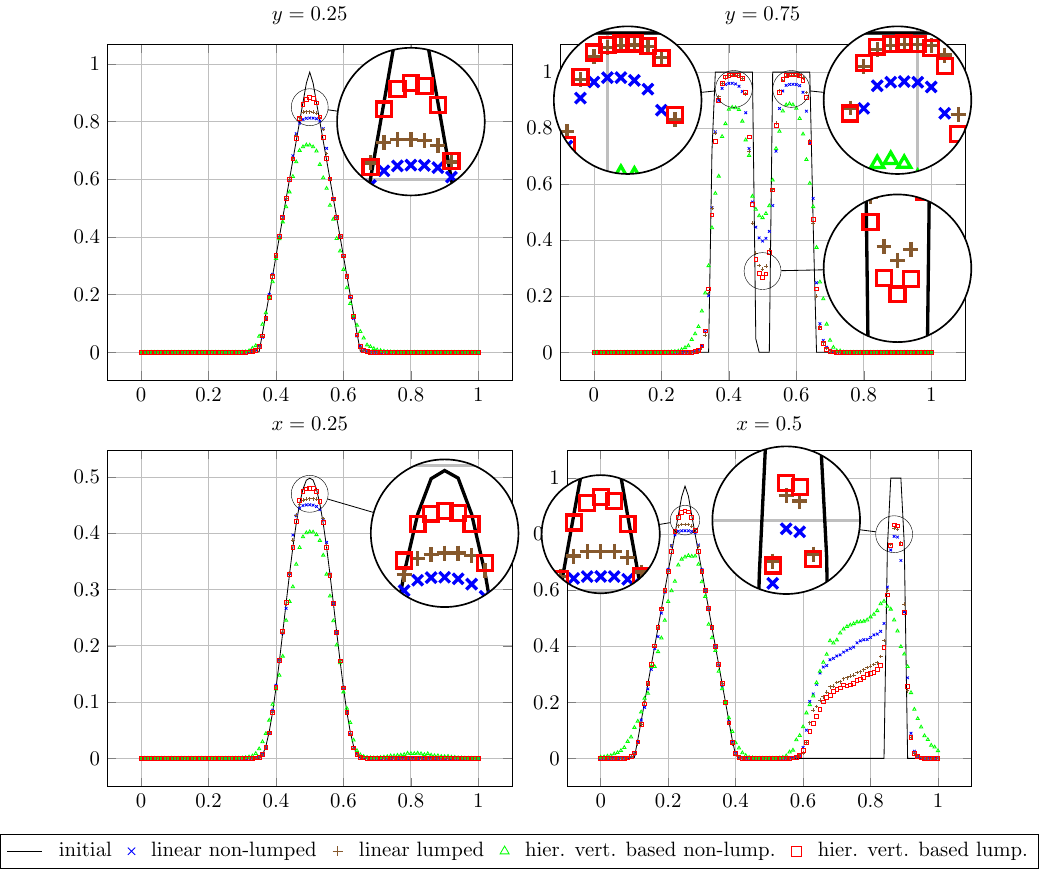}
\caption{Intersection lines of the DG~solution with $p=2$ at end time $t_\mathrm{end}=2\pi$ with and without the selectively lumped time-derivative, as described in Section.~\ref{sec:slope-lim:time}. The large deviation of the solutions from the initial data in the lower right plot are due to fill-in of the slot in the cylinder.}
\label{fig:solid-body:lumped-vs-non-lumped}
\end{figure}

First, we would like to emphasize the huge improvement of the solution when applying Kuzmin's selectively lumped time-stepping scheme (cf.~Sec.~\ref{sec:slope-lim:time}), which reduces numerical diffusion and peak clipping, as visible in the intersection lines in Fig.~\ref{fig:solid-body:lumped-vs-non-lumped} and the 3D~visualization (Figs.~\ref{fig:solid-body:3d:linear-non-lumped}--\ref{fig:solid-body:3d:kuzmin}) and results in a~smaller error (see Tab.~\ref{tab:solid-body}).
Without this technique, the implicit coupling of the derivatives renders the higher order limiter inferior to the linear limiter.
\par
\begin{figure}[ht!]
\centering
\includegraphics[width=.7\textwidth]{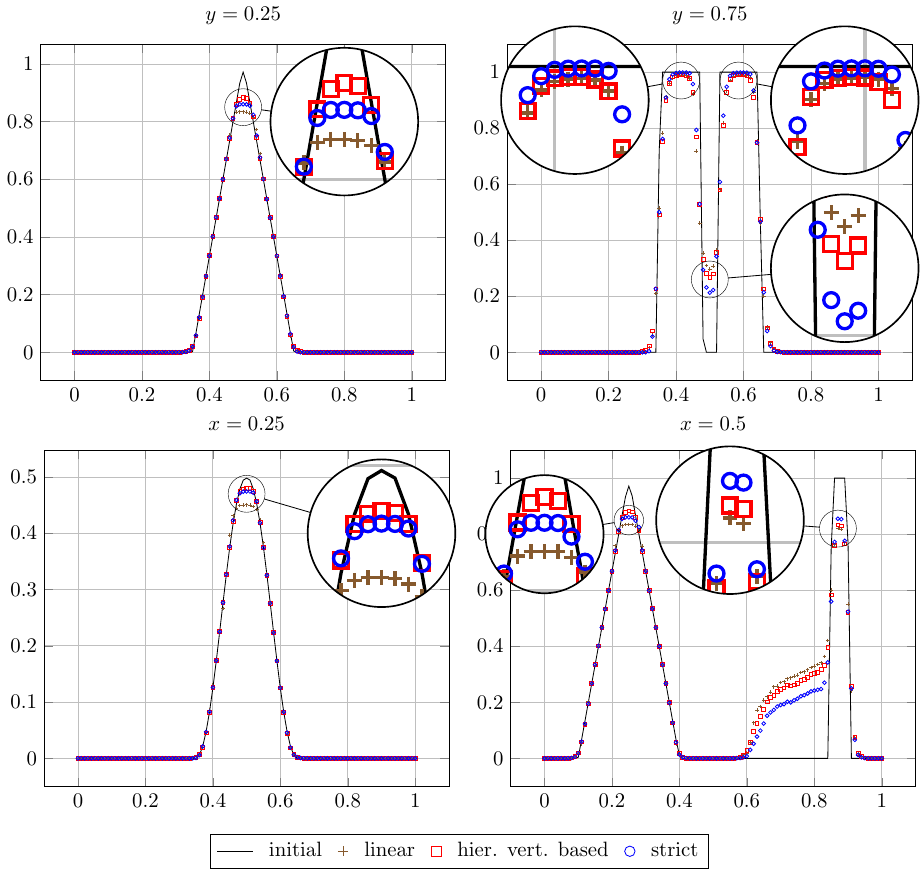}
\caption{Intersection lines of the DG~solution with $p=2$ at end time~$t_\mathrm{end}=2\pi$ when applying the lumped time-stepping scheme (cf.~Sec.~\ref{sec:slope-lim:time}). Compared are solutions with the different limiters described in Section~\ref{sec:slope-lim}.}
\label{fig:solid-body:lumped}
\end{figure}

The results for the different slope limiters (see Fig.~\ref{fig:solid-body:3d:linear}--\ref{fig:solid-body:3d:strict}) when applying the lumped time-stepping scheme are similar with all limiters producing errors in the same range (see Tab.~\ref{tab:solid-body}).
The linear vertex based limiter exhibits significantly stronger peak clipping than the higher order limiters in the intersection lines of Fig.~\ref{fig:solid-body:lumped}.
When looking at the slotted cylinder the strict limiter outperforms the others, producing less fill-in in the slot and better preserving the shape of the cylinder.
Also, the total error is lowest for the strict limiter.
\par
For polynomial degrees $p\ge 2$ the solution might still violate the bounds along the edges when using only control points in the vertices of the element for the slope limiting procedure.
When evaluating the numerical fluxes, these values are transported to the neighboring elements and can lead to cell averages lying outside of the initial bounds in the next time step.
Table~\ref{tab:solid-body} shows that all limiter types suffer from this, however these small violations of the bounds usually smooth out over time due to the effects of numerical diffusion and do not introduce numerical problems.
The linear limiter is the only that preserves the upper bound but comes at the price of stronger peak-clipping.
To effectively restrict the solution to these bounds for $p=2$, additional control points on the edges are necessary, which leads to a~significant increase in numerical diffusion, hence this is not a~suitable technique.
\begin{table}[ht!]
\small
\begin{tabularx}{\linewidth}{@{}l|cC|CCC@{}}
\toprule
& \multicolumn{2}{c|}{\textit{non-lumped}} & \multicolumn{3}{c}{\textit{lumped}} \\
& linear & hier. vertex-based & linear & hier. vertex-based & strict \\
\midrule
$\min_k c_h(\vec{x}_{kc})$       &      0.00 &--1.11e--5 &--3.95e--10 &--6.18e--5 &--6.42e--4 \\
$\min_{k,i} c_h(\vec{x}_{ki})$   &--2.93e--2 &--2.47e--2 &--9.53e--3  &--9.69e--3 &--9.79e--16 \\
$\min_{k,j} c_h(\vec{x}_{e,kj})$ &--3.32e--8 &--4.28e--3 &--1.18e--7 &--4.28e--3 &--5.35e--3 \\
$\max_k c_h(\vec{x}_{kc})$       &   1.00000 &   1.00000 &   1.00000 &   1.00000 &   1.00142 \\
$\max_{k,i} c_h(\vec{x}_{ki})$   &   1.00000 &   1.00915 &   1.00000 &   1.01760 &   1.00000 \\
$\max_{k,j} c_h(\vec{x}_{e,kj})$ &   1.00000 &   1.00014 &   1.00000 &   1.00297 &   1.00296 \\
$\|c_h(0) - c^0\|_{L^2(\Omega)}$             & 3.73e--2  &  3.66e--2 &  3.73e--2 &  3.66e--2 &  3.66e--2 \\
$\|c_h(t_\mathrm{end}) - c^0\|_{L^2(\Omega)}$& 8.18e--2  &  1.15e--1 &  7.38e--2 &  7.40e--2 &  7.07e--2 \\
\bottomrule
\end{tabularx}
\caption{Minimum and maximum values in centroids~$\vec{x}_{kc}$, nodes~$\vec{x}_{ki}$, and edge midpoints~$\vec{x}_{e,kj}$ with $p=2$ over the entire simulation time $J=(0,2\pi)$. Two last rows give the $L^2$-errors of the slope-limited projection of initial data and the end solution when compared to the analytical initial data.}
\label{tab:solid-body}
\end{table}

\section{Register of Routines}\label{sec:routines}
We list here all routines of our implementation that were added since the first paper~\cite{FESTUNG1} in alphabetic order.  
For the reason of compactness, we waive the check for correct function arguments, e.\,g., by means of routines as~\code{assert}.  
However, it is strongly recommended to catch exceptions if the code is to be extended. 
The argument~\code{g}~is always a~\code{struct} representing the triangulation~$\setT_h$, the argument~\code{N}~is always the number of local basis functions~$N$.  
A~script that demonstrates the application of the presented routines is given in~\code{mainAdvection.m}.
\par
Additionally, we recommend using Laurent Sorbers slightly faster implementation of the Kronecker product~\cite{kron} for \code{sparse} matrices from \Matlab\;File Exchange\footnote{\url{http://www.mathworks.com/matlabcentral/fileexchange/}} to speed up the computation.
Note that this implementation is only suitable for \code{logical} and \code{double} matrices.

\begin{lstlisting}[title={%
%\begin{lstlisting}[title={%
\code{dataDisc = applySlopeLimiterDisc(g, dataDisc, markV0TbdrD, dataV0T, globM, globMDiscTaylor, type)}
implements the slope limiting operator $\vecc{\Phi}$, as given in Eq.~\eqref{eq:slope-lim:op} with the chosen limiter type given as a~string in \code{type}. Parameter \code{dataDisc} is the representation matrix in modal DG~basis, \code{markV0TbdrD} is a~logical matrix marking all Dirichlet boundary nodes for which boundary data is given in \code{dataV0T}.
}%
%]
%...
%\end{lstlisting}
]
function dataDisc = applySlopeLimiterDisc(g,dataDisc,markV0TbdrD,dataV0T,globM,globMDiscTaylor,type)
dataTaylor = projectDataDisc2DataTaylor(dataDisc, globM, globMDiscTaylor);
dataTaylor = applySlopeLimiterTaylor(g, dataTaylor, markV0TbdrD, dataV0T, type);
dataDisc = projectDataTaylor2DataDisc(dataTaylor, globM, globMDiscTaylor);
end

\end{lstlisting}

\begin{lstlisting}[title={%
%\begin{lstlisting}[title={%
\code{dataTaylor = applySlopeLimiterTaylor(g, dataTaylor, markV0TbdrD, dataV0T, type)}
implements the slope limiting operator $\vecc{\Phi}^\mathrm{Taylor}$, as described in Sec.~\ref{sec:slope-lim:time}. Parameter \code{dataTaylor} is the representation matrix in Taylor basis, the other parameters are the same as for \code{applySlopeLimiterDisc}.
}%
%]
%...
%\end{lstlisting}
]
function dataTaylorLim = applySlopeLimiterTaylor(g, dataTaylor, markV0TbdrD, dataV0T, type)
switch type
  case 'linear'
    dataTaylorLim = applySlopeLimiterTaylorLinear(g, dataTaylor, markV0TbdrD, dataV0T);
  case 'hierarch_vert'
    dataTaylorLim = applySlopeLimiterTaylorHierarchicalVertex(g, dataTaylor, markV0TbdrD, dataV0T);
  case 'strict'
    dataTaylorLim = applySlopeLimiterTaylorStrict(g, dataTaylor, markV0TbdrD, dataV0T);
  otherwise
    error('Unknown limiter type');
end
end

\end{lstlisting}

\begin{lstlisting}[title={%
%\begin{lstlisting}[title={%
\code{dataTaylorLim = applySlopeLimiterTaylorHierarchicalVertex(g, dataTaylor, markV0TbdrD, dataV0T)}
applies the hierarchical vertex-based limiter as described in Sec.~\ref{sec:slope-lim:kuzmin}, with input parameters as for \code{applySlopeLimiterTaylor}.
}%
%]
%...
%\end{lstlisting}
]
function dataTaylorLim = applySlopeLimiterTaylorHierarchicalVertex(g, dataTaylor, markV0TbdrD, dataV0T)
global gPhiTaylorV0T
[K, N] = size(dataTaylor); p = (sqrt(8*N+1)-3)/2;
dataTaylorLim = zeros(size(dataTaylor)); dataTaylorLim(:, 1) = dataTaylor(:, 1);
alpha = zeros(K, 1);
for ord = p : -1 : 1
  alphaOrd = ones(K, 1);
  indDOF = ord*(ord+1)/2 + 1 : (ord+1)*(ord+2)/2;
  for i = 1 : ord
    mult = bsxfun(@plus, [ord - i, i - 1], multiindex(1));
    ind = mult2ind(mult);
    valV0T = computeFuncDiscAtPoints(dataTaylor(:, ind), gPhiTaylorV0T(:, :, 1:3));
    if ord > 1
      alphaTmp = computeVertexBasedLimiter(g, dataTaylor(:, ind(1)), valV0T, markV0TbdrD, valV0T);
    else
      alphaTmp = computeVertexBasedLimiter(g, dataTaylor(:, ind(1)), valV0T, markV0TbdrD, dataV0T);
    end
    alphaOrd = min(alphaOrd, alphaTmp);
  end
  alpha = max(alpha, alphaOrd);
  dataTaylorLim(:, indDOF) = bsxfun(@times, alpha, dataTaylor(:, indDOF));
end
end 

\end{lstlisting}

\begin{lstlisting}[title={%
%\begin{lstlisting}[title={%
\code{dataTaylorLim = applySlopeLimiterTaylorLinear(g, dataTaylor, markV0TbdrD, dataV0T)}
applies the linear vertex-based limiter as described in Sec.~\ref{sec:slope-lim:linear}, with input parameters as for \code{applySlopeLimiterTaylor}.
}%
%]
%...
%\end{lstlisting}
]
function dataTaylorLim = applySlopeLimiterTaylorLinear(g, dataTaylor, markV0TbdrD, dataV0T)
global gPhiTaylorV0T
alphaE = computeVertexBasedLimiter(g, dataTaylor(:, 1), computeFuncDiscAtPoints(dataTaylor(:, 1:3), gPhiTaylorV0T(:,:,1:3)), markV0TbdrD, dataV0T);
dataTaylorLim = [dataTaylor(:,1), bsxfun(@times, alphaE, dataTaylor(:, 2:3)), bsxfun(@times, alphaE == 1, dataTaylor(:, 4:end))];
end 

\end{lstlisting}

\begin{lstlisting}[title={%
%\begin{lstlisting}[title={%
\code{dataTaylorLim = applySlopeLimiterTaylorStrict(g,dataTaylor, markV0TbdrD, dataV0T)}
applies the stricter form of the hierarchical vertex-based limiter as described in Sec.~\ref{sec:slope-lim:strict}, with input parameters as for \code{applySlopeLimiterTaylor}.
}%
%]
%...
%\end{lstlisting}
]
function dataTaylorLim = applySlopeLimiterTaylorStrict(g, dataTaylor, markV0TbdrD, dataV0T)
global gPhiTaylorV0T
[K, N] = size(dataTaylor); p = (sqrt(8*N+1)-3)/2;
dataTaylorLim = dataTaylor;
for ord = p : -1 : 1
  alphaOrd = ones(K, 1);
  indDOF = ord * (ord + 1) / 2 + 1 : N;
  for i = 1 : ord
    pReconstruction = p - ord + 1;
    mult = bsxfun(@plus, [ord - i, i - 1], multiindex(pReconstruction)); ind = mult2ind(mult);
    valV0T = computeFuncDiscAtPoints(dataTaylorLim(:, ind), gPhiTaylorV0T(:, :, 1 : size(ind, 1)));
    if ord > 1
      alphaTmp = computeVertexBasedLimiter(g, dataTaylorLim(:, ind(1)), valV0T, markV0TbdrD, valV0T);
    else
      alphaTmp = computeVertexBasedLimiter(g, dataTaylorLim(:, ind(1)), valV0T, markV0TbdrD, dataV0T);
    end
    alphaOrd = min(alphaOrd, alphaTmp);
  end
  dataTaylorLim(:, indDOF) = bsxfun(@times, alphaOrd, dataTaylorLim(:, indDOF));
end
end 

\end{lstlisting}

\begin{lstlisting}[title={%
%\begin{lstlisting}[title={%
\code{ret = assembleMatEdgePhiPhiValUpwind(g, refEdgePhiIntPhiInt, refEdgePhiIntPhiExt, valOnQuad)} assembles the matrix $\vecc{R}$ according to Sec.~\ref{sec:assembly:globR}, containing edge integrals of products of two basis functions multiplied with for each quadrature point specified values, stored in \code{valOnQuad}, and where the upwind-sided value w.\,r.\,t.~\code{valOnQuad} is chosen. The input arguments \code{refEdgePhiIntPhiInt} and \code{refEdgePhiIntPhiExt} provide the local matrices~$\hat{\vecc{R}}^\mathrm{diag}$ and~$\hat{\vecc{R}}^\mathrm{offdiag}$ (multidimensional arrays), respectively.
}%
%]
%...
%\end{lstlisting}
]
function ret = assembleMatEdgePhiPhiValUpwind(g, refEdgePhiIntPhiIntOnQuad, refEdgePhiIntPhiExtOnQuad, valOnQuad)
K = g.numT; N = size(refEdgePhiIntPhiIntOnQuad, 1);
ret = sparse(K*N, K*N);
p = (sqrt(8*N+1)-3)/2;
qOrd = 2*p+1;  [~, W] = quadRule1D(qOrd);
for nn = 1 : 3
  Rkn = g.areaE0T(:, nn);
  for r = 1 : length(W)  
    ret = ret + kron(spdiags(W(r) .* Rkn .* valOnQuad(:, nn, r) .* (valOnQuad(:, nn, r) > 0), 0, K, K), refEdgePhiIntPhiIntOnQuad(:, :, nn, r));
  end
  for np = 1 : 3  
    RknTimesVal = sparse(K*N, N);
    for r = 1 : length(W)
      RknTimesVal = RknTimesVal + kron(W(r) .* Rkn .* valOnQuad(:, nn, r) .* sparse(valOnQuad(:, nn, r) < 0), refEdgePhiIntPhiExtOnQuad(:, :, nn, np, r));
    end
    ret = ret + kronVec(g.markE0TE0T{nn, np}, RknTimesVal);
  end 
end 
end 
\end{lstlisting}
 
\begin{lstlisting}[title={%
%\begin{lstlisting}[title={%
\code{ret = assembleMatElemDphiPhiFuncDiscVec(g, refElemDphiPhiPhi, dataDisc1, dataDisc2)} assembles two matrices, each containing integrals of products of a~basis function with a (spatial) derivative of a basis function and with a~component of a~discontinuous coefficient function whose coefficients are specified in~\code{dataDisc1} and~\code{dataDisc2}, respectively. The matrices are returned in a~$2\times1$ \code{cell}~variable. This corresponds to the matrices~$\vecc{G}^m$, $m\in\{1,2\}$ according to Sec.~\ref{sec:assembly:globG}. The input argument~\code{refElemDphiPhiPhi} stores the local matrices~$\hat{\vecc{G}}$ (multidimensional array) as defined in~\eqref{eq:hatG} and can be computed by~\code{integrateRefElemDphiPhiPhi}~\cite{FESTUNG1}. The coefficients of the projection of the algebraic diffusion coefficient $d$ into the broken polynomial space are stored in the input arguments~\code{dataDisc1} and~\code{dataDisc2} and can be computed by~\code{projectFuncCont2dataDisc}.
}%
%]
%...
%\end{lstlisting}
]
function ret = assembleMatElemDphiPhiFuncDiscVec(g, refElemDphiPhiPhi, dataDisc1, dataDisc2)
[K, N] = size(dataDisc1);
ret = cell(2, 1);  ret{1} = sparse(K*N, K*N);  ret{2} = sparse(K*N, K*N);
for l = 1 : N
  ret{1} = ret{1} + kron(spdiags(dataDisc1(:,l) .* g.B(:,2,2), 0,K,K), refElemDphiPhiPhi(:,:,l,1)) ...
                  - kron(spdiags(dataDisc1(:,l) .* g.B(:,2,1), 0,K,K), refElemDphiPhiPhi(:,:,l,2));
  ret{2} = ret{2} - kron(spdiags(dataDisc2(:,l) .* g.B(:,1,2), 0,K,K), refElemDphiPhiPhi(:,:,l,1)) ...
                  + kron(spdiags(dataDisc2(:,l) .* g.B(:,1,1), 0,K,K), refElemDphiPhiPhi(:,:,l,2));
end 
end 
\end{lstlisting}
 
\begin{lstlisting}[title={%
%\begin{lstlisting}[title={%
\code{ret = assembleMatElemPhiDiscPhiTaylor(g, N)} assembles the matrix $\vecc{M}^\mathrm{DG,Taylor}$ according to Sec.~\ref{sec:taylor:basis}, which corresponds to the basis transformation matrix with one basis function from each, modal and Taylor basis. It is required for the transformation between modal and Taylor basis in the routines \code{projectDataDisc2DataTaylor} and \code{projectDataTaylor2DataDisc}.
}%
%]
%...
%\end{lstlisting}
]
function ret = assembleMatElemPhiDiscPhiTaylor(g, N)
global gPhi2D
p = (sqrt(8*N+1)-3)/2;  qOrd = max(2*p+1, 1);  [Q1,Q2,W] = quadRule2D(qOrd);
K = g.numT; ret = sparse(K*N, K*N);
for i = 1 : N
  for j = 1 : N
    intPhiIPhiJ =  ( repmat(gPhi2D{qOrd}(:, i)', [K 1]) .* phiTaylorRef(g, j, Q1, Q2) ) * W';
    ret = ret + sparse(i : N : K*N, j : N : K*N, 2 * g.areaT .* intPhiIPhiJ, K*N, K*N );
  end 
end 
end 
\end{lstlisting}
 
\begin{lstlisting}[title={%
%\begin{lstlisting}[title={%
\code{ret = assembleMatElemPhiTaylorPhiTaylor(g, N)} assembles the mass matrix in Taylor basis $\vecc{M}^\mathrm{Taylor}$.
}%
%]
%...
%\end{lstlisting}
]
function ret = assembleMatElemPhiTaylorPhiTaylor(g, N)
p = (sqrt(8*N+1)-3)/2;  qOrd = max(2*p, 1);  
[Q1, Q2, W] = quadRule2D(qOrd);
K = g.numT;
ret = sparse(K*N, K*N);
for i = 1 : N
  for j = 1 : N
    intPhiIPhiJ =  ( phiTaylorRef(g, i, Q1, Q2) .* phiTaylorRef(g, j, Q1, Q2) ) * W';
    ret = ret + sparse(i : N : K*N, j : N : K*N, 2 * g.areaT .* intPhiIPhiJ, K*N, K*N );
  end 
end 
end 
\end{lstlisting}
 
\begin{lstlisting}[title={%
%\begin{lstlisting}[title={%
\code{ret = assembleVecEdgePhiIntFuncContVal(g, markE0Tbdr, funcCont, valOnQuad, N)} assembles a~vector containing integrals of products of a~basis function with a~continuous function and a~given value that is provided in each quadrature point on each edge for all triangles. This corresponds to the contributions of Dirichlet boundaries~$\vec{K}_\mathrm{D}$ to the right-hand side of~\eqref{eq:spacediscretesystem} according to Sec.~\ref{sec:assembly:globKD}. \code{markE0Tbdr} marks the boundary edges on which the vector should be assembled, \code{funcCont} is a function handle and \code{valOnQuad} is the value provided in each quadrature point, as computed by \code{computeFuncContNuOnQuadEdge}.
}%
%]
%...
%\end{lstlisting}
]
function ret = assembleVecEdgePhiIntFuncContVal(g, markE0Tbdr, funcCont, valOnQuad, N)
global gPhi1D
K = g.numT;  p = (sqrt(8*N+1)-3)/2; qOrd = 2*p+1;  [Q, W] = quadRule1D(qOrd);
Q2X1 = @(X1,X2) g.B(:,1,1)*X1 + g.B(:,1,2)*X2 + g.coordV0T(:,1,1)*ones(size(X1));
Q2X2 = @(X1,X2) g.B(:,2,1)*X1 + g.B(:,2,2)*X2 + g.coordV0T(:,1,2)*ones(size(X1));
ret = zeros(K, N);
for n = 1 : 3
  [Q1, Q2] = gammaMap(n, Q);  funcOnQuad = funcCont(Q2X1(Q1, Q2), Q2X2(Q1, Q2));
  Kkn = markE0Tbdr(:, n) .* g.areaE0T(:,n);
  for i = 1 : N
    integral = (funcOnQuad .* squeeze((valOnQuad(:, n, :) < 0) .* valOnQuad(:, n, :))) * ( W' .* gPhi1D{qOrd}(:,i,n));
    ret(:,i) = ret(:,i) + Kkn .* integral;
  end 
end 
ret = reshape(ret',K*N,1);
end 
\end{lstlisting}

\begin{lstlisting}[title={%
%\begin{lstlisting}[title={%
\code{ret = computeFuncContNuOnQuadEdge(g, funcCont1, funcCont2, qOrd)} 
  assembles a~three-dimensional array with the normal velocity~$\vec{u}\cdot\vec{\nu}_{ki}$ evaluated in all quadrature points of all edges of each triangle.
}%
%]
%...
%\end{lstlisting}
]
function ret = computeFuncContNuOnQuadEdge(g, funcCont1, funcCont2, qOrd)
K = g.numT;  [Q, W] = quadRule1D(qOrd);
ret = zeros(K, 3, length(W));
for n = 1 : 3
  [Q1, Q2] = gammaMap(n, Q);
  ret(:,n,:) = bsxfun(@times,g.nuE0T(:,n,1),funcCont1(g.mapRef2Phy(1,Q1,Q2),g.mapRef2Phy(2,Q1,Q2)))+...
               bsxfun(@times,g.nuE0T(:,n,2),funcCont2(g.mapRef2Phy(1,Q1,Q2),g.mapRef2Phy(2,Q1,Q2)));
end 
end
\end{lstlisting}

\begin{lstlisting}[title={%
%\begin{lstlisting}[title={%
\code{valV0T = computeFuncContV0T(g, funcCont)} 
assembles a~matrix containing the function \code{funcCont} evaluated in each node of each triangle.
}%
%]
%...
%\end{lstlisting}
]
function valV0T = computeFuncContV0T(g, funcCont)
valV0T = zeros(g.numT,3);
for n = 1 : 3
    valV0T(:, n) = funcCont(g.coordV0T(:, n, 1), g.coordV0T(:, n, 2));
end
end
\end{lstlisting}

\begin{lstlisting}[title={%
%\begin{lstlisting}[title={%
\code{ret = computeFuncDiscAtPoints(funcDisc, phiAtPoints)} 
assembles a~matrix containing the values of a~discrete function with representation matrix stored in \code{funcDisc} evaluated in all points, for which the evaluated basis functions are given in \code{phiAtPoints}.
}%
%]
%...
%\end{lstlisting}
]
function ret = computeFuncDiscAtPoints(funcDisc, phiAtPoints)
nPoints = size(phiAtPoints, 2); K = size(funcDisc, 1);
ret = zeros(K, nPoints);
for i = 1 : nPoints
  ret(:, i) = sum(funcDisc .* squeeze(phiAtPoints(:, i, :)), 2);
end 
end 
\end{lstlisting}

\begin{lstlisting}[title={%
%\begin{lstlisting}[title={%
\code{minMaxV0T = computeMinMaxV0TElementPatch(g, valCentroid, markV0TbdrD, dataV0T)}
determines a~matrix with bounds~$c_{ki}^\mathrm{max}$,~$c_{ki}^\mathrm{min}$ of Eq.~\eqref{eq:slope-lim:lin:bounds} for each vertex of each triangle, as required by \code{computeVertexBasedLimiter}.
}%
%]
%...
%\end{lstlisting}
]
function minMaxV0T = computeMinMaxV0TElementPatch(g, valCentroid, markV0TbdrD, dataV0T)
minMaxV0T = cell(2,1); 
minMaxV0T{1} = zeros(g.numT, 3); minMaxV0T{2} = zeros(g.numT, 3);
shiftCentroidPos = abs(min(valCentroid)) + 1;
shiftCentroidNeg = abs(max(valCentroid)) + 1;
valCentroidPos = valCentroid + shiftCentroidPos;
valCentroidNeg = valCentroid - shiftCentroidNeg;
valD = NaN(g.numT, 3); valD(markV0TbdrD) = dataV0T(markV0TbdrD);
for i = 1 : 3
  markNbV0T = g.markV0TV0T{i, 1} | g.markV0TV0T{i, 2} | g.markV0TV0T{i, 3}; 
  if exist('OCTAVE_VERSION','builtin')
    markNbV0T = markNbV0T + 0 * speye(size(markNbV0T, 1), size(markNbV0T, 2));
  end
  minMaxV0T{1}(:,i)=min(min(bsxfun(@times,markNbV0T,valCentroidNeg'),[],2)+shiftCentroidNeg,valD(:,i));
  minMaxV0T{2}(:,i)=max(max(bsxfun(@times,markNbV0T,valCentroidPos'),[],2)-shiftCentroidPos,valD(:,i));
end 
end 
\end{lstlisting}

\begin{lstlisting}[title={%
%\begin{lstlisting}[title={%
\code{computeTaylorBasesV0T(g, N)}
evaluates the Taylor basis functions~$\Phi_{kj}$ in all vertices of all triangles~$\vec{x}_{ki}$ and stores them in a~global multidimensional array.
}%
%]
%...
%\end{lstlisting}
]
function computeTaylorBasesV0T(g, N)
global gPhiTaylorV0T
gPhiTaylorV0T = zeros(g.numT, 3, N);
for n = 1 : 3
  for i = 1 : N
    gPhiTaylorV0T(:, n, i) = phiTaylorPhy(g, i, g.coordV0T(:, n, 1), g.coordV0T(:, n, 2));
  end
end
end
\end{lstlisting}

\begin{lstlisting}[title={%
%\begin{lstlisting}[title={%
\code{alphaE = computeVertexBasedLimiter(g, valCentroid, valV0T, markV0TbdrD, dataV0T)}
computes a~vector with correction factors~$\alpha_{ke}$ (cf.~Eq.~\eqref{eq:slope-lim:lin:corr}) for all elements. Centroid values~$c_{kc}$ are given in \code{valCentroid}, values of the unconstrained reconstruction~$c_{ki}$ are specified in \code{valV0T} and \code{markV0TbdrD} is a~logical matrix marking all Dirichlet boundary nodes for which boundary data is given in \code{dataV0T}.
}%
%]
%...
%\end{lstlisting}
]
function alphaE = computeVertexBasedLimiter(g, valCentroid, valV0T, markV0TbdrD, dataV0T)
minMaxV0T = computeMinMaxV0TElementPatch(g, valCentroid, markV0TbdrD, dataV0T);
diffV0TCentroid = valV0T - repmat(valCentroid, [1 3]);
diffMinCentroid = minMaxV0T{1} - repmat(valCentroid, [1 3]);
diffMaxCentroid = minMaxV0T{2} - repmat(valCentroid, [1 3]);
tol = 1.e-8;
markNeg = diffV0TCentroid < diffMinCentroid + tol;
markPos = diffV0TCentroid > diffMaxCentroid - tol;
alphaEV0T = ones(g.numT,3);
alphaEV0T(markNeg) = max(0, min(1, diffMinCentroid(markNeg) ./ (diffV0TCentroid(markNeg) - tol) ) );
alphaEV0T(markPos) = max(0, min(1, diffMaxCentroid(markPos) ./ (diffV0TCentroid(markPos) + tol) ) );
alphaE = min(alphaEV0T,[],2);
end 
\end{lstlisting}

\begin{lstlisting}[title={%
%\begin{lstlisting}[title={%
\code{ret = integrateRefEdgePhiIntPhiExtPerQuad(N)} computes a~multidimensional array of functionals in the quadrature points on the edges of the reference triangle~$\hat{T}$ that consist of all permutations of two basis functions of which one belongs to a~neighboring element that is transformed using~$\mapEE$. This corresponds to the local matrix~$\hat{\vecc{R}}^\mathrm{offdiag}$ as given in~\eqref{eq:hatRoffdiag}.
}%
%]
%...
%\end{lstlisting}
]
function ret = integrateRefEdgePhiIntPhiExtPerQuad(N)
global gPhi1D gThetaPhi1D
p = (sqrt(8*N+1)-3)/2;  qOrd = 2*p+1;  [~, W] = quadRule1D(qOrd);
ret = zeros(N,N,3,3,length(W)); 
for nn = 1 : 3 
  for np = 1 : 3
      for i = 1 : N
        for j = 1 : N
          ret(i, j, nn, np, :) = gPhi1D{qOrd}(:,i,nn) .* gThetaPhi1D{qOrd}(:,j,nn,np);
        end 
      end 
  end 
end 
end
\end{lstlisting}

\begin{lstlisting}[title={%
%\begin{lstlisting}[title={%
\code{ret = integrateRefEdgePhiIntPhiIntPerQuad(N)} computes a~multidimensional array of functionals in the quadrature points on the edges of the reference triangle~$\hat{T}$ that consist of all permutations of two basis functions. This corresponds to the local matrix~$\hat{\vecc{R}}^\mathrm{diag}$ as given in~\eqref{eq:hatRdiag}.
}%
%]
%...
%\end{lstlisting}
]
function ret = integrateRefEdgePhiIntPhiIntPerQuad(N)
global gPhi1D
p = (sqrt(8*N+1)-3)/2;  qOrd = max(2*p+1,1);  [~, W] = quadRule1D(qOrd);
ret = zeros(N, N, 3, length(W)); 
for n = 1 : 3 
    for i = 1 : N
      for j = 1 : N
        ret(i,j,n,:) = gPhi1D{qOrd}(:,i,n).* gPhi1D{qOrd}(:,j,n);
      end 
    end 
end 
end 
\end{lstlisting}

\begin{lstlisting}[title={%
%\begin{lstlisting}[title={%
\code{K = kronVec(A,B)} computes the result of $\IR^{m_b \times n_a n_b} \ni \vecc{K} = \vecc{A} \otimes_\mathrm{V} \vecc{B}$ as given in Eq.~\eqref{eq:kronVec}.
}%
%]
%...
%\end{lstlisting}
]
function K = kronVec(A, B)
[ma,na] = size(A);
[mb,nb] = size(B);
mc = mb / ma;
if ~issparse(A) && ~issparse(B) 
  A = reshape(A, [1 ma 1 na]);
  B = reshape(B, [mc ma nb 1]);
  K = reshape(bsxfun(@times, A, B), [mb na*nb]);
else
  [i2, j2, v2] = find(kron(A, ones(mc, 1)));
  ik = repmat(i2, [1 nb]);
  jk = bsxfun(@plus, nb * (j2 - 1), 1 : nb);
  sk = bsxfun(@times, v2, B(i2, :));
  K = sparse(ik, jk, sk, mb, na * nb);
end 
end 
\end{lstlisting}
\begin{lstlisting}[title={%
%\begin{lstlisting}[title={%
\code{mainAdvection.m} {This is the main script to solve~\eqref{eq:model} which can be used as a template for further modifications. Modifiable parameters are found in Lines~5--16, the problem data (initial condition, velocity, right-hand side and boundary data) is specified in Lines~35--42.}
}%
%]
%...
%\end{lstlisting}
]
function mainAdvection()
more off 
tic 
hmax        = 2^-6;              
p           = 2;                 
ordRK       = min(p+1,3);        
tEnd        = 2*pi;              
numSteps    = 3142;              
isVisGrid   = false;             
isVisSol    = true;              
isSlopeLim  = true;              
typeSlopeLim = 'hierarch_vert';  
outputFrequency = 100;           
outputBasename  = ['solution_' typeSlopeLim]; 
outputTypes     = cellstr(['vtk';'tec']);
diary([outputBasename '.log'])
assert(p >= 0 && p <= 4        , 'Polynomial order must be zero to four.'     )
assert(ordRK >= 1 && ordRK <= 3, 'Order of Runge Kutta must be zero to three.')
assert(hmax > 0                , 'Maximum edge length must be positive.'      )
assert(numSteps > 0            , 'Number of time steps must be positive.'     )
assert(~isSlopeLim || p > 0    , 'Slope limiting only available for p > 0.'   )
g = domainSquare(hmax); 
if isVisGrid,  visualizeGrid(g);  end
K           = g.numT;                      
N           = nchoosek(p + 2, p);          
tau         = tEnd/numSteps;               
markE0Tint  = g.idE0T == 0;                
markE0TbdrD = ~markE0Tint;                 
markV0TbdrD = ismember(g.V0T, g.V0E(g.E0T(markE0TbdrD),:)); 
G = @(x1, x2, x1_0, x2_0) (1/0.15) * sqrt((x1-x1_0).^2 + (x2-x2_0).^2);
c0Cont = @(x1, x2) ((x1-0.5).^2+(x2-0.75).^2 <= 0.0225 & (x1<=0.475|x1>=0.525|x2>=0.85)) + ...
                    (1-G(x1, x2, 0.5, 0.25)) .* ((x1 - 0.5).^2 + (x2 - 0.25).^2 <= 0.0225) + ...
                   0.25*(1+cos(pi*G(x1, x2, 0.25, 0.5))).*((x1 - 0.25).^2 + (x2 - 0.5).^2 <= 0.0225);
fCont = @(t,x1,x2) zeros(size(x1));
u1Cont = @(t,x1,x2) 0.5 - x2;
u2Cont = @(t,x1,x2) x1 - 0.5;
cDCont = @(t,x1,x2) zeros(size(x1));
computeBasesOnQuad(N);
if isSlopeLim, computeTaylorBasesV0T(g, N); end
hatM              = integrateRefElemPhiPhi(N);
hatG              = integrateRefElemDphiPhiPhi(N);
hatRdiagOnQuad    = integrateRefEdgePhiIntPhiIntPerQuad(N);
hatRoffdiagOnQuad = integrateRefEdgePhiIntPhiExtPerQuad(N);
globM  = assembleMatElemPhiPhi(g, hatM);
if isSlopeLim
  globMTaylor     = assembleMatElemPhiTaylorPhiTaylor(g, N);
  globMDiscTaylor = assembleMatElemPhiDiscPhiTaylor(g, N);
  globMCorr       = spdiags(1./diag(globMTaylor), 0, K*N, K*N) * globMTaylor;
end 
cDisc = projectFuncCont2DataDisc(g, c0Cont, 2*p+1, hatM);
if isSlopeLim
  cDV0T = computeFuncContV0T(g, @(x1, x2) cDCont(0, x1, x2));
  cDisc = applySlopeLimiterDisc(g, cDisc, markV0TbdrD, cDV0T, globM, globMDiscTaylor, typeSlopeLim);
end 
fprintf('L2 error w.r.t. the initial condition: 
if isVisSol
  cLagrange = projectDataDisc2DataLagr(cDisc);
  visualizeDataLagr(g, cLagrange, 'u_h', outputBasename, 0, outputTypes)
end
fprintf('Starting time integration from 0 to 
for nStep = 1 : numSteps
  [t, omega] = rungeKuttaSSP(ordRK, tau, (nStep - 1) * tau);
  cDiscRK = cell(length(omega)+1, 1); cDiscRK{1} = reshape(cDisc', [K*N 1]);
  for rkStep = 1 : length(omega)
    fDisc  = projectFuncCont2DataDisc(g, @(x1,x2) fCont(t(rkStep),x1,x2),  2*p, hatM);
    u1Disc = projectFuncCont2DataDisc(g, @(x1,x2) u1Cont(t(rkStep),x1,x2), 2*p, hatM);
    u2Disc = projectFuncCont2DataDisc(g, @(x1,x2) u2Cont(t(rkStep),x1,x2), 2*p, hatM);
    vNormalOnQuadEdge = computeFuncContNuOnQuadEdge(g, @(x1,x2) u1Cont(t(rkStep),x1,x2), @(x1,x2) u2Cont(t(rkStep),x1,x2), 2*p+1); 
    globG = assembleMatElemDphiPhiFuncDiscVec(g, hatG, u1Disc, u2Disc);
    globR = assembleMatEdgePhiPhiValUpwind(g, hatRdiagOnQuad, hatRoffdiagOnQuad, vNormalOnQuadEdge);
    globKD = assembleVecEdgePhiIntFuncContVal(g, markE0TbdrD, @(x1,x2) cDCont(t(rkStep),x1,x2), vNormalOnQuadEdge, N);
    globL = globM * reshape(fDisc', K*N, 1);
    sysA = -globG{1} - globG{2} + globR;
    sysV = globL - globKD;
    cDiscDot = globM \ (sysV - sysA * cDiscRK{rkStep});
    if isSlopeLim
      cDiscDotTaylor = projectDataDisc2DataTaylor(reshape(cDiscDot, [N K])', globM, globMDiscTaylor);
      cDiscDotTaylorLim = applySlopeLimiterTaylor(g,cDiscDotTaylor,markV0TbdrD,NaN(K,3),typeSlopeLim);
      cDiscDotTaylor = reshape(cDiscDotTaylorLim', [K*N 1]) + globMCorr * reshape((cDiscDotTaylor - cDiscDotTaylorLim)', [K*N 1]);
      cDiscDot = reshape(projectDataTaylor2DataDisc(reshape(cDiscDotTaylor, [N K])', globM, globMDiscTaylor)', [K*N 1]);
    end
    cDiscRK{rkStep + 1} = omega(rkStep) * cDiscRK{1} + (1 - omega(rkStep)) * (cDiscRK{rkStep} + tau * cDiscDot);
    if isSlopeLim
      cDV0T = computeFuncContV0T(g, @(x1, x2) cDCont(t(rkStep), x1, x2));
      cDiscRK{rkStep + 1} = reshape(applySlopeLimiterDisc(g, reshape(cDiscRK{rkStep + 1}, [N K])', markV0TbdrD, cDV0T, globM, globMDiscTaylor, typeSlopeLim)', [K*N 1]);
    end 
  end 
  cDisc = reshape(cDiscRK{end}, N, K)';
  if isVisSol && mod(nStep, outputFrequency) == 0
    cLagrange = projectDataDisc2DataLagr(cDisc);
    visualizeDataLagr(g, cLagrange, 'u_h', outputBasename, nStep, outputTypes);
  end
end 
if isVisSol
  cLagrange = projectDataDisc2DataLagr(cDisc);
  visualizeDataLagr(g, cLagrange, 'u_h', outputBasename, nStep, outputTypes);
end
fprintf('L2 error w.r.t. the initial condition: 
fprintf('Total computation time: 
diary off
end 
\end{lstlisting}

\begin{lstlisting}[title={%
%\begin{lstlisting}[title={%
\code{ind = mult2ind(a)} computes the linear index $I(\vec{a})$ corresponding to a two-dimensional multi-index $\vec{a}$ as defined in~\eqref{eq:multiindex:lin}.
}%
%]
%...
%\end{lstlisting}
]
function ind = mult2ind(a)
  p = sum(a, 2);  N = p .* (p + 1) / 2;
  ind = N + 1 + a(:, 2);
end
\end{lstlisting}

\begin{lstlisting}[title={%
%\begin{lstlisting}[title={%
\code{mult = multiindex(p)} computes all two-dimensional multi-indices involved in the representation of a~polynomial solution of degree~$p$ and returns them in a~$N\times2$ array.
}%
%]
%...
%\end{lstlisting}
]
function mult = multiindex(p)
mult = zeros(p * (p+1) / 2, 2);
mult(1,:) = [0, 0];
for ord = 1 : p
  offset = ord * (ord+1) / 2;
  for i = 1 : ord + 1
    mult(offset + i, :) = mult(1, :) + [ord - i + 1, i - 1];
  end
end
end
\end{lstlisting}

\begin{lstlisting}[title={%
%\begin{lstlisting}[title={%
\code{ret = phiTaylorPhy(g, i, X1, X2)} {evaluates the $i$th basis function~$\phi_i$ on each triangle~$T \in \setT_h$ (cf.~Sec.~\ref{sec:taylor:basis}) at points specified by a~list of $n$ $x^1$~coordinates~\code{X1} $\in\IR^{K\times n}$ and $x^2$~coordinates~\code{X2} $\in\IR^{K\times n}$.}
}%
%]
%...
%\end{lstlisting}
]
function ret = phiTaylorPhy(g, i, X1, X2)
qOrd = ceil((sqrt(8*i+1)-3)/2);
[Q1, Q2, W] = quadRule2D(qOrd);
Q2X1 = @(X1, X2) g.B(:, 1, 1) * X1 + g.B(:, 1, 2) * X2 + g.coordV0T(:, 1, 1) * ones(size(X1));
Q2X2 = @(X1, X2) g.B(:, 2, 1) * X1 + g.B(:, 2, 2) * X2 + g.coordV0T(:, 1, 2) * ones(size(X1));
R = length(W); K = g.numT; numP = size(X1, 2);
dX1 = repmat(2 ./ (max(g.coordV0T(:,:,1),[],2) - min(g.coordV0T(:,:,1),[],2)), [1 numP]);
dX2 = repmat(2 ./ (max(g.coordV0T(:,:,2),[],2) - min(g.coordV0T(:,:,2),[],2)), [1 numP]);
switch i
  case 1 
    ret = ones(K, numP);
  case 2 
    ret = (X1-repmat(g.baryT(:,1), [1 numP]).*dX1;
  case 3 
    ret = (X2-repmat(g.baryT(:,2), [1 numP]).*dX2;
  case 4 
    ret = ( 0.5*(X1-repmat(g.baryT(:,1), [1 numP]) ).^2-...
          repmat( ( Q2X1(Q1,Q2)-repmat(g.baryT(:,1), [1 R]) ).^2 * W', [1 numP])).*(dX1.*dX1);
  case 5 
    ret = ( ( X1-repmat(g.baryT(:,1), [1 numP]) ) .* ( X2-repmat(g.baryT(:,2), [1 numP]))-...
                   repmat( 2 * ( ( Q2X1(Q1,Q2)-repmat(g.baryT(:,1), [1 R]) ) .* ...
                          ( Q2X2(Q1,Q2)-repmat(g.baryT(:,2), [1 R]) ) ) * W', [1 numP])).*(dX1.*dX2);
  case 6 
    ret = ( 0.5*(X2-repmat(g.baryT(:,2), [1 numP]) ).^2-...
          repmat( ( Q2X2(Q1,Q2)-repmat(g.baryT(:,2), [1 R]) ).^2 * W', [1 numP])).*(dX2.*dX2);
  case 7 
    ret = ( ( X1-repmat(g.baryT(:,1), [1 numP]) ).^3 / 6-...
          repmat( ( Q2X1(Q1,Q2)-repmat(g.baryT(:,1), [1 R]) ).^3 * W' / 3, [1 numP])).*dX1.^3;
  case 8 
    ret = ( 0.5*(X1-repmat(g.baryT(:,1), [1 numP]) ).^2 .* ( X2-repmat(g.baryT(:,2), [1 numP]))-...
          repmat( ( ( Q2X1(Q1,Q2)-repmat(g.baryT(:,1), [1 R]) ).^2 .* ...
                    ( Q2X2(Q1,Q2)-repmat(g.baryT(:,2), [1 R]) ) ) * W', [1 numP])).*(dX1.^2.*dX2);
  case 9 
    ret = ( 0.5*(X1-repmat(g.baryT(:,1), [1 numP]) ) .* ( X2-repmat(g.baryT(:,2), [1 numP]) ).^2-...
          repmat( ( ( Q2X1(Q1,Q2)-repmat(g.baryT(:,1), [1 R]) ) .* ...
                    ( Q2X2(Q1,Q2)-repmat(g.baryT(:,2), [1 R]) ).^2 ) * W', [1 numP])).*(dX1.*dX2.^2);
  case 10 
    ret = ( ( X2-repmat(g.baryT(:,2), [1 numP]) ).^3 / 6-...
          repmat( ( Q2X2(Q1,Q2)-repmat(g.baryT(:,2), [1 R]) ).^3 * W' / 3, [1 numP])).*dX2.^3;
  case 11 
    ret = ( ( X1-repmat(g.baryT(:,1), [1 numP]) ).^4 / 24-...
          repmat( ( Q2X1(Q1,Q2)-repmat(g.baryT(:,1), [1 R]) ).^4 * W' / 12, [1 numP])).*dX1.^4;
  case 12 
    ret = ( ( X1-repmat(g.baryT(:,1), [1 numP]) ).^3 .* ( X2-repmat(g.baryT(:,2), [1 numP]) ) / 6-...
          repmat( ( ( Q2X1(Q1,Q2)-repmat(g.baryT(:,1), [1 R]) ).^3 .* ...
                    ( Q2X2(Q1,Q2)-repmat(g.baryT(:,2), [1 R]) ) ) * W' / 3, [1 numP])).*(dX1.^3.*dX2);
  case 13 
    ret = ( ( X1-repmat(g.baryT(:,1), [1 numP]) ).^2 .* ( X2-repmat(g.baryT(:,2), [1 numP])).^2 / 4-...
          repmat( ( (Q2X1(Q1,Q2)-repmat(g.baryT(:,1), [1 R])).^2 .* ...
                    (Q2X2(Q1,Q2)-repmat(g.baryT(:,2), [1 R])).^2 ) * W'/2, [1 numP])).*(dX1.*dX2).^2;
  case 14 
    ret = ( ( X1-repmat(g.baryT(:,1), [1 numP]) ) .* ( X2-repmat(g.baryT(:,2), [1 numP])).^3 / 6-...
          repmat( ( (Q2X1(Q1,Q2)-repmat(g.baryT(:,1), [1 R]) ) .* ...
                    (Q2X2(Q1,Q2)-repmat(g.baryT(:,2), [1 R])).^3) * W'/3, [1 numP])).*(dX1.*dX2.^3);
  case 15 
    ret = ( ( X2-repmat(g.baryT(:,2), [1 numP]) ).^4 / 24-...
          repmat( ( Q2X2(Q1,Q2)-repmat(g.baryT(:,2), [1 R]) ).^4 * W' / 12, [1 numP])).*dX2.^4;
end 
end

\end{lstlisting}

\begin{lstlisting}[title={%
%\begin{lstlisting}[title={%
\code{ret = phiTaylorRef(g, i, hatX1, hatX2)} {evaluates the $i$th basis function~$\phi_i$ on each triangle~$T \in \setT_h$ (cf.~Sec.~\ref{sec:taylor:basis}) at points specified by a~list of $\hat{x}^1$~coordinates~\code{hatX1} and $\hat{x}^2$~coordinates~\code{hatX2}.}
}%
%]
%...
%\end{lstlisting}
]
function ret = phiTaylorRef(g, i, hatX1, hatX2)
Q2X1 = @(X1, X2) g.B(:, 1, 1) * X1 + g.B(:, 1, 2) * X2 + g.coordV0T(:, 1, 1) * ones(size(X1));
Q2X2 = @(X1, X2) g.B(:, 2, 1) * X1 + g.B(:, 2, 2) * X2 + g.coordV0T(:, 1, 2) * ones(size(X1));
ret = phiTaylorPhy(g, i, Q2X1(hatX1, hatX2), Q2X2(hatX1, hatX2));
end

\end{lstlisting}

\begin{lstlisting}[title={%
%\begin{lstlisting}[title={%
\code{dataTaylor = projectDataDisc2DataTaylor(dataDisc, globMDisc, globMDiscTaylor)}
converts the representation matrix in the DG\,/\,modal basis to the respective representation matrix in a~Taylor basis, both of size~\mbox{$K\times N$}. It solves Eq.~\eqref{eq:taylor:trafo} for $\vecc{C}^\mathrm{Taylor}$.
}%
%]
%...
%\end{lstlisting}
]
function dataTaylor = projectDataDisc2DataTaylor(dataDisc, globMDisc, globMDiscTaylor)
[K, N] = size(dataDisc);
dataTaylor = reshape(globMDiscTaylor \ ( globMDisc * reshape(dataDisc', [K*N 1]) ), [N K])';
end 
\end{lstlisting}

\begin{lstlisting}[title={%
%\begin{lstlisting}[title={%
\code{dataDisc = projectDataTaylor2DataDisc(dataTaylor, globMDisc, globMDiscTaylor)}
converts the representation matrix in the Taylor basis to the respective representation matrix in a~ DG\,/\,modal basis, both of size~\mbox{$K\times N$}. It solves Eq.~\eqref{eq:taylor:trafo} for $\vecc{C}^\mathrm{DG}$.
}%
%]
%...
%\end{lstlisting}
]
function dataDisc = projectDataTaylor2DataDisc(dataTaylor, globMDisc, globMDiscTaylor)
[K, N] = size(dataTaylor);
dataDisc = reshape(globMDisc \ (globMDiscTaylor * reshape(dataTaylor', [K*N 1]) ), [N K])';
end 
\end{lstlisting}

\begin{lstlisting}[title={%
%\begin{lstlisting}[title={%
\code{[t, omega] = rungeKuttaSSP(ord, tau, nStep)}
provides a~list of time levels $t^{(i)} = t^n + \delta_i \Delta t^n$ and weights~$\omega_i$ for the $n$-th time step \code{nStep} with time step size \code{tau} according to Sec.~\ref{sec:timediscretization}. The order of the Runge-Kutta method is given as parameter~\code{ord}.
}%
%]
%...
%\end{lstlisting}
]
function [t, omega] = rungeKuttaSSP(ord, tau, t0)
switch ord
  case 1
    omega = 0;             t = t0;
  case 2
    omega = [0, 0.5];      t = t0 + [0, 1] * tau;
  case 3
    omega = [0, 3/4, 1/3]; t = t0 + [0, 1, 0.5] * tau;
end
end 
\end{lstlisting}

\section{Conclusion and Outlook}\label{sec:conclusion}

The second installment in the present paper series on implementing a~\MatOct~toolbox 
introduced performance optimized techniques for dealing with linear advection 
operators, higher order Runge--Kutta time discretizations, and a~range of slope limiters
designed to support general order DG~discretizations.
Our future work plans include nonlinear advection operators and coupled systems of PDEs as
well as multi-physics applications with corresponding coupling mechanisms.

\subsection*{\bf Acknowledgments}
The work of B.~Reuter was supported by the German Research Foundation (DFG) under grant AI 117/1-1.

\section*{Index of notation}
\noindent
\begin{footnotesize}
\begin{tabularx}{\linewidth}{@{}lX@{}}\toprule
\textbf{Symbol}    & \textbf{Definition}\\\midrule
$\overline{\,\cdot\,}$  & Integral mean, $\overline{v} \coloneqq \frac{1}{|T|}\int_T v(\vec{x})\,\dd\vec{x}$, where $v:T\rightarrow \IR$.\\
$\diag(\vecc{A},\vecc{B})$ & $\coloneqq \begin{bmatrix}\vecc{A}&\quad\\\quad&\vecc{B}\end{bmatrix}$, block-diagonal matrix with blocks~$\vecc{A}$, $\vecc{B}$.\\
$\card{\mathcal{M}}$ & Cardinality of a~set~$\mathcal{M}$.\\
$\vec{a}\cdot\vec{b}$& $\coloneqq \sum_{m=1}^2a_mb_m$, Euclidean scalar product in~$\IR^2$.\\
$\grad$             & $\coloneqq \transpose{[\partial_{x^1}, \partial_{x^2}]}$, spatial gradient in the physical domain~$\Omega$.\\
$\circ$             & Composition of functions or Hadamard product.\\
$\otimes$           & Kronecker product.\\
$c$                 & Concentration (scalar-valued unknown).\\
$c^0$               & Concentration prescribed at initial time~$t=0$.\\
$c_\mathrm{D}$      & Concentration prescribed on the inflow boundary.\\
$\vec{C}$           & $\in\IR^{KN}$, representation vector of~$c_h\in\IP_p(\setT_h)$ with respect to $\{\varphi_{kj}\}$.\\
$\vec{C}^\mathrm{Taylor}$ & $\in\IR^{KN}$, representation vector of~$c_h\in\IP_p(\setT_h)$ with respect to $\{\phi_{kj}\}$.\\
$\delta_\mathrm{[condition]}$ & $\coloneqq \{1~\text{if condition is true, 0~otherwise}\}$, Kronecker delta.\\
$\vec{e}_m$         & $m$th unit vector.\\
$E_{kn}$, $\hat{E}_n$          & $n$th edge of the physical triangle~$T_k$, $n$th edge of the reference triangle~$\hat{T}$.\\
$\setV$,\; $\setE$,\; $\setT$   & Sets of vertices, edges, and triangles.\\
$\setE_{\Omega}$,\; $\setE_{\partial\Omega}$ & Set of interior edges, set of boundary edges.\\
$f$                 & Source\,/\,sink (scalar-valued coefficient function).\\ 
$\vec{F}_k$         & Affine mapping from $\hat{T}$ to $T_k$.\\
$h$                 & Mesh fineness of $\setT_h$.\\
$h_T$               & $\coloneqq \diam(T)$, diameter of triangle~$T\in\setT_h$.\\
$J$                 & $\coloneqq (0,t_\mathrm{end})$, open time interval.\\
$K$                 & $\coloneqq \card{\setT_h}$, number of triangles.\\
$\vec{\nu}$         & Unit normal on~$\partial \Omega$ pointing outward of~$\Omega$.\\
$\vec{\nu}_{T}$     & Unit normal on~$\partial T$ pointing outward of~$T$.\\
$\vec{\nu}_k$       & $\coloneqq\vec{\nu}_{T_k}$.\\
$N=N_p$             & $\coloneqq (p+1)(p+2)/2$, number of local degrees of freedom of~$\IP_p(T)$.\\
$\omega_r$          & Quadrature weight associated with~$\hat{\vec{q}}_r$.\\
$\Omega$,\; $\partial\Omega$          & spatial domain in two dimensions, boundary of $\Omega$.\\
$\partial\Omega_\mathrm{in}$,\; $\partial\Omega_\mathrm{out}$& inflow- and outflow boundaries,  $\partial\Omega = \partial\Omega_\mathrm{in}\cup\partial\Omega_\mathrm{out}$.\\
$p$                 & $= (\sqrt{8N+1}-3)/2$, polynomial degree.\\
$\vphi_{ki}$,\;  $\hat{\vphi}_i$      & $i$th hierarchical basis function on~$T_k$, $i$th hierarchical basis function on~$\hat{T}$.\\
$\phi_{ki}$         & $i$th Taylor basis function on~$T_k$.\\
$\IP_p(T)$          & Space of polynomials on~$T\in\setT_h$ of degree at most~$p$.\\
$\IP_p(\setT_h)$    & $\coloneqq \{ w_h:\overline{\Omega}\rightarrow \IR\,;\forall T\in\setT_h,\,   {w_h}|_T\in\IP_p(T)\}$.\\
$\vecc{\Phi}$       & Slope limiting operator with respect to $\{\varphi_{kj}\}$.\\
$\vecc{\Phi}^\mathrm{Taylor}$ & Slope limiting operator with respect to $\{\phi_{kj}\}$.\\
$\hat{\vec{q}}_r$   & $r$th quadrature point in~$\hat{T}$.\\
$R$                 & Number of quadrature points.\\
$\IR^+$,\; $\IR_0^+$ & Set of (strictly) positive real numbers, set of non-negative real numbers.\\
\end{tabularx}

\begin{tabularx}{\linewidth}{@{}lX@{}}
$t$                 & Time variable.\\
$t^n$               & $n$th time level.\\
$t_\mathrm{end}$    & End time.\\
$\mapEE_{n^-n^+}$   & Mapping from $\hat{E}_{n^-}$ to $\hat{E}_{n^+}$.\\
$\Delta t^n$        & $\coloneqq t^{n+1} - t^n$, time step size.\\
%
$T_k$,\; $\partial T_k$               & $k$th physical triangle, boundary of~$T_k$.\\
$\hat{T}$           & Bi-unit reference triangle.\\
$\vec{u}$           & Velocity (vector-valued coefficient function).\\
$\vec{x}$           & $=\transpose{[x^1,x^2]}$, space variable in the physical domain~$\Omega$.\\
$\hat{\vec{x}}$     & $=\transpose{[\hat{x}^1, \hat{x}^2]}$, space variable in the reference triangle~$\hat{T}$.\\
$\vec{x}_{k\mathrm{c}}$& Centroid of the element $T_k\in\setT_h$.\\
$\vec{x}_{ki}$      & $i$th vertex of the physical triangle~$T_k$.\\
\bottomrule
\end{tabularx}
\end{footnotesize}

\bibliography{FESTUNG}

\begin{thebibliography}{10}
\expandafter\ifx\csname url\endcsname\relax
  \def\url#1{\texttt{#1}}\fi
\expandafter\ifx\csname urlprefix\endcsname\relax\def\urlprefix{URL }\fi
\expandafter\ifx\csname href\endcsname\relax
  \def\href#1#2{#2} \def\path#1{#1}\fi

\bibitem{FESTUNG}
F.~Frank, B.~Reuter, V.~Aizinger,
  \href{http://www.math.fau.de/FESTUNG}{{FESTUNG}---{T}he {F}inite {E}lement
  {S}imulation {T}oolbox for {UN}structured {G}rids} (2016).
\newblock \href {http://dx.doi.org/10.5281/zenodo.46069}
  {\path{doi:10.5281/zenodo.46069}}.
\newline\urlprefix\url{http://www.math.fau.de/FESTUNG}

\bibitem{FESTUNGGithub}
F.~Frank, B.~Reuter, \href{https://github.com/FESTUNG}{{FESTUNG}: {T}he
  {F}inite {E}lement {S}imulation {T}oolbox for {UN}structured {G}rids} (2016).
\newline\urlprefix\url{https://github.com/FESTUNG}

\bibitem{FESTUNG1}
F.~Frank, B.~Reuter, V.~Aizinger, P.~Knabner, {FESTUNG}:
  {A~MATLAB\,/\,GNU~Octave}~toolbox for the discontinuous {G}alerkin method,
  {P}art {I}: {D}iffusion operator, Computers \& Mathematics with Applications
  70~(1) (2015) 11--46.
\newblock \href {http://dx.doi.org/10.1016/j.camwa.2015.04.013}
  {\path{doi:10.1016/j.camwa.2015.04.013}}.

\bibitem{CockburnShu1998}
B.~Cockburn, C.~Shu, The local discontinuous {G}alerkin method for
  time-dependent convection--diffusion systems, SIAM Journal on Numerical
  Analysis 35~(6) (1998) 2440--2463.
\newblock \href {http://dx.doi.org/10.1137/S0036142997316712}
  {\path{doi:10.1137/S0036142997316712}}.

\bibitem{ReedHill1973}
H.~Reed, T.~R. Hill, Triangular mesh methods for the neutron transport
  equation, Tech. Rep. LA-UR-73-479, Los Alamos Scientific Laboratory, NM
  (1973).

\bibitem{Johnson1986}
C.~Johnson, J.~Pitk{\"a}ranta, An {A}nalysis of the {D}iscontinuous {G}alerkin
  {M}ethod for a {S}acalar {H}yperbolic {E}quation, Mathematics of Computation
  46~(173) (1986) 1--26.

\bibitem{Kuzmin2012}
D.~Kuzmin, Slope limiting for discontinuous {G}alerkin approximations with a
  possibly non-orthogonal {T}aylor basis, International Journal for Numerical
  Methods in Fluids 71~(9) (2013) 1178--1190.
\newblock \href {http://dx.doi.org/10.1002/fld.3707}
  {\path{doi:10.1002/fld.3707}}.

\bibitem{GottliebShu1998}
S.~Gottlieb, C.-W. Shu, Strong stability-preserving high-order time
  discretization methods, Math. Comp. 67~(221) (1998) 73--85.
\newblock \href {http://dx.doi.org/10.1090/S0025-5718-98-00913-2}
  {\path{doi:10.1090/S0025-5718-98-00913-2}}.

\bibitem{GottliebShu2001}
S.~Gottlieb, C.-W. Shu, E.~Tadmor, Strong stability-preserving high-order time
  discretization methods, SIAM Review 43~(1) (2001) 89--112.

\bibitem{CockburnShuRKDG21989}
B.~Cockburn, C.-W. Shu, {TVB} {R}unge-{K}utta local projection discontinuous
  {G}alerkin finite element method for conservation laws. {II}. {G}eneral
  framework, Mathematics of computation 52~(186) (1989) 411--435.

\bibitem{Krivodonova2004}
L.~Krivodonova, J.~Xin, J.-F. Remacle, N.~Chevaugeon, J.~E. Flaherty, Shock
  detection and limiting with discontinuous galerkin methods for hyperbolic
  conservation laws, Applied Numerical Mathematics 48~(3) (2004) 323--338.

\bibitem{Tu2005}
S.~Tu, S.~Aliabadi, A slope limiting procedure in discontinuous galerkin finite
  element method for gasdynamics applications, International Journal of
  Numerical Analysis and Modeling 2~(2) (2005) 163--178.

\bibitem{Michoski2011}
C.~Michoski, C.~Mirabito, C.~Dawson, D.~Wirasaet, E.~Kubatko, J.~Westerink,
  Adaptive hierarchic transformations for dynamically p-enriched slope-limiting
  over discontinuous galerkin systems of generalized equations, Journal of
  Computational Physics 230~(22) (2011) 8028 -- 8056.
\newblock \href {http://dx.doi.org/10.1016/j.jcp.2011.07.009}
  {\path{doi:10.1016/j.jcp.2011.07.009}}.

\bibitem{Yang2009}
M.~Yang, Z.-J. Wang, A parameter-free generalized moment limiter for high-order
  methods on unstructured grids, Adv. Appl. Math. Mech 1~(4) (2009) 451--480.

\bibitem{Zhang2012}
X.~Zhang, Y.~Xia, C.-W. Shu, Maximum-principle-satisfying and
  positivity-preserving high order discontinuous {G}alerkin schemes for
  conservation laws on triangular meshes, J. Sci. Comput. 50~(1) (2012) 29--62.
\newblock \href {http://dx.doi.org/10.1007/s10915-011-9472-8}
  {\path{doi:10.1007/s10915-011-9472-8}}.

\bibitem{Kuzmin2010}
D.~Kuzmin, A vertex-based hierarchical slope limiter for adaptive discontinuous
  {G}alerkin methods, Journal of Computational and Applied Mathematics 233~(12)
  (2010) 3077--3085, {F}inite Element Methods in Engineering and Science
  (FEMTEC 2009).
\newblock \href {http://dx.doi.org/10.1016/j.cam.2009.05.028}
  {\path{doi:10.1016/j.cam.2009.05.028}}.

\bibitem{YangWang2009}
M.~Yang, Z.-J. Wang, A parameter-free generalized moment limiter for high-order
  methods on unstructured grids, in: 47th AIAA Aerospace Sciences Meeting
  Including The New Horizons Forum and Aerospace Exposition, AIAA-2009-605.
\newblock \href {http://dx.doi.org/10.2514/6.2009-605}
  {\path{doi:10.2514/6.2009-605}}.

\bibitem{Michalak2008}
K.~Michalak, C.~Ollivier-Gooch, Limiters for unstructured higher-order accurate
  solutions of the euler equations, in: 46th AIAA Aerospace Sciences Meeting
  and Exhibit, AIAA-2008-776.
\newblock \href {http://dx.doi.org/10.2514/6.2008-776}
  {\path{doi:10.2514/6.2008-776}}.

\bibitem{Luo2008}
H.~Luo, J.~D. Baum, R.~L{\"o}hner, A discontinuous {G}alerkin method based on
  a~{T}aylor basis for the compressible flows on arbitrary grids, Journal of
  Computational Physics 227~(20) (2008) 8875 -- 8893.
\newblock \href {http://dx.doi.org/10.1016/j.jcp.2008.06.035}
  {\path{doi:10.1016/j.jcp.2008.06.035}}.

\bibitem{Aizinger2011}
V.~Aizinger, A geometry independent slope limiter for the discontinuous
  {G}alerkin method, in: E.~Krause, Y.~Shokin, M.~Resch, D.~Kröner, N.~Shokina
  (Eds.), Computational Science and High Performance Computing IV, Vol. 115 of
  Notes on Numerical Fluid Mechanics and Multidisciplinary Design, Springer
  Berlin Heidelberg, 2011, pp. 207--217.
\newblock \href {http://dx.doi.org/10.1007/978-3-642-17770-5_16}
  {\path{doi:10.1007/978-3-642-17770-5_16}}.

\bibitem{BarthJespersen1989}
T.~Barth, D.~Jespersen, The design and application of upwind schemes on
  unstructuredmeshes, in: Proc. AIAA 27th Aerospace Sciences Meeting, Reno,
  1989.

\bibitem{Cools2003}
R.~Cools, An encyclopaedia of cubature formulas, Journal of Complexity 19~(3)
  (2003) 445--453.
\newblock \href {http://dx.doi.org/10.1016/S0885-064X(03)00011-6}
  {\path{doi:10.1016/S0885-064X(03)00011-6}}.

\bibitem{CockburnShu1998b}
B.~Cockburn, C.-W. Shu, The {R}unge--{K}utta discontinuous {G}alerkin method
  for conservation laws~{V}: multidimensional systems, J. Comput. Phys. 141~(2)
  (1998) 199--224.
\newblock \href {http://dx.doi.org/10.1006/jcph.1998.5892}
  {\path{doi:10.1006/jcph.1998.5892}}.

\bibitem{LeVeque1996}
R.~J. Leveque, High-resolution conservative algorithms for advection in
  incompressible flow, SIAM Journal on Numerical Analysis 33~(2) (1996)
  627--665.
\newblock \href {http://dx.doi.org/10.2307/2158391}
  {\path{doi:10.2307/2158391}}.

\bibitem{kron}
L.~Sorber,
  \href{http://de.mathworks.com/matlabcentral/fileexchange/28889-kronecker-product}{{K}ronecker
  product}, {MATLAB} Central File Exchange. Retrieved October 30, 2015 (2010).
\newline\urlprefix\url{http://de.mathworks.com/matlabcentral/fileexchange/28889-kronecker-product}

\end{thebibliography}
\bibliographystyle{elsarticle-num}

\end{document}